\newif\iffinal
\finalfalse	
\finaltrue	

\documentclass[reqno,12pt]{amsart}
\usepackage{amsmath,amsthm,amsfonts,amssymb,amstext,upref}
\usepackage[foot]{amsaddr}
\usepackage[mathscr]{euscript}
\iffinal\else\usepackage[notref,notcite]{showkeys}\fi
\usepackage{hyperref}
\usepackage{graphicx}
\usepackage{color}
\usepackage[usenames,dvipsnames]{pstricks}
\usepackage{epsfig}
\usepackage{pst-grad} 
\usepackage{pst-plot} 
\usepackage[textsize=small]{todonotes}       

\usepackage{pdfsync}

\usepackage{fourier}
\usepackage[margin=1.05in]{geometry}

\usepackage{enumerate}
\newenvironment{enumeratei}{\begin{enumerate}[\upshape (i)]}{\end{enumerate}}

\newenvironment{enumerateA}{\begin{enumerate}[\upshape (A)]}{\end{enumerate}}

\numberwithin{equation}{section}
\theoremstyle{plain}

\newtheorem{theorem}{Theorem}[section]
\newtheorem*{theorem*}{Theorem}
\newtheorem{lemma}[theorem]{Lemma}

\newtheorem{proposition}[theorem]{Proposition}

\newtheorem*{proposition*}{Proposition}
\newtheorem*{thm2.4*}{Theorem 2.4*}

\theoremstyle{definition}
\newtheorem{definition}[theorem]{Definition}

\newtheorem{remark}[theorem]{Remark}

\newcommand{\E}{\mathbb{E}}
\newcommand{\PR}{\mathbb{P}}
\newcommand{\eps}{\varepsilon}
\newcommand{\poi}{\mathcal{P}}
\newcommand{\C}{\mathcal{C}}
\newcommand{\R}{\mathbb{R}}
\newcommand{\Z}{\mathbb{Z}}
\newcommand{\I}{\mathbb{I}}
\newcommand{\EL}{\mathcal{L}}

\newcommand{\cC}{\mathcal{C}}
\newcommand{\cD}{\mathcal{D}}\newcommand{\cE}{\mathcal{E}}
\newcommand{\cH}{\mathcal{H}}\newcommand{\cI}{\mathcal{I}}
\newcommand{\cJ}{\mathcal{J}}\newcommand{\cL}{\mathcal{L}}

\newcommand{\cQ}{\mathcal{Q}}\newcommand{\cR}{\mathcal{R}}

\newcommand{\cW}{\mathcal{W}}\newcommand{\cX}{\mathcal{X}}

\newcommand{\mvtwo}{\boldsymbol{2}}\newcommand{\mvthree}{\boldsymbol{3}}

\newcommand{\mveight}{\boldsymbol{8}}
\newcommand{\mvA}{\boldsymbol{A}}\newcommand{\mvB}{\boldsymbol{B}}\newcommand{\mvC}{\boldsymbol{C}}
\newcommand{\mvD}{\boldsymbol{D}}

\newcommand{\mvd}{\boldsymbol{d}}

\newcommand{\mvn}{\boldsymbol{n}}
\newcommand{\mvr}{\boldsymbol{r}}

\newcommand{\mvx}{\boldsymbol{x}}

\newcommand{\mvdelta}{\boldsymbol{\delta}}

\newcommand{\fJ}{\mathfrak{J}}

\newcommand{\bI}{\mathbb{I}}

\newcommand{\bP}{\mathbb{P}}\newcommand{\bR}{\mathbb{R}}

\newcommand{\bZ}{\mathbb{Z}}

\DeclareMathOperator{\var}{Var}

\DeclareMathOperator{\inn}{in}

\begin{document}

\title[CLT for MSTs]{Minimal Spanning Trees and Stein's Method}

\author[Chatterjee]{Sourav Chatterjee$^1$}
\address{$^1$Department of Statistics, Stanford University, 390 Serra Mall, Sequoia Hall, Stanford, CA, USA}
\author[Sen]{Sanchayan Sen$^2$}
\address{$^2$Department of Mathematics and Statistics, McGill University, 805 Sherbrooke Street West, Montr\'{e}al, QC, Canada}
\email{souravc@stanford.edu, sanchayan.sen1@gmail.com }


\begin{abstract}
Kesten and Lee \cite{kesten} proved that the total length of a minimal
spanning tree on certain random point configurations in $\R^d$ satisfies a central limit
theorem. They also raised the question: how to make these results quantitative?
Error estimates in central limit theorems satisfied by many other standard functionals studied in geometric probability
are known, but techniques employed to tackle the problem for those functionals do not apply directly to the minimal spanning tree.
Thus the problem of determining the convergence rate in the central limit theorem for Euclidean
minimal spanning trees has remained open. In this work, we establish bounds on
the convergence rate for the Poissonized version of this problem by using a variation of Stein's method. We also
derive bounds on the convergence rate for the analogous problem in the setup
of the lattice $\Z^d$.\\

The contribution of this paper is twofold. First, we develop a general technique to compute convergence rates in central limit theorems satisfied by minimal spanning trees on sequences of weighted graphs, including minimal spanning trees on Poisson points inside a sequence of growing cubes. Secondly, we present a way of quantifying the Burton-Keane argument for the uniqueness of the infinite open cluster. The latter is interesting in its own right and based on a generalization of our technique, Duminil-Copin, Ioffe and Velenik \cite{duminil-copin} have recently obtained bounds on probability of two-arm events in a broad class of translation-invariant percolation models.

\end{abstract}

\maketitle


\section{Introduction}\label{sec:intro}

Consider a finite, connected weighted graph $(V,E,w)$
where $(V,E)$ is the underlying graph and $w:E\to [0,\infty)$
is the weight function. A spanning tree of $(V,E)$ is a tree which is a
connected subgraph of $(V,E)$ with vertex set $V$. A minimal spanning tree (MST)
$T$ of $(V,E,w)$ satisfies
\[
\sum_{e\in T}w(e)=\min\bigg\{\sum_{e\in T'}w(e):\ T'\text{ is a spanning tree of }(V,E)\bigg\}.
\]
In this paper, whenever $(V,E)$ is a graph on some random point configuration in $\R^d$, the weight function will map every edge to its Euclidean length.

Minimal spanning trees and other related functionals are of great interest
in geometric probability. For an account of law of large numbers and related asymptotics for these
functionals, see e.g~\cite{aldous, alexanderI, avramI, beardwood, steeleI, steeleII}.
One of the early successes in the direction of proving distributional convergence of such functionals came with the paper of Avram and Bertsimas \cite{avram} in 1993 where the authors proved
central limit theorems (CLT) for three such functionals, namely the
lengths of the $k$-th nearest neighbor graph, the Delaunay triangulation,
and the Voronoi diagram on Poisson point configurations in $[0,1]^2$.
Central limit theorems for minimal spanning trees were first proven
by Kesten and Lee \cite{kesten} and by Alexander \cite{alexander} in 1996. This was a long-standing open question at the time of its solution.
In \cite{kesten}, the CLT for the total weight of an MST on
both the complete graph on Poisson points inside $[0,n^{1/d}]^d$
and the complete graph on $n$ i.i.d. uniformly distributed points inside
$[0,1]^d$ were established when $d\geq 2$. (Their results included the case of more general
weight functions and not just Euclidean distances.) Alexander \cite{alexander}
proved the CLT for the Poissonized problem in two dimensions.
Later certain other CLTs related to MSTs were proven in
\cite{leeI} and \cite{leeII}.

Studies related to Euclidean MSTs
in several other directions were undertaken in
\cite{leeIII, bhatt, penroseI, penroseII, penroseIII}.
An account of the structural properties of minimal spanning forests
(in both Euclidean and non-Euclidean setting) can be found in
\cite{alexander+molchanov, alexander-forest, lyons, haggstrom} and the references therein.
For an account of the scaling limit of minimal spanning trees,
see e.g.~\cite{aizenman, camia,  petegabor}.

Minimal spanning trees on the complete graph and on the hypercube have been studied extensively as well and we refer the reader to \cite{K1, K2, K3, K4, K5} for such results. In the recent preprint \cite{addario-berry}, existence of a scaling limit of the minimal spanning tree on the complete graph viewed as a metric space has been established.
Our primary focus in this paper, however, will be on minimal spanning trees on Poisson points and subsets of $\Z^d$,

The methods of \cite{alexander} and \cite{kesten} cannot be used to
get bounds on the rate of convergence to normality in the CLT
for Euclidean MSTs.
Indeed, Kesten and Lee remark that
\begin{quote}
{\it ``... [A] drawback of our approach is that it is not quantitative. Further ideas are needed to obtain an error estimate in our central limit theorem.''}
\end{quote}
A general method for tackling such a problem is to show that the function of interest satisfies certain ``stabilizing" properties \cite{penroseIV}. In \cite{penrose-yukich} (see also \cite{leeI}), it was shown that
Euclidean MSTs do satisfy a stabilizing property but there was no quantitative bound
on how fast this stabilization occurs. Quoting Penrose and Yukich \cite{penroseIV}
\begin{quote}
{\it ``Some functionals, such as those defined in terms of the minimal spanning tree, satisfy a weaker form of stabilization but are not known to satisfy exponential stabilization. In these cases univariate and multivariate central limit theorems hold ... but our [main theorem] does not apply and explicit rates of convergence are not known.''}
\end{quote}
This poses the major difficulty in obtaining an error estimate in the CLT and the problem has remained open since the work of Kesten and Lee.

In this paper we use a variation of Stein's method, given by approximation theorems from
\cite{chatterjee, lachieze}, to connect the problem of bounding the convergence rate
in this CLT to the problem of getting upper bounds on the probabilities
of certain events in the setup of continuum percolation driven by a Poisson process
and thus obtaining an error estimate in this CLT (Theorem \ref{thm:poissonmst}).
Using a similar approach, we also obtain error estimates in the CLT
for the total weights of the MSTs on subgraphs of $\Z^d$ under
various assumptions on the edge weights (Theorem~\ref{thm:latticemst}).
In Theorem \ref{thm:general graphs}, we present a general CLT satisfied by the MSTs
on subgraphs of a vertex-transitive graph.
The percolation theoretic estimates used in the proofs are given in
Section \ref{sec:burton-keane}. Our techniques for proving these percolation theoretic estimates are of independent interest.

This paper is organized as follows. In Section \ref{sec:main-results}, we state our results about convergence rates in CLTs satisfied by MSTs. In Section \ref{sec:stein's-method-survey}, we give a brief survey of literature on Stein's method and state the theorems used for Gaussian approximation. In Section \ref{sec:notations}, we introduce the necessary notation. In Section \ref{sec:burton-keane}, we state the percolation theoretic estimates we will be using. In Section \ref{sec:outline}, we briefly discuss the idea in the proof and how to connect the problem of getting convergence rates in the CLT to a problem in percolation. Section \ref{sec:mst} lists some properties and preliminary results about minimal spanning trees. Sections \ref{sec:proofs}--\ref{sec:general graphs} are devoted to proofs of the central limit theorems and the percolation theoretic estimates.

\section{Main Results}\label{sec:main-results}
We summarize our main results in this section.
Define the distance $\mathcal{D}(\mu_1,\mu_2)$ between two probability measures $\mu_1$ and $\mu_2$ on $\R$
by the sup norm of the difference between their distribution functions, or equivalently
\begin{align}\label{eqn:def-kolmogorov}
\mathcal{D}(\mu_1,\mu_2):=\sup_{x\in\R}|\mu_1(-\infty,x]-\mu_2(-\infty,x]|.
\end{align}
This metric is sometimes called the `Kolmogorov distance'. A bound on the Kolmogorov distance between two probability measures is sometimes called a `Berry-Esseen bound'.

Recall also that the Kantorovich-Wasserstein distance between
two probability measures $\mu_1$ and $\mu_2$ on $\R$ is given by
\begin{align}\label{eqn:def-wasserstein}
\mathcal{W}(\mu_1,\mu_2):=\sup\bigg\{\bigg|\int f\ d\mu_1 - \int f\ d\mu_2\bigg|: f\text{ Lipschitz with }\|f\|_{\text{Lip}}\leq 1\bigg\}.
\end{align}
Convergence in this metric implies weak convergence.

Our result on Euclidean minimal spanning trees is the following.
\begin{theorem}\label{thm:poissonmst}
Let $\poi$ be a Poisson process with intensity one in $\R^d$. Let $(V_n,E_n,w_n)$
be the complete graph on $\poi\cap [-n,n]^d$ with each edge weighted by
its Euclidean length.  Let $\mu_n$ be the law of $(M_n-\E(M_n))/\sqrt{\mathrm{Var}(M_n)}$,
where $M_n$ is the total weight of an MST of $(V_n,E_n,w_n)$. Let $\gamma$ denote the standard normal distribution on $\R$.\\
(i) When $d=2$, there exist positive constants $\xi$ and $c_1$ such that for every $n\geq 1$,
\begin{equation}\label{eqn:theorem 1 d=2}
\max\big\{\cW(\mu_n,\gamma),\ \cD(\mu_n,\gamma)\big\}\leq c_1n^{-\xi}.
\end{equation}
(ii) When $d\geq 3$, for every $p>1$ and every $n\geq 2$,
\begin{equation}\label{eqn:theorem 1 d>2}
\max\big\{\cW(\mu_n,\gamma),\ \cD(\mu_n,\gamma)\big\}\leq c_2\left(\log n\right)^{-\frac{d}{4p}}
\end{equation}
for a positive constant $c_2$ depending only on $p$ and $d$.
\end{theorem}

\begin{remark}
If $\poi_{\lambda}$ is a Poisson process with intensity $\lambda>0$ in $\R^d$ and
$M_n(\lambda)$ is the weight of a minimal spanning tree of the complete graph on
$\poi_{\lambda}\cap[-n/\lambda^{\frac{1}{d}},n/\lambda^{\frac{1}{d}}]^d$, then $(M_n(\lambda)-\E M_n(\lambda))/\sqrt{\mathrm{Var}(M_n(\lambda))}$
is distributed as $\mu_n$ where $\mu_n$ is as defined in the statement of Theorem \ref{thm:poissonmst}.
For this reason, it is enough to consider only Poisson processes with intensity one.
\end{remark}

Our next theorem deals with the case of minimal spanning trees on subsets of $\mathbb{Z}^d$.
To state the theorem conveniently, we first make a definition. In what follows,
$p_c=p_c(\Z^d)$ denotes the critical probability of bond percolation in $\Z^d$ (see, e.g., \cite{grimmett, bollobas}).

\begin{definition}\label{def:property-ABCD}
A probability measure $\mu$ on $[0,\infty)$ satisfies
\begin{enumerateA}
\item {\bf Property $\mvA_{\mvdelta}$} (for some $\delta>0$) if $\mu$ has unbounded support and $\int_{0}^{\infty} x^{4+\delta}\mu(dx)<\infty$;
\item {\bf Property $\mvB$} if $\mu$ has bounded support;
\item {\bf Property $\mvC$} if either $\mu[0,x]=p_c(\Z^d)$ for some unique $x\in\R$, or $\mu[0,x)=p_c(\Z^d)$ for some unique $x\in\R$;
\item {\bf Property $\mvD$} if $\mu[0,x]>p_c(\Z^d)>\mu[0,x)$ for some $x\in \R$.
\end{enumerateA}
\end{definition}

\begin{theorem}\label{thm:latticemst}
Let $d\geq 2$ and assume that the edges of the lattice $\Z^d$ have been given i.i.d. nonnegative weights having some non-degenerate distribution $\mu$.
Let $M_n$ denote the total weight of an MST of the weighted
subgraph of $\Z^d$ within the cube $[-n,n]^d$, and let $\nu_n$ be the distribution of
$(M_n-\E(M_n))/\sqrt{\mathrm{Var}(M_n)}$. Let $\gamma$ be the standard normal distribution on $\R$.
\begin{enumeratei}
\item If $\mu$ satisfies either Property $B$ or Property $A_{\delta}$ for some $\delta>0$, then
for every $n\geq 2$,
\begin{equation}\label{eqn:d=2}
\cW(\nu_n,\gamma)\leq
\eps_n
(\log n)^{\frac{1}{4(1+3\xi)}}\big/n^{\frac{1}{6(1+2\xi)}},
\end{equation}
where
\begin{align}
\xi=
\left\{
\begin{array}{l}
1/\delta, \text{ if }\mu\text{ satisfies Property }A_{\delta},\\
0,\text{ if }\mu\text{ satisfies Property B},
\end{array}
\right.\nonumber
\end{align}
and $\eps_n\to 0$ if $\mu$ satisfies Property $C$, and is a bounded sequence otherwise.

If $\mu$ satisfies either Property $B$ or Property $A_{\delta}$ for some $\delta\geq 2$, then \eqref{eqn:d=2} holds if we replace $\cW(\nu_n,\gamma)$ by $\cD(\nu_n,\gamma)$.

\medskip

\item If $\mu$ satisfies Property $D$ and either Property $B$ or Property $A_{\delta}$ for some $\delta>0$, then for every $\eta<d/2$,
\begin{equation}\label{eqn:L3}
\cW(\nu_n,\gamma)
\leq
c_3 n^{-\eta}\text{ for }n\geq 1,
\end{equation}
where $c_3$ is a positive constant depending on $\mu$, $d$ and $\eta$.

If $\mu$ satisfies Property D and either Property $B$ or Property $A_{\delta}$ for some $\delta\geq 2$, then \eqref{eqn:L3} holds if we replace $\cW(\nu_n,\gamma)$ by $\cD(\nu_n,\gamma)$.
\end{enumeratei}

\end{theorem}

\begin{remark}
It is very likely that the bounds are sub-optimal.
However, the question of optimal error bounds is probably very difficult. Improving the bounds stated in Theorems \ref{thm:poissonmst} and \ref{thm:latticemst} can be thought of as an independent problem in percolation (see Remark \ref{rem:two-arm}).
\end{remark}

Our approach can be used to give
a simple proof of asymptotic normality of the total weight of the minimal spanning
tree under a very general assumption on the underlying graph. We present this result
in the following theorem. {\it The advantage of this approach is that
we can get a convergence rate in the central limit theorem whenever
we can prove the percolation theoretic estimates analogous to the ones
used in the proofs of Theorem \ref{thm:poissonmst} and Theorem \ref{thm:latticemst}.}

Before stating the theorem, let us recall the definition of a vertex-transitive graph.
A graph $G=(V, E)$ is said to be vertex-transitive if for any $v_1, v_2\in V$, there exists a graph automorphism $f$ of $G$ such that $f(v_1)=v_2$.

For a graph $G=(V,E)$ and a vertex $v\in V$, we will write $S_G(v,r)$ to denote the
subgraph of $G$ spanned by the set of
all vertices $v'\in V$ such that $d_G(v',v)\leq r$ where $d_G$ denotes the graph distance of $G$.
\begin{theorem}\label{thm:general graphs}
Let $G=(V,E)$ be a\\
(I) connected, infinite, locally finite, vertex-transitive graph.

Consider a sequence of finite connected subgraphs $G_n=(V_n,E_n)$ such that\\
(II) $|V_n|\to\infty$, and \\
(III) $|\big\{v\in V_n:\ S_G(v,r)\not\subset G_n\big\}|=o(|V_n|)$ for every $r>0$.

Consider i.i.d. nonnegative weights associated with the edges of $G$ where the weights follow
some non-degenerate distribution $\mu$ that satisfies either Property $B$ or Property $A_{\delta}$ for some $\delta>0$.
Let $M_n$ be the total weight of a minimal spanning tree of $G_n$. Then\\
(i) $\mathrm{Var}(M_n)=\Theta(|V_n|)$ and\\
(ii) $(M_n-\E (M_n))/\sqrt{\mathrm{Var}(M_n)}\stackrel{d}{\to} Z$, where $Z$ follows a $N(0,1)$ distribution.
\end{theorem}

\begin{remark}\label{rem:amenable}
Note that $G$ in Theorem \ref{thm:general graphs} is necessarily amenable (because of Conditions (II) and (III)).
\end{remark}

\section{Stein's method}\label{sec:stein's-method-survey}
In 1972, Charles Stein \cite{stein72} proposed a radically different approach to proving convergence to normality.
Stein's observation was that the standard normal distribution is the only probability distribution that satisfies the equation
\begin{equation*}
\E (Zf(Z)) = \E f'(Z)
\end{equation*}
for all absolutely continuous $f$ with a.e.~derivative $f'$ such that $\E|f'(Z)|<\infty$. From this, one might expect that if $W$ is a random variable that satisfies the above equation in an approximate sense, then the distribution of $W$ should be close to the standard normal distribution.
The key to Stein's implementation of his idea was the method of exchangeable pairs, devised by Stein in \cite{stein72}. A notable success story of Stein's method was authored by Bolthausen~\cite{bolthausen} in 1984, when he used a sophisticated version of the method of exchangeable pairs to obtain an error bound in a famous combinatorial central limit theorem of Hoeffding. Stein's 1986 monograph \cite{stein86} was the first book-length treatment of Stein's method. After the publication of \cite{stein86}, the field was given a boost by the popularization of the method of dependency graphs by Baldi and Rinott \cite{baldi}, a striking application to the number of local maxima of random functions by Baldi, Rinott and Stein \cite{baldi2}, and central limit theorems for random graphs by Barbour, Karo\'nski and Ruci\'nski \cite{barbour2}, all in 1989.

The new surge of activity that began in the late eighties continued through the nineties, with important contributions coming from Barbour \cite{barbour0} in 1990, who introduced the diffusion approach to Stein's method; Avram and Bertsimas \cite{avram} in 1993, who applied Stein's method to solve an array of important problems in geometric probability; Goldstein and Rinott \cite{goldstein2} in 1996, who developed the method of size-biased couplings for Stein's method, improving on earlier insights of Baldi, Rinott and Stein \cite{baldi2}; Goldstein and Reinert \cite{goldstein1} in 1997, who introduced the method of zero-bias couplings; and Rinott and Rotar \cite{rinott} in 1997, who solved a well known open problem related to the antivoter model using Stein's method. Sometime later, in 2004, Chen and Shao \cite{chen2} did an in-depth study of the dependency graph approach, producing optimal Berry-Ess\'een type error bounds in a wide range of problems. The 2003 monograph of Penrose \cite{penrosebook} gave extensive applications of the dependency graph approach to problems in geometric probability.

A new version of Stein's method with potentially wider applicability was introduced for discrete systems \cite{chatterjee}, and a corresponding continuous version in \cite{chatterjee4}. This new approach was used to solve a number of questions in geometric probability in \cite{chatterjee}, random matrix central limit theorems in \cite{chatterjee4}, and number theoretic central limit theorems in \cite{chatterjeeI}. The main result of \cite{chatterjee} gives convergence rates in terms of the Kantorovich-Wasserstein distance. Very recently, this approach has been generalized in \cite[Theorem 4.2]{lachieze} to give convergence rates in the Kolmoogorov distance. These two results are our main tools for normal approximation.

As mentioned before in Section \ref{sec:intro}, MSTs on Poisson points exhibit a stabilization property; but no tail bound on the radius of stabilization (in the sense of \cite{penrose-yukich}) is known. If such a tail bound were known, then there would be a number of ways of obtaining a convergence rate in the CLT satisfied by MSTs on Poisson points (for example, using the results of \cite{chatterjee} or \cite{last} or \cite{penroseIV}). However, \cite[Theorem 2.2]{chatterjee} and \cite[Theorem 4.2]{lachieze} allow us to circumvent this problem and instead reduce the problem to finding upper bounds on probability of two-arm events. We will state these theorems in the following section.

\subsection{Main approximation theorems}\label{sec:stein's-method}
To state the theorems, we need some notation; we will use them repeatedly in this paper.

Let $\mathcal{X}$ be a Polish space. For every $A\subset [n]:=\{1,\hdots,n\}$, define the ``replacement" operator $\cR^A:\cX^n\times\cX^n\to\cX^n$ as follows: for $y=(y_1,\hdots,y_n)\in\cX^n$, and $y'=(y_1',\hdots,y_n')\in\cX^n$, the $i$-th component of $\cR^A(y, y')$ is given by
\begin{align}
\big(\cR^A(y, y')\big)_i=
\left\{
\begin{array}{l}
y_i', \text{ if }i\in A,\\
y_i,\text{ if }i\notin A.
\end{array}
\right.\nonumber
\end{align}
Suppose $f:\cX^n\to\bR$ is a measurable function. For $j\in[n]$, define $\Delta_j f:\cX^n\times\cX^n\to\bR$ by
\[\Delta_j f(y, y'):=f(y)-f\big(\cR^{\{j\}}(y, y')\big).\]

Let $X_1,\hdots,X_n$ be independent $\cX$ valued random variables and set $X=(X_1,\hdots,X_n)$. Let $X'=(X_1',\hdots,X_n')$
be an independent copy of $X$.
To simplify notation, we will write $X^A$ to denote the random vector $\cR^A(X, X')$. We will simply write $X^j$ instead of $X^{\{j\}}$. With this convention, for every $A\subset[n]$,
\[\Delta_j f(X^A, X')=f(X^A)-f\big(X^{A\cup\{j\}}\big).\]
For every $A\subset [n]$, let
\[T_A:=\sum_{j\notin A}\Delta_j f(X, X')\Delta_j f(X^A, X'),\
\text{ and }\
T_A':=\sum_{j\notin A}\Delta_j f(X, X')|\Delta_j f(X^A, X')|.\]
Finally define
\[ T=\frac{1}{2}\sum_{A\subsetneq [n]}\frac{T_A}{\dbinom{n}{|A|}(n-|A|)},\
\text{ and }\
T'=\frac{1}{2}\sum_{A\subsetneq [n]}\frac{T_A'}{\dbinom{n}{|A|}(n-|A|)}.\]
Recall the definitions of the Kantorovich-Wasserstein distance (see \eqref{eqn:def-wasserstein}) and the Kolmogorov distance (see \eqref{eqn:def-kolmogorov}).
\begin{theorem}\label{thm:chatterjee}(\cite[Theorem 2.2]{chatterjee})
Let all terms be defined as above and let $W=f(X)$ with $\sigma^2:=\mathrm{Var}(W)<\infty$. Then $\E T=\sigma^2$ and
\begin{equation}\label{eqn:chatterjee}
\mathcal{W}(\mu,\gamma)\leq\frac{1}{\sigma^2}\big[\mathrm{Var}(\E(T|W))\big]^{1/2}+\frac{1}{2\sigma^3}\sum_{j=1}^n\E|\Delta_j f(X, X')|^3,
\end{equation}
where $\mu$ is the law of $(W-\E W)/\sigma.$
\end{theorem}
\begin{theorem}\label{thm:lrp}(\cite[Theorem 4.2]{lachieze})
Let all terms be defined as above and let $W=f(X)$ with $\sigma^2:=\mathrm{Var}(W)<\infty$. Then
\begin{align}\label{eqn:lrp}
\cD(\mu,\gamma)
&\leq\frac{1}{\sigma^2}\big[\mathrm{Var}(\E(T|X))\big]^{1/2}
+\frac{1}{\sigma^2}\big[\mathrm{Var}(\E(T'|X))\big]^{1/2}\\
&\hskip20pt+\frac{1}{4\sigma^3}\sum_{j=1}^n\big(\E|\Delta_j f(X, X')|^6\big)^{1/2}
+\frac{\sqrt{2\pi}}{16\sigma^3}\sum_{j=1}^n\E|\Delta_j f(X, X')|^3,\notag
\end{align}
where $\mu$ is the law of $(W-\E W)/\sigma.$
\end{theorem}

Note that
\[\var\big(\E(T|W)\big)\leq \mathrm{Var}(T),\ \text{ and }\ \var\big(\E(T|X)\big)\leq \mathrm{Var}(T),\]
and
\begin{align}\label{eqn:var expansion term bounded}
\mathrm{Var}(T)&=\frac{1}{4}\mathrm{Var}\Bigg[\sum_{A\subsetneq [n]}\sum_{j\in [n]\setminus A}\frac{\Delta_j f(X)\Delta_j f(X^A)}{\dbinom{n}{|A|}(n-|A|)}\Bigg]\nonumber\\
&=\frac{1}{4}
\sum_{\substack{A\subsetneq [n]\\j\in [n]\setminus A}}\sum_{\substack{A'\subsetneq [n]\\j'\in [n]\setminus A'}}
\frac{\mathrm{Cov}\bigg(\Delta_j f(X)\Delta_j f(X^A),\ \Delta_{j'} f(X)\Delta_{j'} f(X^{A'})\bigg)}
{\dbinom{n}{|A|}(n-|A|)\dbinom{n}{|A'|}(n-|A'|)}.
\end{align}
We will make repeated use of this identity.

The expression of the upper bound in Theorem \ref{thm:lrp} is very similar to the bound in Theorem \ref{thm:chatterjee}.
We will give detailed proofs of bounds in the Kantorovich-Wasserstein distance using Theorem \ref{thm:chatterjee}, and then briefly sketch how to adapt the proof using Theorem \ref{thm:lrp} to get a bound of the same order in the Kolmogorov distance.


\section{Notation}\label{sec:notations} We will use some notation
frequently throughout this paper. For convenience, we collect them together in this section.

\subsection{Euclidean setup}\label{sec:notation-continuum}
If $x$ is a point in $\R^d$ and $A\subset\R^d$, then we define
$x+A:=\{x+y : y\in A\}$. If $r>0$, $S_{\R^d}(x,r)$ will denote the closed $L^2$ ball
of radius $r$ centered at $x$, and $B_{\R^d}(x,r)$ will denote the closed $L^{\infty}$ ball of radius $r$ centered at $x$,
i.e., $B_{\R^d}(x,r)=x+[-r,r]^d$. When $x$ is the origin, we will simply write
$B_{\R^d}(r)$ instead of $B_{\R^d}(0,r)$. For any cube $B$, we refer to its center as $c(B)$.
We will denote by $d_{\R^d}(\cdot,\cdot)$, the metric induced by the $L^2$ norm in $\R^d$.
When the underlying space is clear from the context, we will drop the subscript ${\R^d}$ and simply write
$S(\cdot,\cdot)$, $B(\cdot,\cdot)$, and $d(\cdot,\cdot)$.

For a finite subset $X$ of $\R^d$, $M_{\bR^d}(X)$ will denote the sum of edge weights
of the minimal spanning tree on the complete graph on $X$ having Euclidean distance as edge weights.
When the ambient space is clear, we will drop the subscript and simply write $M(X)$.

For $A\subset\R^d$ and $r>0$, we define
\[A^{(r)}:=\big\{x\in\R^d: d_{\R^d}(x,A)\leq r\big\}.\]
(With this notation $S_{\R^d}(x, r)=\{x\}^{(r)}$.) Let us also define
\[A_{(r)}:=\big\{x\in A: d_{\R^d}(x,\partial A)\leq r\big\}.\]
Let $\poi$ be a Poisson process in $\R^d$ and let $A$ be a subset of $\R^d$.
Then $\C\subset\poi\cap A$ will be called an $r$-cluster
in $A$ (or just $r$-cluster if $A$ is clear) if $\C^{(r)}$ is a connected component of
$(\poi\cap A)^{(r)}$; $\C^{(r)}$ should be thought of as the region occupied by the cluster $\C$.
We say that two $r$-clusters $\C_1$ and $\C_2$  in $A\subset \R^d$ are disjoint if
$\C_1^{(r)}$ and $\C_2^{(r)}$ are.
We emphasize that the occupied regions must be disjoint in $\R^d$, and
it is not enough to have their restrictions to $A$ to be disjoint. We will write {\it{configuration}}
to mean a locally finite subset of $\R^d$. For $A\subset\R^d$, $\mathfrak{X}(A)$ will denote the space
of all locally finite subsets of $A$.

For two compact sets $K_1,\ K_2\subset\R^d$ with $K_1\subset K_2$, a positive integer $k$ and a positive real $r$,
we write $K_1\underset{r}{\overset{k}{\longleftrightarrow}}K_2$ if
there exists a collection of $k$ disjoint $r$-clusters $\C_1,\hdots,\C_k$ in $K_2\setminus K_1$ such that
\[\C_j\cap K_1^{(r)}\neq \emptyset\text{ and }\C_j\cap (K_2)_{(2r)}\neq \emptyset,\text{ for } j=1,\hdots,k.\]
For $x\in\R^d$ and $b>a>0$, we call $\{B(x, a)\underset{r}{\overset{2}{\longleftrightarrow}}B(x, b)\}$ a two-arm event at level $r$.

We will write $K_1\underset{r}{\overset{k}{\longrightarrow}}K_2$, if
there exists a collection of $k$ pairwise disjoint $r$-clusters $\C_1,\hdots,\C_k$ in $(K_2\setminus K_1)$ such that
\begin{align}
\C_j\cap K_1^{(2r)}\neq \emptyset\text{ and }\C_j\cap (K_2)_{(2r)}\neq \emptyset,\text{ for } j=1,\hdots,k.
\end{align}


\subsection{Discrete setup}\label{sec:notation-discrete}
Consider a graph $G=(V,E)$. Recall from Section \ref{sec:main-results} that $d_G(\cdot,\cdot)$ denotes the graph distance on $G$, and
\[S_G(v,r):=\big\{v'\in V\ :\  d_G(v',v)\leq r\big\}.\]
Assume that each $e\in E$ has a nonnegative weight $x_e$ attached to it. Let
$\mvx=(x_e:e\in E)$. Then for any finite connected subgraph $H=(V_1, E_1)$ of $G$,
$M_G(H, \mvx)$ will denote the total weight of an MST on the weighted graph $H$, where  $e_1\in E_1$ has weight $x_{e_1}$.
When the underlying graph $G$ is clear, we will drop the subscripts and simply write $d(\cdot,\cdot)$, $S(v,r)$, and $M(H, \mvx)$.

For any $e\in E$, $G-e$ will denote the graph $(V,E-e)$. If $G_i=(V_i, E_i)$, $i=1, 2$ are two subgraphs of $G$, then $G_1\cap G_2$ will denote the subgraph $(V_1\cap V_2, E_1\cap E_2)$.

When working with
the lattice $\Z^d$, $B_{\Z^d}(x,r)$ will denote the set of all lattice points
inside $x+[-r,r]^d$ and $B_{\Z^d}(r)$ will stand for $B_{\Z^d}(0,r)$.
We will simply write $B(x,r)$ and $B(r)$ when the ambient space is clear from the context.

For a subset $V$ of $\Z^d$, let $G(V)$ denote the
subgraph of $\Z^d$ induced by $V$.
We will sometimes make abuse of notation
by referring to $G(V)$ as $V$. With this convention
$B_{\Z^d}(x,r)$ will sometimes mean $G(B_{\Z^d}(x,r))$ and the meaning
will be clear from the context.
For a cube $Q$ in $\Z^d$, $\partial^{\inn} Q$ will denote the ``inner vertex boundary" of $Q$, i.e., the set of all vertices
in $Q$ that are adjacent to at least one vertex not in $Q$.

For $p\in[0,1]$, consider i.i.d. Bernoulli$(p)$ random variables $\{X_e\}_{e\in\Z^d}$ associated with
edges of $\Z^d$, i.e. $\PR(X_e=1)=p=1-\PR(X_e=0)$. We call an edge $e$
open (resp. closed) at level $p$ if $X_e=1$ (resp. $X_e=0$).
Given a subgraph $G=(V,E)$ of $\Z^d$ and $V'\subset V$,
we say that $V'$ forms a $p$-cluster in $G$ if there is a path consisting of
open edges in $E$ between any two vertices in $V'$ and $V'$ is a maximal subset
of $V$ in this regard.

For two cubes $Q_1\subset Q_2$ in $\Z^d$, denote by $Q_2-Q_1$ the
subgraph $(V,E)$ of $Q_2$ with
\[E=\{\text{all edges in }Q_2\text{ except
the ones with both endpoints in }Q_1\}\text{ and}\]
\[V=\{v:\ v\text{ is an endpoint of }e\text{ for some }e\in E\}.\]
For two cubes $Q_1\subset Q_2$ in $\Z^d$ and $p\in[0,1]$, $Q_1\underset{p}{\overset{k}{\leftrightsquigarrow}}Q_2$
will mean that there exist at least $k$ disjoint $p$-clusters in $Q_2-Q_1$ that
intersect both $\partial^{\inn} Q_1$ and $\partial^{\inn} Q_2$.
If $Q_1$, $Q_2$, $Q_3$ are cubes in $\Z^d$ such that
(i) $Q_1\subset Q_2\cap Q_3$, and (ii) $\partial^{\inn} Q_2$ has a vertex in $Q_3$, then we will write
``$Q_1\underset{p}{\overset{k}{\leftrightsquigarrow}}Q_2\text{ in }Q_3$" if
there exist $k$ disjoint $p$-clusters in $(Q_2-Q_1)\cap Q_3$
each intersecting $\partial^{\inn} Q_1$ and $\partial^{\inn} Q_2$.

\begin{figure}[htbp]
\centering
\begin{minipage}{.5\textwidth}
  \centering
  \includegraphics[trim=4cm 16.5cm 8.2cm 2.6cm, clip=true, angle=0, scale=0.5]{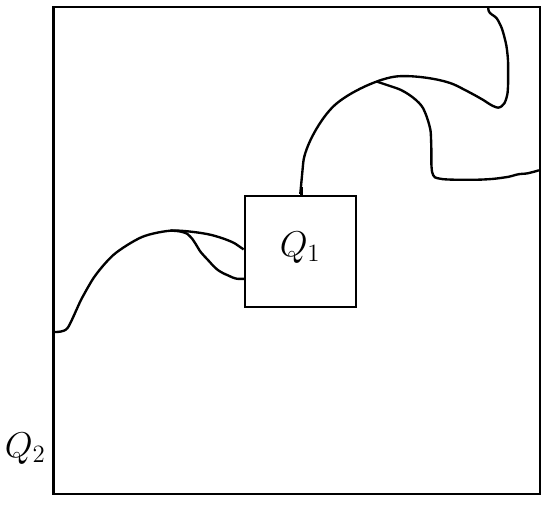}
  \caption{$Q_1\underset{p}{\overset{2}{\leftrightsquigarrow}}Q_2$}
  \label{fig:twoarm}
\end{minipage}%
\begin{minipage}{.5\textwidth}
  \centering
  \includegraphics[trim=4cm 15.5cm 3.5cm 2.6cm, clip=true, angle=0, scale=0.5]{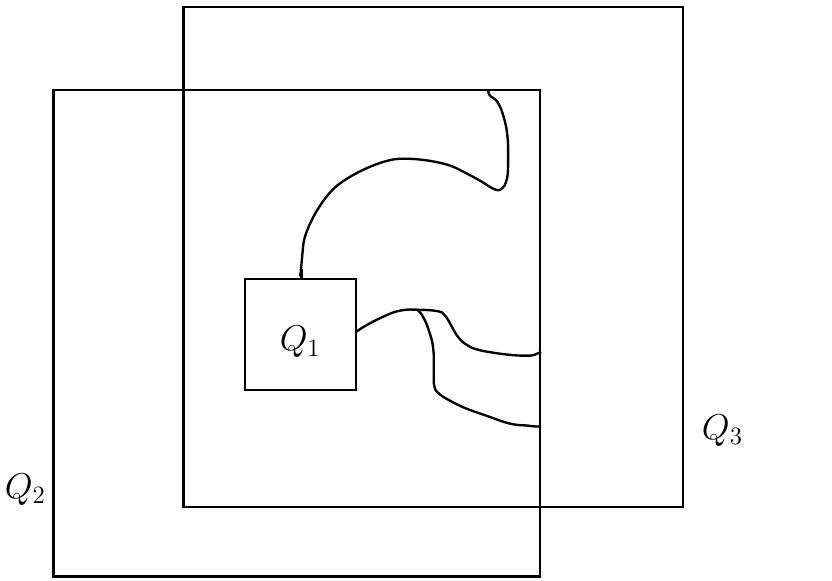}
  \caption{$Q_1\underset{p}{\overset{2}{\leftrightsquigarrow}}Q_2\text{ in }Q_3$}
  \label{fig:twoarm1}
\end{minipage}
\end{figure}


For an edge $\{x,y\}$ in $\Z^d$ and a cube $Q$ containing both $x$ and $y$,
$\{x,y\}\underset{p}{\overset{2}{\leftrightsquigarrow}}Q$ will mean that the $p$-clusters in $Q$ containing
$x$ and $y$ are disjoint and that they both intersect $\partial^{\inn} Q$. Similarly
we can define $\{x,y\}\underset{p}{\overset{2}{\leftrightsquigarrow}}Q-\{x,y\}$  to be the event that
the $p$-clusters in $Q-\{x,y\}$ containing
$x$ and $y$ intersect $\partial^{\inn} Q$ and are disjoint.

Assume that $\{x,y\}$ is an edge of $\Z^d$ and $n\geq 2$. Analogous to the continuum setup, we call $\{\{x,y\}\underset{p}{\overset{2}{\leftrightsquigarrow}}B_{\Z^d}(x, n)\}$ or $\{B_{\Z^d}(x, 1)\underset{p}{\overset{2}{\leftrightsquigarrow}}B_{\Z^d}(x, n)\}$ a two-arm event at level $p$.

\subsection{Convention about constants}
To ease notation, most constants in this paper will be denoted by
$c$, $c'$, $C$ etc. and their values may change from line to line. These
constants may depend on parameters like the dimension
and often we will not mention this dependence explicitly; none of these constants will depend
on the quantity ``$n$,'' used to index infinite sequences.
Specific constants will have a subscript as for example $c_1$, $c_2$ etc.

\section{Two-arm event: Quantification of the Burton-Keane argument}\label{sec:burton-keane}
The key ingredients in the proofs of Theorem \ref{thm:poissonmst} and Theorem
\ref{thm:latticemst} are some percolation theoretic estimates
which are of independent interest. We state them in the following lemmas.
\begin{lemma}\label{lem:critical connectivity}
Assume $d\geq 3$ and let $\poi$ be a Poisson process having intensity one in $\R^d$. Let $0<r_1<r_2<\infty$.
Then there exist constants $c_6$ and $c_7$ depending only on $r_1,\ r_2$ and $d$ such that for
every $r\in[r_1, r_2]$, every $n\geq 2$ and every $a\in(1/2, (\log\log n)^{1/(d-1/2)})$, we have
\begin{equation}\label{eqn:critical connectivity}
\PR\big(B_{\bR^d}(a)\underset{r}{\overset{2}{\longleftrightarrow}}B_{\bR^d}(n)\big)\leq \frac{c_6\exp(c_7 a^{d-1})}{(\log n)^\frac{d}{2}}.
\end{equation}
The same bound holds if we replace $B_{\bR^d}(a)$ by $B_{\bR^d}(a)^{(r)}$ or $B_{\bR^d}(a)^{(r)}\cup S_{\bR^d}(x,r)$ for some $x\in B_{\bR^d}(a)^{(r)}$.
\end{lemma}
The proof of this lemma is given in Section \ref{sec:proofs}.
Lemma \ref{lem:critical connectivity} deals with the case $d\geq 3$. The case $d=2$ is simpler and will be handled in Lemma \ref{lem:->connectivity bound}. The next lemma states a similar result for the lattice case.

\begin{lemma}(\cite[Proposition 5.3]{cerf})\label{lem:cerf}
Consider the lattice $\Z^d$ where $d\geq 2$ and let $e_1,\hdots,e_{2d}$ be as in Lemma \ref{lem:critical connectivityl lattice}. Then for any $0<p_1<p_2<1$, there exists a constant $c_9$ depending only on $p_1,\ p_2$ and $d$ such that for any $p\in[p_1, p_2]$ and $n\geq 2$,
\begin{equation}\label{eqn:cerf}
\PR\big(\{0,e_i\}\underset{p}{\overset{2}{\leftrightsquigarrow}}B_{\bZ^d}(n)\big)\leq c_9\bigg(\frac{\log n}{n}\bigg)^{1/2}, \text{ for }1\leq i\leq 2d.
\end{equation}
The same bound holds if we replace the edge $\{0,e_i\}$ by the cube $B_{\bZ^d}(1)$.
\end{lemma}

\begin{remark}\label{rem:only matters at criticality}
Let $r_c=r_c(d)$ be the critical radius for continuum percolation in $\R^d$ driven by
a Poisson process with intensity one (see, e.g., \cite{alexander} or \cite[Chapter 8]{bollobas}).
Note that we can actually get an exponentially decaying bound in (\ref{eqn:critical connectivity})
when $r_2< r_c$. It is also possible to prove exponential decay in (\ref{eqn:critical connectivity}) if $r_1>r_c$. So
the bound in (\ref{eqn:critical connectivity}) is really useful when $r_c\in(r_1,r_2)$.

The same is true for Lemma \ref{lem:cerf}. Exponential decay in \eqref{eqn:cerf}
is standard when $p_c(\Z^d)\notin [p_1,p_2]$.
\end{remark}

\begin{remark}
Proposition 5.3 of \cite{cerf} was actually proved for site percolation on $\Z^d$. However, the proof can be easily generalized to bond percolation. Also, the bound given in Proposition 5.3 of \cite{cerf} is of the form $O(\log n/\sqrt{n})$, but it is straightforward to modify the proof to get a bound of the form $O(\sqrt{(\log n/n)})$. Indeed, in Section 5 of \cite{cerf}, we can modify the definition of the event $\cE$ as follows:
\[\cE:=\big\{\forall C\in\cC,\
\big|h\big(\overline{C}\cap\Lambda(n)\big)\big|<\alpha(\log n)^{1/2}\big|\overline{C}\cap\Lambda(n)\big|^{1/2}\big\},\]
where $\alpha>0$ is a large constant. Then it will follow that
\[\bP(\cE^c)\leq 2|\Lambda(n)|^2\exp\big(-2\alpha^2(\log n) p^2(1-p)^2\big).\]
We can choose $\alpha$ sufficiently large and follow the rest of analysis in \cite{cerf} to get a bound of the form $O(\sqrt{(\log n/n)})$.
\end{remark}

\begin{remark}\label{rem:two-arm}
In the proof of Theorem \ref{thm:latticemst}, we need a bound on the probability of two-arm events \textit{which is uniform in $p$ over an open interval containing $p_c(\Z^d)$}. Lemma \ref{lem:cerf} serves this purpose. It is, however, possible that the estimate in Lemma \ref{lem:cerf} is sub-optimal. In \cite{cerf}, Cerf improves the bound given in Lemma \ref{lem:cerf} but only at $p=p_c$. In the recent preprint \cite{duminil-copin}, the authors prove a bound of the form $O(1/n)$ for bond percolation in $\Z^2$ (in fact, their result is true for the more general random cluster model), and Kozma and Nachmias \cite{kozma} prove a bound of the form $O(1/n^4)$ for bond percolation in $\Z^d$ when $d\geq 19$ but again, these bounds hold only at $p=p_c$. For site percolation on the triangular lattice, a bound of the form $O(n^{-5/4+o(1)})$ is known to hold at criticality \cite{smirnov}, but an analogous result is not known for the square lattice $\Z^2$.

 To the best of our knowledge the bound in \eqref{eqn:cerf} is the best known estimate valid uniformly over an interval around $p_c$. Any improvement over Lemma \ref{lem:cerf} can be used in the proof of Theorem \ref{thm:latticemst} to get better bounds in \eqref{eqn:d=2}. Similarly, any improvement over Lemma \ref{lem:critical connectivity} will yield a sharper upper bound in \eqref{eqn:theorem 1 d>2}.
\end{remark}

\begin{remark}\label{rem:quantification-of-burton-keane}
The arguments used in the proof of Lemma \ref{lem:critical connectivity} can be used in the lattice setup to get the following result.
\begin{lemma}\label{lem:critical connectivityl lattice}
Consider the lattice $\Z^d$ where $d\geq 3$. Denote the vertices
adjacent to the origin by $e_1,\hdots,e_{2d}$. Then for any $0<p_1<p_2<1$,
there exists a constant $c_8$ depending only on $p_1,\ p_2$ and $d$ such that for any $p\in[p_1, p_2]$ and $n\geq 2$,
\begin{equation}\label{eqn:critical connectivity lattice}
\PR\big(\{0,e_i\}\underset{p}{\overset{2}{\leftrightsquigarrow}}B_{\bZ^d}(n)\big)\leq c_8 (\log n)^{-\frac{d}{2}}, \text{ for }1\leq i\leq 2d.
\end{equation}
The same bound
holds if we replace the edge $\{0,e_i\}$ by the cube $B_{\bZ^d}(1)$.
\end{lemma}

The proof of this lemma is outlined briefly in Appendix \ref{Appendix}. Lemma \ref{lem:critical connectivity} and Lemma \ref{lem:critical connectivityl lattice} may be seen as quantifications of the statement that the infinite open cluster is unique. This uniqueness theorem was first proved by Aizenman, Kesten and Newman \cite{akn} for percolation on lattices. A very elegant proof was given by Burton and Keane~\cite{burtonkeane}, which has now become the standard textbook proof of the theorem. Unlike the original argument of Aizenman, Kesten and Newman, the Burton--Keane argument admits a wide array of applications and generalizations due to its simplicity and robustness.

The AKN argument is known to have a quantitative version {\it in the lattice setup} (Lemma \ref{lem:cerf}), while the Burton--Keane argument, due to its use of translation-invariance, is not expected to be quantifiable. The argument used in the proofs of Lemma \ref{lem:critical connectivity} and Lemma \ref{lem:critical connectivityl lattice} show that it is actually possible to quantify the Burton--Keane argument. Thus the technique used in the proofs of Lemma \ref{lem:critical connectivity} and Lemma \ref{lem:critical connectivityl lattice} is expected to have wider applicability in other contexts, where the Burton--Keane argument works but the AKN argument does not. As mentioned earlier, using a generalization of the arguments used in the proof of Lemma \ref{lem:critical connectivityl lattice}, Duminil-Copin, Ioffe and Velenik \cite{duminil-copin} have recently obtained bounds on the probability of two-arm events in a broad class of translation-invariant percolation models on $\Z^d$. Due to this recent development, we have included a brief sketch of the proof of Lemma \ref{lem:critical connectivityl lattice} in Appendix \ref{Appendix} even though in the proof of Theorem \ref{thm:latticemst} we will use Lemma \ref{lem:cerf} which gives a sharper bound.
\end{remark}

\section{Two standard facts about minimal spanning trees}\label{sec:two-standard-facts}
We collect two well-known facts about minimal spanning trees in this section.

\subsection{Minimax property of paths in MST}
\begin{lemma}\label{lem:mst minimax criterion}
Consider a  finite, connected and weighted graph $G=(V,E,w)$. Let $T$ be a minimal spanning tree of $G$. Then
any path $(x_0,\hdots,x_n)$ with $x_i\in V$ and $\{x_i,x_{i+1}\}\in T$
satisfies
\[\max_i\ w(\{x_i,x_{i+1}\})\leq \max_j\ w(\{x_j',x_{j+1}'\})\]
for any path $(x_0',\hdots,x_m')$ with $\{x_j',x_{j+1}'\}\in E$
and $x_0=x_0'$ and $x_n=x_m'$.
\end{lemma}
\noindent\textbf{Proof:}
This is just a restatement of \cite[Lemma 2]{kesten}.\hfill$\blacksquare$

In words, Lemma \ref{lem:mst minimax criterion} states that any path in the MST is minimax, i.e., for any two vertices $x$ and $y$, the path in the MST that connects $x$ and $y$ minimizes the maximum edge-weight among all paths in the graph that connect $x$ and $y$.

\subsection{Add and delete algorithm}\label{sec:add-and-delete-algo}
We now state an algorithm from \cite{kesten} for constructing an MST on a connected graph starting from an MST on a connected subgraph.

\begin{enumeratei}
\item {\bf Addition of an edge:}
Suppose $G_1=(V,E_1,w)$ is a finite connected weighted graph and $G_0=(V,E_0,w)$ is a connected subgraph of $G_1$ such that $E_1=E_0\cup\{e_0\}$, i.e., $G_1$ has the same vertex set and one extra edge $e_0$. Suppose $T_0$ is an MST on $G_0$. Consider the graph $T_0\cup\{e_0\}$, i.e., add the edge $e_0$ to $T_0$. Then $T_0\cup\{e_0\}$ has a unique cycle $C$. Let $e$ be an edge in $C$ such that $w(e)=\max_{e'\in C}w(e')$, and set $T_1=T_0\cup\{e_0\}\setminus e$. (Thus, we are removing an edge in $C$ that has the maximal edge-weight in $C$.)

\item {\bf Addition of a vertex:}
Suppose $G_1=(V_1,E_1,w)$ is a finite connected weighted graph and $G_0=(V_0,E_0,w)$ is a connected subgraph of $G_1$ such that $V_1=V_0\cup\{v_0\}$ and $E_1=E_0\cup\{e_0\}$. (Thus $G_1$ has one extra vertex $v_0$ and one extra edge $e_0$. Since $G_1$ is connected, $v_0$ is necessarily an endpoint of $e_0$.) Suppose $T_0$ is an MST on $G_0$. Set $T_1=T_0\cup\{e_0\}$.

\end{enumeratei}
\begin{proposition}\label{prop:add-and-delete-algo}(\cite[Proposition 2]{kesten})
The tree $T_1$ constructed in (i) or (ii) is an MST on $G_1$.
\end{proposition}

We can start from an MST on a connected graph and use the add and delete algorithm inductively to construct an MST on any larger finite connected graph.

\section{Outline of proof}\label{sec:outline}
We briefly sketch here the main ideas in the proof. For simplicity, let us consider the case where {\it the edges of $\Z^d$ have been weighted by i.i.d. Uniform$[0, 1]$ random variables}. Let $X_f$ denote the weight associated with an edge $f$ of $\bZ^d$, and let $X=(X_f:f\text{ is an edge of }B_{\Z^d}(n))$.
Heuristically, we expect $M(B_{\Z^d}(n),X)$ to satisfy a CLT if the change in $M(B_{\Z^d}(n),X)$ due to the
replacement of $X_f$ by an independent identically distributed observation
$X_f'$ ``is not observed far away from $f$.'' A quantitative formulation
of this vague statement will give us a convergence rate in the CLT.

To this end, fix $\alpha\in(0,1)$ and take an edge $e=\{x_1, x_2\}$ in $B_{\Z^d}(n)$
such that $d(x_1,\partial^{\inn} B_{\Z^d}(n))\geq \lceil n^{\alpha}\rceil$.
Let $X'$ be an independent copy of $X$. Recall the notation $X^e$ from Section \ref{sec:stein's-method}.
Define
\[\Delta_e M=M\big(B_{\Z^d}(n), X\big)-M\big(B_{\Z^d}(n), X^e\big)\ \text{ and }\
\tilde{\Delta}_e M=M\big(B_{\Z^d}(x_1,n^{\alpha}), X\big)-M\big(B_{\Z^d}(x_1,n^{\alpha}), X^e\big).\]
Then an application of Theorem \ref{thm:chatterjee}
reduces the problem to getting an upper bound on $\E|\Delta_e M-\tilde{\Delta}_e M|$.
The actual calculations are given in Section \ref{sec:actual proof of latticemst}.
This is the precise formulation of the heuristics explained above.

Noting that
\[\Delta_e M=\big[M\big(B_{\Z^d}(n),X\big)-M\big(B_{\Z^d}(n)-e, X\big)\big]
-\big[M\big(B_{\Z^d}(n),X^e\big)-M\big(B_{\Z^d}(n)-e, X^e\big)\big],\]
and a similar identity holds for $\tilde{\Delta}_e M$, it is easily seen that getting a bound on $\E|\Delta_e M-\tilde{\Delta}_e M|$
amounts to proving an upper bound on $\E|\delta_e M|$, where
\[\delta_e M:=\big[M\big(B_{\Z^d}(n),X\big)-M\big(B_{\Z^d}(n)-e, X\big)\big]
-\big[M\big(B_{\Z^d}(x_1,n^{\alpha}), X\big)-M\big(B_{\Z^d}(x_1,n^{\alpha})-e, X\big)\big].\]
It follows from Proposition \ref{prop:add-and-delete-algo} that
\[M\big(B_{\Z^d}(n), X\big)-M\big(B_{\Z^d}(n)-e, X\big)=X_e-\max\{X_e, Y\},\text{ and}\]
\[M\big(B_{\Z^d}(x_1,n^{\alpha}), X\big)-M\big(B_{\Z^d}(x_1,n^{\alpha})-e, X\big)=X_e-\max\{X_e,\tilde{Y}\},\]
where $Y$ (resp. $\tilde{Y}$) is the maximum weight associated with the edges
in the path, $\Gamma_1$ (resp. $\Gamma_2$) connecting $x_1$ and $x_2$ in an MST
of $B_{\Z^d}(n)-e$ (resp. $B_{\Z^d}(x_1,n^{\alpha})-e)$. Thus, $\E|\delta_e M|\leq\E|\tilde{Y}-Y|$.

By the minimax property of paths in MST (Lemma \ref{lem:mst minimax criterion}), $(\tilde{Y}-Y)$ is always nonnegative. Further,
\begin{align}\label{eqn:11}
\E(\tilde{Y}-Y)=\int_0^1\PR(Y<u<\tilde Y)\ du.
\end{align}
Note that $\{\bI_{X_f\leq u}: f\text{ is an edge of } B_{\Z^d}(n)\}$ is a collection of i.i.d. Bernoulli$(u)$ random variables. Declare the edge $f$ to be open at level $u$ if $X_f\leq u$, and consider the corresponding $u$-clusters.
On the set $\{Y<u<\tilde Y\}$, the $u$-clusters in $B_{\Z^d}(x_1, n^{\alpha})-e$ containing $x_1$ and $x_2$ are disjoint (since $\tilde Y>u$). However, $x_1$ and $x_2$ are connected in $B_{\Z^d}(n)-e$ by a path open at level $u$ (since $Y<u$). Hence the $u$-clusters in
$B_{\Z^d}(x_1, n^{\alpha})-e$ containing $x_1$ and $x_2$ both intersect $\partial^{\inn} B_{\Z^d}(x_1, n^{\alpha})$ . (In this case, part of $\Gamma_1$ lies outside $B_{\Z^d}(x_1,n^{\alpha})$; see Figure \ref{fig:twoarm}.) Thus
\[\PR(Y<u<\tilde Y)\leq\PR\left(e\underset{u}{\overset{2}{\leftrightsquigarrow}}B_{\Z^d}(x_1, n^{\alpha})-e\right).\]
We can now use estimates on probability of two-arm events to bound $\E(\tilde Y-Y)$.
Thus, for any small positive $\eps$, the integrand in \eqref{eqn:11} is bounded by
$c(\log(n)/n)^{1/2}$ for $u \in (p_c-\eps, p_c+\eps)$ (Lemma \ref{lem:cerf}), and
benefits from the exponential decay when $u\notin(p_c-\eps, p_c+\eps)$.

\begin{figure}[htb]
\begin{center}
\includegraphics[trim=2cm 12.8cm 2cm 2cm, clip=true, angle=0, scale=0.5]{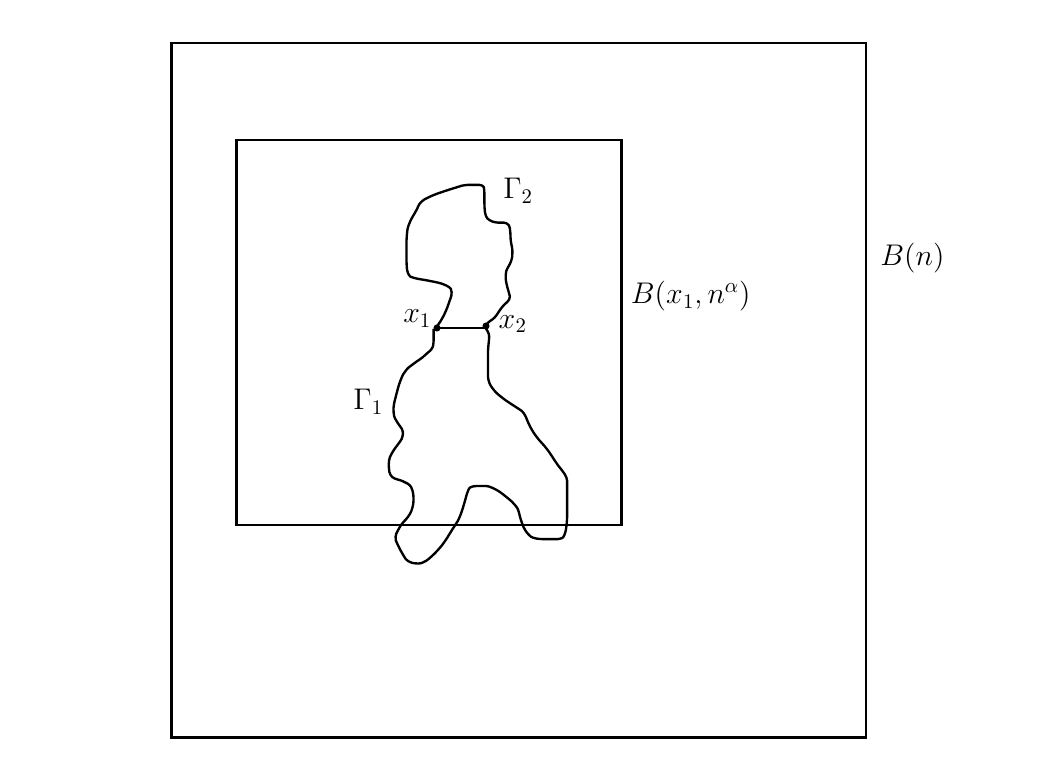}
\caption{The minimax paths connecting $x_1$ and $x_2$ when $\tilde{Y}>Y$.}
\label{fig:twoarm}
\end{center}
\end{figure}

For Euclidean MST, we start by dividing $B_{\R^d}(n)$ into cubes $\{Q\in\mathcal{Q}\}$ with
disjoint interiors having side length $s\in [1,2]$.
Consider a Poisson process $\poi$ in $\bR^d$ of intensity one and let $X_Q:=\poi\cap Q$ for any cube $Q$.
Set $X=(X_Q: Q\in\cQ)$, and let $X'$ be an independent copy of $Q$.
Consider a cube $Q_0\in\cQ$ with $d(c(Q_0), \partial B_{\R^d}(n))\geq n^{\alpha}$.
In line with the notation in Section \ref{sec:stein's-method},
$X^{Q_0}$ denotes the configuration in $B_{\R^d}(n)$ when the configuration inside $Q_0$ is $X_{Q_0}'$, and the configuration in
$B_{\R^d}(n)\setminus Q_0$ is given by $\cup_{Q\in\cQ\setminus Q_0} X_Q$. Similar to the discrete case, our aim then is to get
a bound on $\E|\Delta_{Q_0} M_n-\tilde{\Delta}_{Q_0} M_n|$, where
\[\Delta_{Q_0} M_n=M_{\bR^d}(X)-M_{\bR^d}(X^{Q_0}),\text{ and}\]
\[\tilde{\Delta}_{Q_0} M_n=M_{\bR^d}\big(X\cap B_{\R^d}(c(Q_0),n^{\alpha})\big)
-M_{\bR^d}\big(X^{Q_0}\cap B_{\R^d}(c(Q_0),n^{\alpha})\big).\]
This can also be reduced to getting a bound on the probability of
the two-arm event in the setup of continuum percolation. However,
since all possible edges between points are permitted, this step requires
a little work. We achieve this by introducing the concept of a ``wall'' (Definition \ref{def:wall})
and then using the add and delete algorithm. We will omit the details
of these steps from the proof sketch.

\section{Some results about Euclidean minimal spanning trees}\label{sec:mst}
In this section, the underlying space will always be $\R^d$, and we will simply write $B(\cdot,\cdot)$, $d(\cdot,\cdot)$, and $M(\cdot)$ instead of $B_{\R^d}(\cdot,\cdot)$, $d_{\R^d}(\cdot,\cdot)$, and $M_{\R^d}(\cdot)$.

When dealing with Euclidean minimal spanning trees, we would like
to have a criterion which ensures that if we fix a small cube, then
there are no ``long'' edges in the MST with one endpoint inside that cube.
Kesten and Lee \cite{kesten} used the idea of a ``separating set" to meet this purpose. (We will not define separating sets since we do not use them in this paper.) We generalize their ideas to define a ``wall" (see Definition \ref{def:wall} below). The reason behind this is that using the notion of separating sets in our proof will yield a weaker convergence rate than the one stated in Theorem \ref{thm:poissonmst}.

\begin{definition}\label{def:wall}
Suppose that $b>a$ are positive numbers and $x\in\R^d$ and let $K$ be a cube containing $B(x,a)$.
Further assume that $K\cap\partial B(x,b)\neq\emptyset$.
We say that a subset $\mathfrak{W}$ of $\R^d$ contains a $K$-wall around $B(x,a)$
in $B(x,b)$ if the following holds :
\begin{align}
&\text{For any }
p_1\in\partial B(x,a)\text{ and }p_2\in K\cap\partial B(x,b),\text{ the set }\nonumber\\
& K\cap\mathfrak{W}\cap S\left(p_1,3d(p_1, p_2)/4\right)\cap S\left(p_2, 3d(p_1, p_2)/4\right)\cap \{B(x,b)\setminus B(x,a)\}\nonumber\\
& \text{is nonempty}.\nonumber
\end{align}
If $B(x,b)\subset K$, we will simply say $\mathfrak{W}$ contains a wall around $B(x,a)$
in $B(x,b)$.
\end{definition}
\begin{figure}[htb]
\begin{center}
\includegraphics[trim=2cm 16.7cm 2cm 4cm, clip=true, angle=0, scale=0.6]{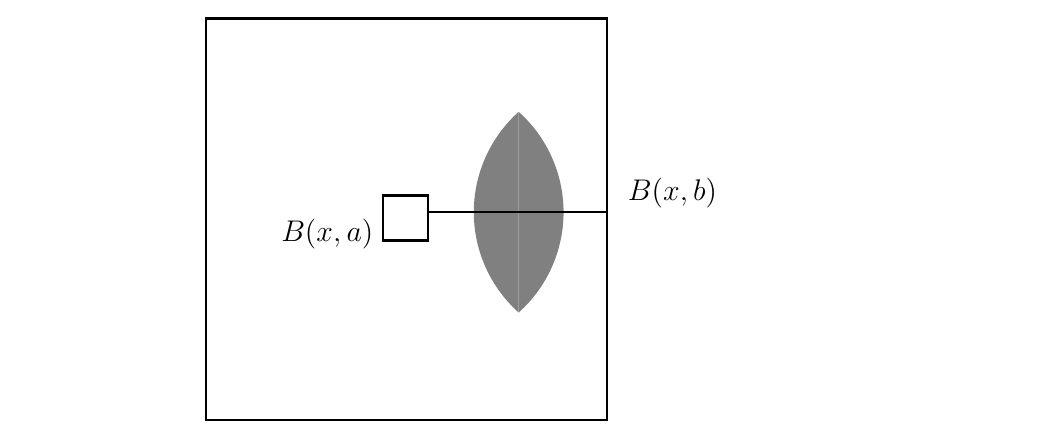}
\caption{For a wall to exist around $B(x,a)$ in $B(x,b)$, the shaded region must contain a point}
\label{fig:wall}
\end{center}
\end{figure}
The importance of this definition will be clear from the following lemma.
\begin{lemma}\label{lem:wall property A}
Let $a, b, x, K$ be as in Definition \ref{def:wall}.
Let $\omega$ be a finite set of points in $K$ and consider the complete
graph $(V,E)$ on $\omega$ with edge weights being the Euclidean length of edges.
If $\omega$ contains a $K$-wall around $B(x,a)$ in $B(x,b)$, then no edge in $E$
with one endpoint in $B(x,a)$ and other endpoint in $B(x,b)^c$ is included
in any MST of $(V,E)$.
\end{lemma}

\noindent\textbf{Proof of Lemma \ref{lem:wall property A}: }
Let $y_1,y_2$ be two points in $\omega$ such that $y_1\in B(x,a)$ and
$y_2\in B(x,b)^c$. Assume that $p_1\in\partial B(x,a)$ and $p_2\in\partial B(x,b)$
are points on the line segment $\overline{y_1y_2}$.
Since $\omega$ contains a $K$-wall around $B(x,a)$ in $B(x,b)$, we can find a point
$z$ such that
\[z\in\left[\omega\cap S\left(p_1, 3d(p_1, p_2)/4\right)\cap S\left(p_2, 3d(p_1, p_2)/4\right)\cap (B(x,b)\setminus B(x,a))\right].\]
Then
\begin{align*}
d(y_1, z)&\leq d(y_1, p_1)+d(p_1, z)\leq d(y_1, p_1)+3d(p_1, p_2)/4 \\
&< d(y_1, p_1)+d(p_1, y_2)=d(y_1, y_2).
\end{align*}
Similarly $d(z, y_2)<d(y_1, y_2)$. Hence, it follows from Lemma \ref{lem:mst minimax criterion} that $\overline{y_1 y_2}$ will not be
included in any minimal spanning tree of $(V,E)$.
\hfill$\blacksquare$

Next we show that a wall exists in a large annulus with high probability.
\begin{lemma}\label{lem:wall tail decay} Let $d\geq 2$ and $x\in\R^d$. As always we let
$\poi$ be a Poisson process of intensity one in $\R^d$. Then for any $a_0>0$,
there exist constants $c$ and $c'$ depending only on $a_0$ and $d$ such that the following holds:
for every $a\leq a_0$ and $b>a$,
\begin{align}
&\PR\big(\poi\textup{ does not contain a $B(n)$-wall around }B(x,a)\text{ in }B(x,b)\big)\leq c\exp(-c'b^d)\nonumber
\end{align}
for any $n$ for which $B(x,a)\subset B(n)$ and $B(n)\cap\partial B(x,b)\neq\emptyset$.
\end{lemma}
\noindent{\bf Proof:}
It suffices to prove the claim for large values of $b$, so let us start with the assumption
$b>4a_0+16$.

Cover $B(n)\cap\partial B(x,b)$ by $(d-1)$ dimensional cubes, $\{Q_i^{1}\}_{i\leq m_1}$
of diameter one. This can be done in a way so that the total number of cubes, $m_1$, is at most
$cb^{d-1}$. Similarly cover $\partial B(x,a)$ by $(d-1)$ dimensional cubes
$\{Q_i^{2}\}_{i\leq m_2}$ of diameter $\min(1, 2a\sqrt{d-1})$ so that the total number of cubes, $m_2$,
is at most $c\max(1, a_0^{d-1})$.

Let $p_1',p_2'$ be two points on $\partial B(x,a)$ and $B(n)\cap\partial B(x,b)$ respectively
and let $z'=(p_1'+p_2')/2$ be the midpoint of $\overline{p_1'p_2'}$. Let $p_1$ and $p_2$
be the centers of the cubes $Q_i^1$ and $Q_j^2$ such that $p_1'\in Q_i^1$ and $p_2'\in Q_j^2$.
Let $z=(p_1+p_2)/2$.

Consider $y'\in S(z',b/8)$. Then $\|y'-z'\|_{\infty}\leq b/8$ and hence
$$\|z'-x\|_{\infty}-b/8\leq \|y'-x\|_{\infty}\leq  \|z'-x\|_{\infty}+b/8.$$
Now,
\begin{align}
\|z'-x\|_{\infty}+\frac{b}{8}&= \bigg\|\frac{1}{2}(p_1'+p_2'-2x)\bigg\|_{\infty}+\frac{b}{8}
\leq\frac{a_0+b}{2}+\frac{b}{8}<b.\nonumber
\end{align}
Also
\[\|z'-x\|_{\infty}-\frac{b}{8}\geq \frac{b-a_0}{2}-\frac{b}{8}>a.\]
Hence $S(z',b/8)\subset B(x,b)\setminus B(x,a)$. Further, if $y\in S(z, b/16)$, then
\begin{align}
d(y,z')\leq \frac{b}{16}+d(z,z')&=\frac{b}{16}+\bigg\|\frac{p_1+p_2}{2}-\frac{p_1'+p_2'}{2}\bigg\|_{L^2}
\leq \frac{b}{16}+1\leq \frac{b}{8}.\nonumber
\end{align}
So $S(z,b/16)\subset S(z',b/8)\subset B(x,b)\setminus B(x,a)$.

If $y'\in S(z',b/8)$, then
\begin{align}
d(y',p_1')&\leq d(y',z')+d(z',p_1')
\leq\frac{b}{8}+\frac{d(p_1', p_2')}{2}
\leq \frac{3d(p_1', p_2')}{4}.\nonumber
\end{align}
The last inequality holds since
\[d(p_1', p_2')\geq b-a\geq b-a_0\geq b/2.\]
By a similar argument $d(y',p_2')\leq 3d(p_1', p_2')/4.$ Hence
\[ S\big(z',b/8\big)\subset S\bigg(p_1',3d\big(p_1', p_2'\big)/4\bigg)
\cap S\bigg(p_2',3d\big(p_1', p_2'\big)/4\bigg)\cap \big(B(x,b)\setminus B(x,a)\big).\]
Letting $\mathfrak{Leb}$ denote the Lebesgue measure, we note that $\mathfrak{Leb}(S(z,b/16)\cap B(n))\geq c' b^d$. So we can conclude that
\begin{align}
&\PR(\poi\text{ does not contain a $B(n)$-wall around }B(x,a)\text{ in }B(x,b))\nonumber\\
&\leq \PR\bigg(\text{For some }
i\leq m_1, j\leq m_2,\ \poi\cap B(n)\cap S\bigg(\frac{p_1+p_2}{2},\frac{b}{16}\bigg)=\emptyset\nonumber\\
&\phantom{mmm}\text{where }p_1\text{ and }p_2\text{ are the centers of }
Q_i^1\text{ and }Q_j^2\text{ respectively}\bigg)\nonumber\\
&\leq c \max(1,a_0^{d-1})b^{d-1}\exp(-c'b^d),\nonumber
\end{align}
where the last inequality follows from union bound. This proves the claim.
\hfill$\blacksquare$

The next lemma puts an upper bound on how much the weight of the MST changes when some points are removed.
\begin{lemma}\label{lem:D(A,B-A) type upper bound} Let $a, b, x, K$ be as in Definition \ref{def:wall}.
Let $\mathcal{A}$ and $\mathcal{B}$ be finite sets of points in $\R^d$ such that
$\mathcal{A}\subset B(x,a)$ and $\mathcal{B}\subset K\setminus B(x,a)$.
If $\mathcal{B}$ contains a $K$-wall around
$B(x,a)$ in $B(x,b)$, then
\[|M(\mathcal{A}\cup\mathcal{B})-M(\mathcal{B})|\leq c|\mathcal{A}|b\]
for some constant $c$ depending only on $d$. If such a wall does not exist, then
\[|M(\mathcal{A}\cup\mathcal{B})-M(\mathcal{B})|\leq c|\mathcal{A}|\mathrm{diameter}(K).\]
\end{lemma}
The proof of Lemma \ref{lem:D(A,B-A) type upper bound} is similar to the proof of \cite[Lemma 7]{kesten}. We include this argument for the reader's convenience. The proof depends on an auxiliary lemma.
\begin{lemma}\label{lem:max-degree}(\cite[Lemma 4]{aldous})
Consider an MST $\mathcal{T}$ on a finite subset $\omega$ of $\R^d$. Then there exists a constant $D_{\max}$ depending only on $d$ such that the degree, in $\mathcal{T}$,  of any point in $\omega$ is bounded by $D_{\max}$.
\end{lemma}
\noindent{\bf Proof of Lemma \ref{lem:D(A,B-A) type upper bound}:}
First we assume that $\mathcal{B}$ contains a $K$-wall around $B(x,a)$ in $B(x,b)$. Then $\mathcal{B}$ has a point, say $p$, in $B(x, b)\setminus B(x, a)$. Thus we can start from an MST on $\mathcal{B}$ and connect the points in $\mathcal{A}$ to $p$ to get a spanning tree on $\mathcal{A}\cup\mathcal{B}$. This gives
\[M(\mathcal{A}\cup\mathcal{B})\leq M(\mathcal{B})+|\mathcal{A}|b\sqrt{d}.\]
To get the other inequality, we start from an MST on $\mathcal{A}\cup\mathcal{B}$ and delete the points in $\mathcal{A}$ and all edges incident to them. By Lemma \ref{lem:wall property A}, each of these edges is contained in $B(x, b)$. By Lemma \ref{lem:max-degree}, we have deleted at most $D_{\max}|\mathcal{A}|$ many edges and this can create at most $(D_{\max}|\mathcal{A}|+1)$ many components. Each of these components has a point in $B(x, b)$. We can then connect these points to get a spanning tree on $\mathcal{B}$. This gives
\[M(\mathcal{B})\leq M(\mathcal{A}\cup\mathcal{B})+D_{\max}|\mathcal{A}|b\sqrt{d}.\]
The proof is similar when a wall does not exist.
\hfill$\blacksquare$

Lemma \ref{lem:D(A,B-A) type upper bound} gives us control over the tails of
$|M_{\bR^d}(\mathcal{A}\cup\mathcal{B})-M_{\bR^d}(\mathcal{B})|$. Using this we can show that all moments
of this quantity are finite when the configuration comes from a Poisson process.
\begin{lemma}\label{lem:moments are bounded}
For $x\in\R^d$, $0<a\leq a_0$ and $n\geq \max(2a_0,1)$
for which $B(x,a)\subset B(n)$, we have
\[\E\bigg(\bigg|M\big(\poi\cap B(n)\big)-M\big(\poi\cap\big[B(n)\setminus B(x,a)\big]\big)\bigg|^q\bigg)\leq C_q\text{ for every }q\geq 1.\]
The constant $C_q$ depends only on $a_0$, $d$ and $q$.
\end{lemma}
\noindent{\bf Proof:}
Define a random variable $Z$ as follows:
if there does not exist a $b\geq a$ such that $\partial B(x,b)\cap B(n)\neq\emptyset$
and $\poi$ contains
a $B(n)$-wall around $B(x,a)$ in $B(x,b)$,
set $Z=2\sqrt{d}n$; otherwise define $Z$ to be the infimum of all such~$b$.
From Lemma \ref{lem:wall tail decay},
\begin{align*}
\E (Z^q)&=
\int_0^{2\sqrt{d}n}qu^{q-1}\PR(Z>u)du
\\
&\leq a_0^q+ c\int_a^{n}qu^{q-1}\exp(-c'u^d)\ du+c(2\sqrt{d}n)^{q}\exp(-c'n^d).\nonumber
\end{align*}
The last expression is bounded by a constant depending only on $a_0$, $d$ and $q$.
Now, from Lemma \ref{lem:D(A,B-A) type upper bound}
\begin{align}
&\E\bigg(\bigg|M\big(\poi\cap B(n)\big)-M\big(\poi\cap\big[B(n)\setminus B(x,a)\big]\big)\bigg|^q\bigg)\nonumber \\
&\hskip20pt\leq c \E\big(Z\cdot|\poi\cap B(x,a)|\big)^q
\leq \frac{c}{2} \E\big[Z^{2q}+(|\poi\cap B(x,a)|)^{2q}\big],\nonumber
\end{align}
and this finishes the proof.\hfill$\blacksquare$

\section{Proofs of percolation estimates in the Euclidean setup}\label{sec:proofs}
In this section, the underlying space will always be $\R^d$, and all Poisson processes will have intensity one.
We will simply write $B(\cdot,\cdot)$ and $d(\cdot,\cdot)$ without referring to the ambient space.
Recall form Remark \ref{rem:only matters at criticality} that $r_c(d)$ denotes the critical radius for continuum percolation in $\R^d$ driven by a Poisson process with intensity one. When the dimension $d$ is clear, we will simply write $r_c$ instead of $r_c(d)$.

Before beginning the proof of Lemma \ref{lem:critical connectivity}, we collect two simple facts in the following lemma.
\begin{lemma}\label{lem:conditional var bound}
(i) Let $X_1,\hdots,X_n$ be independent random variables defined
on $(\Omega,\mathcal{A},\PR)$ taking values in some measurable
space $(\mathcal{X}, \mathcal{S})$. Let $f:\mathcal{X}^n\to \R$ be a bounded measurable function. Then
for any $A_1,\hdots,A_k\subset\{1,\hdots,n\}$ such that $A_i$ are pairwise disjoint,
\begin{equation}\label{eqn:conditional var bound}
\var\big(f(X_1,\hdots,X_n)\big)\geq \sum_{i=1}^k \var\big[\E\big(f(X_1,\hdots,X_n)\mid \{X_j\}_{j\in A_i}\big)\big].
\end{equation}
(ii) If $Y_1$ and $Y_2$ are independent and identically distributed real valued random variables such that $\E(Y_1^2)<\infty$, then
\begin{equation}\label{eqn:basic-var-identity}
\var(Y_1)=\frac{1}{2}\E\big(Y_1-Y_2\big)^2.
\end{equation}
\end{lemma}
\noindent{\bf Proof:}
\eqref{eqn:basic-var-identity} is a basic identity whose proof we will omit. To prove \eqref{eqn:conditional var bound},
without loss of generality, we can assume $\E (f(X_1,\hdots,X_n))=0$.  Let
\begin{align}
H=\big\{g\in L^2(\Omega,\mathcal{A},\PR):\int g=0\big\},\text{ and }
H_i=\big\{g\in H : g\text{ is }\sigma\big(\{X_j\}_{j\in A_i}\big)\text{ measurable}\big\}.\nonumber
\end{align}
Then under the natural inner product, $H$ is a Hilbert space and the $H_i$ are closed orthogonal subspaces of $H$.
\eqref{eqn:conditional var bound} follows upon observing that $\E\big(f(X_1,\hdots,X_n)\mid\{X_j\}_{j\in A_i}\big)$ is the projection of
$f(X_1,\hdots,X_n)$ on $H_i$.
\hfill$\blacksquare$


The following lemma plays a crucial role in the proof of Lemma \ref{lem:critical connectivity}.
\begin{lemma}\label{lem:burton-keane} Let $0<r_1<r_2<\infty$.
Fix two nonnegative numbers $s$ and $t$ such that $s+t>2r_2$.
Then there exist positive constants $c$ and $c'$ depending only on
$r_1,\ r_2$ and the dimension $d$ such that for every $m>100(s+t)$ and $r\in[r_1, r_2]$
\begin{equation}\label{eqn:burton-keane}
\PR\big(B(s)^{(t)}\underset{r}{\overset{3}{\longleftrightarrow}}B(m)\big)\leq c\cdot\exp\left(c'(s+t)\right)/m.
\end{equation}
For $z_1,z_2\in B(s)^{(t)}$, the same bound holds for
$\PR\big(B(s)^{(t)}\cup S(z_1,r)\underset{r}{\overset{3}{\longleftrightarrow}}B(m)\big)$ and
$\PR\big(B(s)^{(t)}\cup S(z_1,r)\cup S(z_2,r)\underset{r}{\overset{3}{\longleftrightarrow}}B(m)\big)$.
\end{lemma}
The proof of Lemma \ref{lem:burton-keane} will be given in Section \ref{sec:proof of burton-keane lemma}. We now proceed with


\subsection{Proof of Lemma \ref{lem:critical connectivity}}
Let us first prove the bound for $\PR(B(a)\underset{r}{\overset{2}{\longleftrightarrow}}B(n))$. The arguments are similar when we replace $B(a)$ by the other sets.
Fix $r\in[r_1,r_2]$. We write $\R^d$ as a union of cubes
\[\R^d=\bigcup_{k\in\Z^d}B_k\text{ where }B_k=2ak+B(a).\]
 Since $\PR(\poi\cap\partial B_k\neq \emptyset\text{ for some }
k\in\Z^d)=0$, we will assume that no Poisson point lies in any of the common interfaces
shared by two cubes.

Consider a sequence $a_n\to\infty$ such that $a_n=o(n)$ but $a_n=\Omega((\log\log n)^2)$ (so that
$a_n$ is large compared to $a$). We will fix the sequence $a_n$ later.
Define
\begin{align}
E=\big\{&\exists\text{ exactly one }r\text{-cluster }\C\text{ in }B(n)\text{ such that }\nonumber\\
    &\C^{(r)}\text{ intersects both }\partial (B(n)_{(r)})\text{ and }\partial B(a_n)\big\}\nonumber.
\end{align}
Let
\[\cL=\big\{k\in\Z^d : B_k\cap B\left(n\right)\neq\emptyset\big\},\
\text{ and }\
\cI=\big\{k\in\Z^d : B_k\cap B\left(a_n/3\right)\neq\emptyset\big\}.\]
Define $f:\prod_{k\in\cL}\mathfrak{X}(B_k)\to\R$ by
\[f\big((\omega_k:k\in\cL)\big)=\I_E\big(\cup_{k\in\cL}\omega_k\big).\]
Write
\[X_k=\poi\cap B_k,\ \text{ and }X=\big(X_k : k\in\cL\big).\]
It then follows from Lemma \ref{lem:conditional var bound} that
\begin{equation}\label{eqn:conditional var bound application}
\var\big(f(X)\big)\geq \sum_{i\in\cI}\var\big[\E\big(f(X)|X_i\big)\big].
\end{equation}

Consider another Poisson process $\poi'$ independent of $\poi$, and set
\[X_k'=\poi'\cap B_k,\ \text{ and }X'=\big(X_k' : k\in\cL\big).\]
Recall the notation $X^j$ from Section \ref{sec:stein's-method}. Define
\begin{align*}
\mathcal{S}_i :=\big\{\omega_i \in\mathfrak{X}(B_i):\
{B_i}_{(2r)}\subset\omega_i^{(r)}\big\},\  \text{ and }\
\mathcal{G}_i:=\big\{\omega_i \in\mathfrak{X}(B_i):\
{B_i}_{(2r)}\cap\omega_i^{(r)}=\emptyset\big\}.
\end{align*}
Then, for any fixed $i\in\mathcal{I}$,
\begin{align}
&\mathrm{Var}\big[\E\big(f(X)|X_i\big)\big] =
\frac{1}{2}\E\bigg[\bigg(\E\big(f(X)\mid X_i\big)-\E\big(f(X^i) \mid X_i'\big)\bigg)^2\bigg]\nonumber \\
&\hskip40pt\geq \frac{1}{2}\E\bigg[\bigg(\E\big(f(X) - f(X^i) \mid X_i, X_i'\big)\bigg)^2
\cdot\I(X_i\in\mathcal{S}_i,X_i'\in\mathcal{G}_i)\bigg],\label{eqn:lower bound on variance}
\end{align}
where the first step uses \eqref{eqn:basic-var-identity} and the fact that $\E\big(f(X)|X_i\big)$ and $\E\big(f(X^i)|X_i'\big)$ are independent and identically distributed.

Consider $i\in\cI$, $\omega\in\mathfrak{X}(B(n)\setminus B_i)$ and $\omega_i'\in\mathcal{G}_i$. Then $\omega^{(r)}$ and $(\omega_i')^{(r)}$ are disjoint. Thus, if $E$ holds when the configuration in $B_i$ is $\omega_i'$ and the configuration in $B(n)\setminus B_i$ is $\omega$, then $E$ continues to hold when $B_i$ is empty and the configuration in $B(n)\setminus B_i$ is $\omega$. Further, if the event $E$ holds with some configuration in $B(n)$, then $E$ continues to hold  with the configuration obtained by adding extra points inside $B(a_n/3)$.
Thus for any $\omega_i\in\mathcal{S}_i$ and $\omega_i'\in\mathcal{G}_i$,
\begin{align*}
&\big\{\omega\in\mathfrak{X}(B(n)\setminus B_i):\ \I_E\left(\omega\cup\omega_i'\right)=1\big\}
\subset \big\{\omega\in\mathfrak{X}(B(n)\setminus B_i):\ \I_E\left(\omega\cup\omega_i\right)=1\big\}.
\end{align*}
Therefore, if $X_i \in \mathcal{S}_i$ and $X_i'\in \mathcal{G}_i$, then
\begin{equation}\label{eqn:positivity}
f(X) - f(X^i)\geq 0.
\end{equation}
Now, for any  $\omega\in\mathfrak{X}(B(n)\setminus B_i)$ for which the event
\begin{align}
A_i:= \big\{& B_i\underset{r}{\overset{2}{\longleftrightarrow}}B(n),\text{ every }r\text{-cluster }
\C\text{ in }B(n)\setminus B_i\text{ for which }\C^{(r)}\nonumber\\
&\text{ intersects both }\partial B(a_n)\text{ and }\partial (B(n)_{(r)})\text{ has a point in } B_i^{(r)}\big\}\nonumber
\end{align}
is true, $\I_E(\omega\cup\omega_i)=1$ when $\omega_i\in\mathcal{S}_i$ and $\I_E(\omega\cup\omega_i')=0$ when $\omega_i'\in\mathcal{G}_i$. Consequently, if $X_i \in \mathcal{S}_i$ and $X_i'\in \mathcal{G}_i$ and $\I_{A_i}(\cup_{k\ne i}X_k)=1$, then
\[
f(X) - f(X^i)=1.
\]
Hence from \eqref{eqn:lower bound on variance} and \eqref{eqn:positivity},
\begin{align}\label{eqn:lower bound on variance I}
\var\big[\E\big(f(X)\mid X_i\big)\big]&\geq \frac{1}{2}\PR(A_i)^2\cdot\PR(X_i\in\mathcal{S}_i)\cdot\PR(X_i'\in\mathcal{G}_i)\\
&\geq \frac{1}{2}\PR(A_i)^2\exp(-c a^{d-1}).\nonumber
\end{align}
The constant depends on $d$ and $r_1$ only.

For $i\in\mathcal{I}$, we also have
\begin{align}\label{eqn:lower bound on variance II}
\PR(A_i)&\geq \PR\big( B_i\underset{r}{\overset{2}{\longleftrightarrow}}B(n);\text{ any }r\text{-cluster }
\C\text{ in } B(n)\setminus B_i\\
&\hskip25pt\text{for which }\C^{(r)}\text{ intersects both }\partial B(c(B_i),2a_n)\nonumber\\
&\hskip25pt\text{and }\partial (B(n)_{(r)})\text{ has a point in }B_i^{(r)}\big)\nonumber\\
&\geq\PR\big( B_i\underset{r}{\overset{2}{\longleftrightarrow}}B(c(B_i),2n);\text{ if }\C\text{ is an }
r\text{-cluster in} \nonumber\\
&\hskip25pt
B\big(c(B_i),2n\big)\setminus B_i\text{ then every connected component}\nonumber\\
&\hskip25pt
\text{of }\bigg(\C\cap B\big(c(B_i),n/2\big)\bigg)^{(r)}
\text{ that intersects both }\nonumber\\
&\hskip27pt
\partial B\big(c(B_i),2a_n\big)\text{ and }\partial \big(B\big(c(B_i),n/2\big)_{(r)}\big)\text{ also intersects }\partial B_i\big).\nonumber
\end{align}
Define the event
\begin{align}
F=&\big\{B_0\underset{r}{\overset{2}{\longleftrightarrow}}B(2n),\text{ if }\C\text{ is an }
r\text{-cluster in } B(2n)\setminus B_0\nonumber\\
&
\text{ then every connected component of }\left(\C\cap B(n/2)\right)^{(r)}\text{ that }\nonumber\\
&
\text{ intersects both }\partial B(2a_n)\text{ and }\partial \left(B\left(n/2\right)_{(r)}\right)
\text{ also intersects }\partial B_0\big\}.\nonumber
\end{align}
From (\ref{eqn:conditional var bound application}), (\ref{eqn:lower bound on variance I}), (\ref{eqn:lower bound on variance II}) and translational invariance, we get
 \begin{equation}\label{eqn:bound on P(F)}
 \PR(F)\leq \frac{c\exp(c' a^{d-1})}{\sqrt{|\mathcal{I}|}}
 \leq \frac{c''\exp(c' a^{d-1})a^{d/2}}{a_n^{d/2}}.
 \end{equation}
Here we have used the fact that $\mathrm{Var}(f(X))\leq 1/4$ and
$|\mathcal{I}|=\Theta((a_n/a)^d)$.

On the event $\{B_0\underset{r}{\overset{2}{\longleftrightarrow}}B(2n)\}\cap F^c$,
we can find two disjoint $r$-clusters $\C_1, \C_2$ in $B(2n)\setminus B_0$ and
an $r$-cluster $\overline{\C}$ (which may be the same as one of the $r$-clusters
$\C_1, \C_2$) in $B(2n)\setminus B_0$ such that
\begin{enumeratei}
\item each of $\C_1$  and $\C_2$ has a point in $B_0^{(r)}$ and a point in $B(2n)_{(2r)}$,
\item there is an $r$-cluster in $B(n/2)\setminus B_0$, call it $\overline{\C}'$, which
 is contained in $\overline{\C}\cap B(n/2)$, such that $\overline{\C}'$ has a point in $B(2a_n)^{(r)}$
 and a point in $B(n/2)_{(2r)}$ but does not have a point in $B_0^{(r)}$.
\end{enumeratei}
So we can find two disjoint $r$-clusters $\C_1'$ and $\C_2'$ in $B(n/2)\setminus B_0$
that are contained in $\C_1\cap B(n/2)$ and $\C_2\cap B(n/2)$ respectively such that $\C_1'$
and $\C_2'$ satisfy the requirements for $\{B_0\underset{r}{\overset{2}{\longleftrightarrow}}B(n/2)\}$
to be true. Further, $\overline{\C}'$ is different from $\C_1'$ and $\C_2'$ since
$\overline{\C}'$ does not have a point in $B_0^{(r)}$. Hence the restrictions of
$\overline{\C}'$, $\C_1'$ and $\C_2'$ to $B(n/2)\setminus B(2a_n)$ will contain
three disjoint $r$-clusters satisfying the requirements for
$\{B(2a_n)\underset{r}{\overset{3}{\longrightarrow}}B(n/2)\}$ to be true.

Hence we have
\begin{equation}\label{eqn:decomposition into two terms}
\PR\big(B_0\underset{r}{\overset{2}{\longleftrightarrow}}B(2n)\big)\leq \PR(F)+\PR\big(B(2a_n)\underset{r}{\overset{3}{\longrightarrow}}B\left(n/2\right)\big).
\end{equation}
All we need now is an upper bound for the second term on the right side. We would like to apply a
Burton-Keane type argument to get a bound for this term.

Assume that $\C_1,\ \C_2$ and $\C_3$ are three disjoint $r$-clusters in $B(n/2)\setminus B(2a_n)$
such that $\C_j^{(r)}$ intersects both $B(n/2)_{(r)}$ and $B(2a_n)^{(r)}$ and let $x_j$ be the
point in $\C_j$ closest to $B(2a_n)$ for $j=1,2,3$.

If $x_j\in B(2a_n)^{(r)}$ for every $j$, then
$B(2a_n)\underset{r}{\overset{3}{\longleftrightarrow}}B\left(n/2\right)$ holds true, and if
$x_j\in B(2a_n)^{(2r)}\setminus B(2a_n)^{(r)}$ for every $j$, then
$B(2a_n)^{(r)}\underset{r}{\overset{3}{\longleftrightarrow}}B\left(n/2\right)$ holds true.

Assume now that the event
\begin{align*}
\bigg\{B(2a_n)\underset{r}{\overset{3}{\longrightarrow}}B\left(n/2\right)\bigg\}
\cap \bigg(\big\{B(2a_n)\underset{r}{\overset{3}{\longleftrightarrow}}B\left(n/2\right)\big\}
\cup\big\{B(2a_n)^{(r)}\underset{r}{\overset{3}{\longleftrightarrow}}B\left(n/2\right)\big\}\bigg)^c
\end{align*}
is true. Then the number of $x_i$'s in $B(2a_n)^{(r)}\setminus B(2a_n)$ is one or two.

Let us assume that $x_1, x_2\in B(2a_n)^{(r)}$ and $x_3\in B(2a_n)^{(2r)}\setminus B(2a_n)^{(r)}$
(the other possibilities can be handled similarly).
We can find a sequence of points  $z_1^{(j)},\hdots,z_{k_j}^{(j)}$ in $\C_j$ for $j=1,2$ such that
\begin{align}
&(i)\ z_1^{(j)}\in B(2a_n)^{(r)}\text{ and }z_i^{(j)}\notin B(2a_n)^{(r)}\text{ if }i\geq 2,
\nonumber\\
&(ii)\ z_{k_j}^{(j)}\in B(n/2)_{(2r)},\nonumber\\
&(iii)\ d(z_i^{(j)}, z_{i+1}^{(j)})\leq 2r\text{ for } 1\leq i\leq k_j-1\text{ and}\nonumber\\
&(iv)\ d(z_i^{(j)}, z_{i'}^{(j)})>2r\text{ whenever }i'\geq i+2.\nonumber
\end{align}
Let $\C_j'(\subset\C_j)$ be the $r$-cluster in $B(n/2)\setminus B(2a_n)^{(r)}$ containing
$\{z_2^{(j)},\hdots,z_{k_j}^{(j)}\}$. Note that
\[
\max_{j=1,2} d(z_1^{(j)},z_2^{(j)})>r,
\]
because otherwise the event
$\{B(2a_n)^{(r)}\underset{r}{\overset{3}{\longleftrightarrow}}B\left(n/2\right)\}$ will be true
(the $r$-clusters $\C_1',\C_2'$ and $\C_3$ will satisfy the requirements).
If $\min_{j=1,2} d(z_1^{(j)},z_2^{(j)})\leq r$ then $E_1(z_1^{(1)})\cup E_1(z_1^{(2)})$ holds,
where
\[E_1(x):=\left\{B(2a_n)^{(r)}\cup S(x,r)
\underset{r}{\overset{3}{\longleftrightarrow}}B\left(n/2\right)\right\}\]
for $x\in B(2a_n)^{(r)}$
and if $\min_{j=1,2} d(z_1^{(j)},z_2^{(j)})> r$ then the event
\[E_2(z_1^{(1)},z_1^{(2)}):=\left\{B(2a_n)^{(r)}\cup S(z_1^{(1)},r)\cup S(z_1^{(2)},r)\underset{r}{\overset{3}{\longleftrightarrow}}B\left(n/2\right)\right\},\]
holds; in each case, $\C_3$ and the appropriate $r$-clusters containing the points $\{z_2^{(j)},\hdots,z_{k_j}^{(j)}\}$
($j=1,2$) satisfying the requirements. Hence
\begin{align}
&\PR\left(B(2a_n)\underset{r}{\overset{3}{\longrightarrow}}B\left(n/2\right)\right)\\
&\leq  \PR\left(B(2a_n)\underset{r}{\overset{3}{\longleftrightarrow}}B\left(n/2\right)\right)
+\PR\left(B(2a_n)^{(r)}\underset{r}{\overset{3}{\longleftrightarrow}}B\left(n/2\right)\right)\nonumber\\
&\quad+\PR\bigg(\exists x,y\in\poi\cap\big(B(2a_n)^{(r)}\setminus B(2a_n)\big)\text{ such that }
x\neq y\text{ and }E_2(x,y)\text{ holds}\bigg)\nonumber\\
&\quad +\PR\bigg(\exists x\in\poi\cap\big(B(2a_n)^{(r)}\setminus B(2a_n)\big)\text{ such that }
E_1(x)\text{ holds}\bigg).\nonumber
\end{align}
This gives
\begin{align}\label{eqn:just before burton-keane part starts}
&\PR\left(B(2a_n)\underset{r}{\overset{3}{\longrightarrow}}B\left(n/2\right)\right)\\
&\leq  \PR\left(B(2a_n)\underset{r}{\overset{3}{\longleftrightarrow}}B\left(n/2\right)\right)
+\PR\left(B(2a_n)^{(r)}\underset{r}{\overset{3}{\longleftrightarrow}}B\left(n/2\right)\right)\nonumber\\
&\hskip10pt+\E\bigg|\poi\cap(B(2a_n)^{(r)}\setminus B(2a_n))\bigg|^2\sup\nolimits_1\PR(E_2(x,y))\nonumber\\
&\hskip20pt +\E\bigg|\poi\cap(B(2a_n)^{(r)}\setminus B(2a_n))\bigg|\ \sup\nolimits_2\PR(E_1(x)),\nonumber
\end{align}
where $\sup_1$ (resp. $\sup_2$) is supremum taken over all $x,y$ (resp. $x$) in
$B(2a_n)^{(r)}\setminus B(2a_n)$. Lemma \ref{lem:burton-keane} helps us in estimating
$\PR(E_2(x,y))$ and $\PR(E_1(x))$.

From \eqref{eqn:bound on P(F)}, (\ref{eqn:decomposition into two terms}),
\eqref{eqn:just before burton-keane part starts} and Lemma \ref{lem:burton-keane} we get
\begin{align}\label{eqn:lemma 4 final step}
\PR(B_0\underset{r}{\overset{2}{\longleftrightarrow}}B(2n))\leq
c\bigg(\exp(c' a^{d-1})\frac{a^{d/2}}{a_n^{d/2}}
+\exp(c''a_n)\frac{a_n^{3d-2}}{n}\bigg).
\end{align}
We choose $a_n$ so that
$c'' a_n=\frac{1}{2} \log n$,
plug this into (\ref{eqn:lemma 4 final step}) and finally replace $n$ by $n/2$
to get (\ref{eqn:critical connectivity}).

If we replace $B(a)$ in \eqref{eqn:critical connectivity} by, say, $K=B(a)^{(r)}\cup S(x,r)$, then define
$B_k:=2(a+2r)k+K$ so that the sets $B_k$ remain disjoint. Define $\mathcal{I}$
as before and think of $f$ as a function of the configurations inside $\{B_k\}_{k\in\mathcal{I}}$
and the configuration in the complement of $\cup_{k\in\mathcal{I}} B_k$. The rest of the
proof can be carried out by following the same arguments as before. This concludes the proof of Lemma \ref{lem:critical connectivity}.
\subsection{Proof of Lemma \ref{lem:burton-keane}}\label{sec:proof of burton-keane lemma}
We start with some auxiliary lemmas. The following lemma is a restatement of Lemma
3.2 in \cite{royII}.
\begin{lemma}\label{lem:rahul roy combinatorial lemma}
Let $R$ be a finite non empty subset of a set $S$. Assume further that\\
(I) for every $r\in R$, there exist pairwise disjoint subsets
(which we call ``branches'') $C_{r}^{(1)},\hdots,C_{r}^{(m_r)}$
of $S$ and a positive integer $k$ such that
\begin{align*}
& (Ia)\ m_r\geq 3,\hskip250pt\\
& (Ib)\ r\notin C_{r}^{(i)}\text{ for }i\leq m_r,\text{ and}\\
& (Ic)\ |C_{r}^{(i)}|\geq k\text{ for }i\leq m_r;
\end{align*}
(II) for all $r,r'\in R$, either
\begin{align*}
& (IIa)\ \left(\cup_{j\leq m_r}C_{r}^{(j)}\cup\{r\}\right)\cap
\left(\cup_{i\leq m_{r'}}C_{r'}^{(i)}\cup\{r'\}\right)=\emptyset\text{ or}\\
& (IIb)\ \left(\cup_{j\leq m_r}C_{r}^{(j)}\cup\{r\}\right)\setminus C_r^{(j_0)}\subset C_{r'}^{(i_0)}\text{ and}\\
&\phantom{(IIb)}\left(\cup_{i\leq m_{r'}}C_{r'}^{(i)}\cup\{r'\}\right)\setminus C_{r'}^{(i_0)}\subset C_r^{(j_0)}\text{ for some }
i_0\leq m_{r'}\text{ and }j_0\leq m_r.
\end{align*}
Then $|S|\geq k|R|$.
\end{lemma}
Let $K\subset B(m)$ be a translate of $B(s)^{(t)}$ where $s, t, m$ are as in the statement of Lemma \ref{lem:burton-keane}.
\begin{figure}[htb]
\begin{center}
\includegraphics[trim=2cm 13.8cm 2cm 2cm, clip=true, angle=0, scale=0.5]{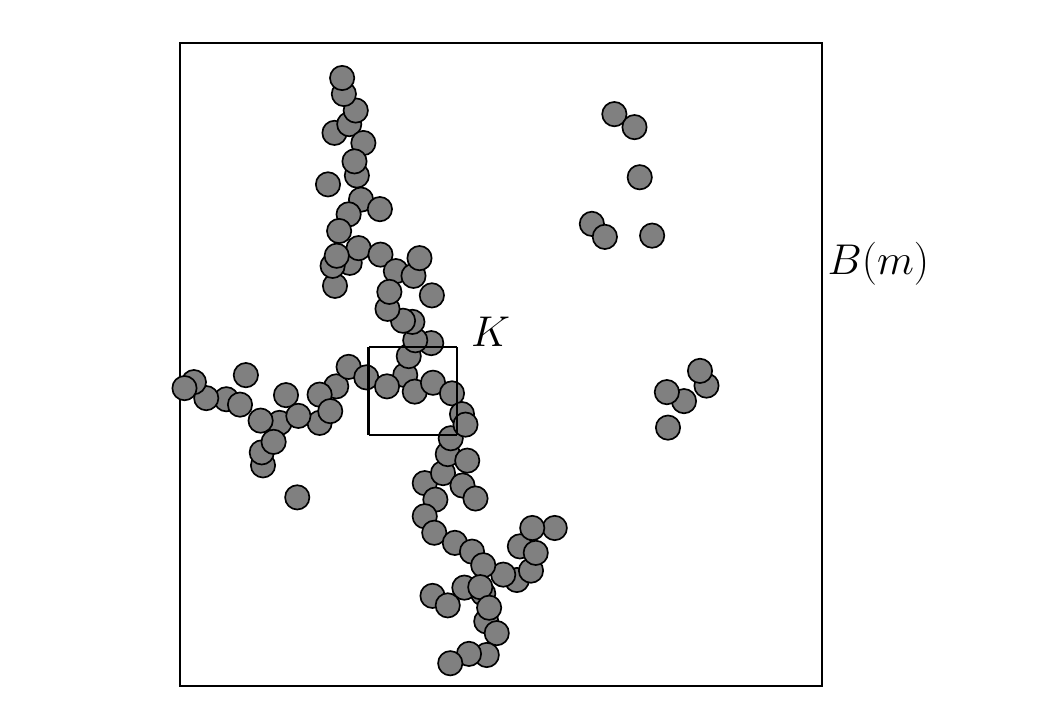}
\caption{$K$ is a trifurcation box in $B(m)$.}
\label{fig:trifurcation box}
\end{center}
\end{figure}
We will say that $K$ is a trifurcation box
in $B(m)$ (in short $``K\text{ T-box in }B(m)"$)
at level $r$ (Figure~\ref{fig:trifurcation box}) if
\begin{align}
& (i)\text{ there is an }r\text{-cluster }\C\text{ in }B(m)\text{ with }\C\cap K\neq\emptyset\text{ and}\nonumber\\
& (ii)\ \C\cap K^c\text{ contains at least three disjoint }r\text{-clusters in }B(m)\setminus K\nonumber\\
& \phantom{(iii)}\text{each having a point in }
B(m)_{(2r)}.\nonumber
\end{align}
Let us define
$$\mathcal{T}:=\{j\in\Z^d: 4(s+t)j+B(s)^{(t)}\subset B(m/4)\}$$
and denote $4(s+t)j+B(s)^{(t)}$ by $K_j$ for $j\in\mathcal{T}$.
Then we have the following
\begin{lemma}\label{lem:burton-keaneI}
There exists a positive constant $c$ depending only on $r_2$ such that
\begin{equation}\label{eqn:burton-keane III}|
\{ \poi \cap B(m /2) \}|\geq cm|\{j\in\mathcal{T}: K_j\text{ T-box in }B(m/2)\}|.
\end{equation}
\end{lemma}
\noindent\textbf{Proof:}
Set $S=\poi\cap B(m/2)$.
If $K_j$ is a trifurcation box in $B(m/2)$ for some $j\in\mathcal{T}$, then there is an $r$-cluster $\C_j$
in $B(m/2)$ such that there is a point $r_j$ in $\C_j\cap K_j$. Further,
$\C_j\cap B(m/2)\setminus K_j$ contains $m_j(\geq 3)$ disjoint $r$-clusters, say
$\C_j^{(1)},\hdots,\C_j^{(m_j)}$, each having a point in  $B(m/2)_{(2r)}$.
Call these clusters the ``branches'' of $r_j$.
Set $R=\{r_j:\ j\in\mathcal{T}, K_j\text{ T-box in }B(m/2)\}$.

For any $r_j,r_{j'}$ in $R$, condition $(IIa)$ of Lemma \ref{lem:rahul roy combinatorial lemma}
holds if $\C_j$ and $\C_{j'}$ are disjoint and condition $(IIb)$ holds otherwise.
Also
$$|\C_{r_j}^{(i)}|\geq\frac{m/4-2r_2}{2r_2}\geq cm$$
for every $r_j\in R$ and $i\leq m_j$.
Hence an application of Lemma \ref{lem:rahul roy combinatorial lemma} yields the result.\hfill$\blacksquare$

We are now ready to prove Lemma \ref{lem:burton-keane}. Note that
\begin{align}\label{eqn:burton-keane I}
\PR\big(B(s)^{(t)}\text{ T-box in }B(m)\big)\geq & \PR\big(B(s)^{(t)}\underset{r}{\overset{3}{\longleftrightarrow}}B(m)\big)\\
&\cdot\PR\left(B(s)^{(t)}\text{ T-box in }B(m)\bigg|B(s)^{(t)}
\underset{r}{\overset{3}{\longleftrightarrow}}B(m)\right).\nonumber
\end{align}
Now, given any $\eta\in\mathfrak{X}(B(m)\setminus B(s)^{(t)})$ for which the event
\[
A:=\big\{B(s)^{(t)}\underset{r}{\overset{3}{\longleftrightarrow}}B(m)\big\}
\]
is true,
we can ensure that the event $\{B(s)^{(t)}\text{ T-box in }B(m)\}$ happens just by placing
enough Poisson points inside $B(s)^{(t)}$ so that
at least three of the $r$-clusters in $B(m)\setminus B(s)^{(t)}$ satisfying the requirements
for $A$ to be true get connected to form a single component. Since this can be done
by placing at least one Poisson point in each of at most $6d^{3/2}(s+t)/r_1$ cubes
(of side length $r_1/\sqrt{d}$) inside $B(s)^{(t)}$,
\[\PR\big(B(s)^{(t)}\text{ T-box in }B(m)\mid B(s)^{(t)}
\underset{r}{\overset{3}{\longleftrightarrow}}B(m)\big)\geq \exp(-c(s+t))\]
for a positive universal constant $c$ depending only on $r_1$ and $d$.
Plugging this into (\ref{eqn:burton-keane I}), we get
\begin{align}\label{eqn:burton-keane II}
\PR\big(B(s)^{(t)}
\underset{r}{\overset{3}{\longleftrightarrow}}B(m)\big)\leq
\exp(c(s+t))
\cdot\PR\big(B(s)^{(t)}\text{ T-box in }B(m)\big).
\end{align}

Taking expectation in (\ref{eqn:burton-keane III}), we get
\begin{align}
\frac{m^{d-1}}{c2^d}\geq
&\sum_{j\in\mathcal{T}}\PR(K_j\text{ T-box in }B(m/2))\nonumber\\
&\hskip10pt\geq \sum_{j\in\mathcal{T}}\PR\big(K_j\text{ T-box in }4(s+t)j+B(m)\big)\nonumber.
\end{align}
By translational invariance and the fact that $|\mathcal{T}|\cdot(s+t)^d=\Theta(m^d)$, we get
\begin{align}\label{eqn:burton-keane IV}
c' m^{d-1}
\geq \frac{m^d}{(s+t)^d}\PR\big(B(s)^{(t)}\text{ T-box in }B(m)\big)
\end{align}
and (\ref{eqn:burton-keane}) follows if we plug this in (\ref{eqn:burton-keane II}).

The same type of arguments work when $B(s)^{(t)}$ is replaced by the other sets,
so we do not repeat them.
\subsection{Estimates in different regimes}\label{sec:proof of poissonmst}
We now collect the estimates on $\PR(B(a)\underset{r}{\overset{2}{\longrightarrow}}B(n))$ in different regimes together in the following lemma.
\begin{lemma}\label{lem:->connectivity bound}
For positive numbers $r_1,\ r_2$ satisfying $r_1<r_c(d)<r_2$ and $n\geq 2$, we have the following estimates.\\
(i) When $d=2$ and $a\in[1/2,\log n]$,
\begin{equation*}
\PR(B(a)\underset{r}{\overset{2}{\longrightarrow}}B(n))\leq\left\{
\begin{array}{l}
c_{10}\exp(-c_{11}n), \text{ if }r\leq r_1,\\
c_{12}/n^{\beta},\text{ if } r_1<r\leq(\log n)^2,
\end{array}
\right.
\end{equation*}
where $c_{10}$ and $c_{11}$ depend only on $r_1$, and $c_{12}$ and $\beta$ are universal positive constants.
\\
(ii) When $d\geq 3$ and $a\in(1/2,(\log\log n)^{1/(d-1/2)})$,
\begin{equation}\label{eqn:->connectivity bound d=3}
\PR(B(a)\underset{r}{\overset{2}{\longrightarrow}}B(n))\leq\left\{
\begin{array}{l}
c_{13}\exp(-c_{14}n) \text{ if }r\leq r_1,\\
c_{15}\frac{\exp(c_{16}a^{d-1})}{\left(\log n\right)^{d/2}}\text{ if }r\in[r_1,r_2],\\
c_{17}\exp(-c_{18} n)\text{ if }r_2\leq r\leq n/8.
\end{array}
\right.
\end{equation}
The constants appearing here depend only on $r_1$, $r_2$ and $d$.
\end{lemma}
\noindent{\bf Proof:}
The proof can be divided into different parts.
\begin{enumerateA}
\item{$\mvr\leq \mvr_1$ \textbf{and} ${\mvd\geq\mvtwo}$:} Note that for any $r>0$ and $d\geq 2$,
\begin{align}\label{eqn:->connectivity bound I}
\{B(a)\underset{r}{\overset{2}{\longrightarrow}}B(n)\}&\subset\{B(a)\underset{r}{\overset{1}{\longrightarrow}}B(n)\}\\
&\subset \{B(a)\underset{r}{\overset{1}{\longleftrightarrow}}B(n)\}
\cup\{B(a)^{(r)}\underset{r}{\overset{1}{\longleftrightarrow}}B(n)\}.\nonumber
\end{align}
That the last inclusion holds can be seen as follows. Consider an $r$-cluster $\C$ in $B(n)\setminus B(a)$ which has a point in both
$B(a)^{(2r)}$ and $B(n)_{(2r)}$ and let $x\in \C$ be the point closest to $B(a)$. If $x\in B(a)^{(r)}$ then
$\{B(a)\underset{r}{\overset{1}{\longleftrightarrow}}B(n)\}$ is true and if
$x\in B(a)^{(2r)}\setminus B(a)^{(r)}$ then
$\{B(a)^{(r)}\underset{r}{\overset{1}{\longleftrightarrow}}B(n)\}$ is true.

For any $r\leq r_1$ and $d\geq 2$,
$\{B(a)\underset{r}{\overset{1}{\longleftrightarrow}}B(n)\}
\subset \{B(a)\underset{r_1}{\overset{1}{\longleftrightarrow}}B(n)\}$ and a similar statement holds
if we replace $B(a)$ by $B(a)^{(r)}$.
If we fix a configuration in $B(n)\setminus B(a)$ (resp. $B(n)\setminus B(a)^{(r)}$) for which
$\{B(a)\underset{r_1}{\overset{1}{\longleftrightarrow}}B(n)\}$ (resp. $\{B(a)^{(r)}\underset{r_1}{\overset{1}{\longleftrightarrow}}B(n)\}$)
holds, we can connect any of the corresponding clusters to the origin by placing at least one Poisson point
in at most $c(a+r_2)/r_1$ many cubes inside $B(a)$
(resp. $B(a)^{(r)}$) each of side length $\min(2a,r_1/\sqrt{d})$. Thus, if $\mu_{\poi}$ is the probability
measure corresponding to a Poisson process of intensity one with an extra point added at the origin, then
\[\mu_{\poi}\bigg(\mathrm{diameter}(\C_0)\geq n\text{ at level }r_1\ \big|\ B(a)\underset{r_1}{\overset{1}{\longleftrightarrow}}B(n)\bigg)\geq c\exp(-c'a),\]
$\C_0$ being the occupied component containing the origin.
A similar inequality holds for
$\mu_{\poi}(\mathrm{diameter}(\C_0)\geq n\text{ at level }r_1|B(a)^{(r)}\underset{r_1}{\overset{1}{\longleftrightarrow}}B(n))$.
Hence,
from \eqref{eqn:->connectivity bound I}, we get
\begin{align}
\PR\big(B(a)\underset{r}{\overset{2}{\longrightarrow}}B(n)\big)& \leq
\PR\big(B(a)\underset{r_1}{\overset{1}{\longleftrightarrow}}B(n)\big)
+\PR\big(B(a)^{(r)}\underset{r_1}{\overset{1}{\longleftrightarrow}}B(n)\big)\nonumber\\
&\hskip20pt \leq c \exp(c'a)\cdot\mu_{\poi}\big(\mathrm{diameter}(\C_0)\geq n\text{ at level }r_1\big)\nonumber\\
&\hskip40pt \leq c \exp(c'a)\exp(-c''n).\nonumber
\end{align}
 The last inequality is just an application
of \cite[Equation (3.60)]{royII}.

\item{${\mvr\in[\mvr_1,\mvr_c]}$ \textbf{and} ${\mvd= \mvtwo}$:} In this case
\begin{align*}
\PR(B(a)\underset{r}{\overset{2}{\longrightarrow}}B(n))&\leq
\PR(B(a)\underset{r_c}{\overset{1}{\longleftrightarrow}}B(n))
+\PR(B(a)^{(r)}\underset{r_c}{\overset{1}{\longleftrightarrow}}B(n))\\
&\hskip30pt\leq c/n^{\theta},\text{ for some }\theta>0.\nonumber
\end{align*}
The last inequality holds because of the following reason. First note that
\begin{align}\label{eqn:vacant-crossing}
g_{\ell}(r_c):=\PR(\exists\text{ a vacant left-right crossing of }[0,\ell]\times[0,3\ell]\text{ at level }r_c)\geq \kappa_0:=\frac{1}{(9e)^{122}}
\end{align}
for every $\ell\geq r_c$. (This is true since otherwise there exists $\ell^{\star}\geq r_c$ for which \eqref{eqn:vacant-crossing} fails. By continuity of the function $g_{\ell^\star}$, we will be able to find $r<r_c$ such that $g_{\ell^\star}(r)< (9e)^{-122}$. Then by Lemma 4.1 of \cite{royII}, the vacant component containing the origin is bounded almost surely which leads to a contradiction since $r<r_c$.)
Now \eqref{eqn:vacant-crossing} together with  Lemma 4.4 of \cite{royII} and the RSW lemma
for vacant crossings (see \cite{royI} or Theorem 4.2 in \cite{royII})
will yield
\[\PR\big(\exists\text{ a vacant left-right crossing of }[0,3\ell]\times[0,\ell]\text{ at level }r_c\big)\geq\delta\]
for a positive constant $\delta$ and every $\ell$ bigger than a fixed threshold $\ell_0$.
It then follows from standard arguments that with probability at least $1-c/n^{\theta}$,
a vacant circuit around $B(a+r_c)$ exists in $B(n)$ at level $r_c$. Hence we get the
desired upper bound on $\PR(B(a)\underset{r}{\overset{2}{\longrightarrow}}B(n))$
for $r\in[r_1,r_c]$.

\item{${\mvr\geq \mvr_c}$ \textbf{and} ${\mvd=\mvtwo}$:} In this case the polynomial decay of $\PR(B(a)\underset{r}{\overset{2}{\longrightarrow}}B(n))$
follows from the existence of occupied ``circuits'' at level $r_c$ around $B(a)$.
The argument for this is also standard. We will give an outline in Appendix \ref{Appendix}.

\item{${\mvr\in [\mvr_1, \mvr_2]}$ \textbf{and} ${\mvd\geq\mvthree}$:}
Fix $r\in[r_1,r_2]$ and assume that $\{B(a)\underset{r}{\overset{2}{\longrightarrow}}B(n)\}$ holds. Take any two
disjoint clusters $\C_1$ and $\C_2$ in $B(n)\setminus B(a)$ each having a point in $B(a)^{(2r)}$ and $B(n)_{(2r)}$
and let $x_j\in\C_j$ be the point closest to $B(a)$. If $x_j\in B(a)^{(r)}$ for $j=1,\ 2$ then the
event $\{B(a)\underset{r}{\overset{2}{\longleftrightarrow}}B(n)\}$ is true, and if
$x_j\in B(a)^{(2r)}\setminus B(a)^{(r)}$ for $j=1,\ 2$ then the event
$\{B(a)^{(r)}\underset{r}{\overset{2}{\longleftrightarrow}}B(n)\}$ is true.

Now, assume that the event
\[\{B(a)\underset{r}{\overset{2}{\longrightarrow}}B(n)\}\cap
\left[\{B(a)\underset{r}{\overset{2}{\longleftrightarrow}}B(n)\}\cup
\{B(a)^{(r)}\underset{r}{\overset{2}{\longleftrightarrow}}B(n)\}\right]^c\]
is true. Then each of the sets $B(a)^{(r)}$ and $B(a)^{(2r)}\setminus B(a)^{(r)}$
contain exactly one of the points $x_1$ and $x_2$.

By arguments similar to the ones leading to (\ref{eqn:just before burton-keane part starts}),
we can show that in this case the event
\[E:=\big\{\exists x\in\poi\cap (B(a)^{(r)}\setminus B(a))\text{ such that }
S(x,r)\cup B(a)^{(r)}\underset{r}{\overset{2}{\longleftrightarrow}} B(n)\big\}\]
is true. For any realization $\eta=\{\eta_1,\hdots,\eta_{\ell}\}$ of $\poi\cap (B(a)^{(r)}\setminus B(a))$, we have
\[E\subset \cup_{j=1}^{\ell} \{S(\eta_j,r)\cup B(a)^{(r)}\underset{r}{\overset{2}{\longleftrightarrow}} B(n)\}.\]
Hence from Lemma \ref{lem:critical connectivity},
\begin{align*}
\PR(E)& \leq \frac{c_6\exp(c_7 a^{d-1})}{(\log n)^\frac{d}{2}}\E\bigg|\poi\cap (B(a)^{(r)}\setminus B(a))\bigg|\\
     &\leq \frac{c\exp(c_7 a^{d-1})a^{d-1}}{(\log n)^\frac{d}{2}}.\nonumber
\end{align*}
From our earlier discussion and another application of Lemma \ref{lem:critical connectivity},
\begin{align*}
\PR(B(a)\underset{r}{\overset{2}{\longrightarrow}}B(n))
&\leq \PR\left(B(a)\underset{r}{\overset{2}{\longleftrightarrow}}B(n)\right)+ \PR\left(B(a)^{(r)}\underset{r}{\overset{2}{\longleftrightarrow}}B(n)\right)+\PR(E)\nonumber\\
&\leq c_{15}\exp(c_{16} a^{d-1})/(\log n)^\frac{d}{2}.
\end{align*}
\item{${\mvr_2\leq \mvr\leq \mvn/\mveight}$ \textbf{and} ${\mvd\geq\mvthree}$:} The exponential decay in this regime can be proven using standard slab technology; see, e.g., the proof of Lemma 10.12 in \cite{penrosebook}. (Lemma 10.12 in \cite{penrosebook} is stated in the setup where each pair of Poisson points are connected if they are at distance at most one and the intensity of the Poisson process  determines sub- or super-criticality. This result translated to our setup where the parameter $r$ varies and the intensity of the Poisson process is kept fixed gives an upper bound for $\PR(B(a)\underset{r}{\overset{2}{\longrightarrow}}B(n))$ for every fixed $r>r_c$; whereas the bound in \eqref{eqn:->connectivity bound d=3} is uniform for $r_2\leq r\leq n/8$. This can be justified as follows. While using slab technology in the supercritical regime, $\PR(B(a)\underset{r}{\overset{2}{\longrightarrow}}B(n))$ is bounded by probability of an event which is decreasing in $r$. Thus the bound on $\PR(B(a)\underset{r_2}{\overset{2}{\longrightarrow}}B(n))$ obtained from the proof of Lemma 10.12 in \cite{penrosebook} works for each $r\in[r_2, n/8]$.) We omit the details.\hfill$\blacksquare$

\end{enumerateA}


\section{Rate of convergence in the CLT for Euclidean MST}\label{sec:proofs-clt}
Our goal in this section is to prove Theorem \ref{thm:poissonmst}.
As before, $d$ will denote the dimension of the ambient space.
Choose an integer $K$ such that $(n-1)/2\geq K\geq (n-2)/4$ and let
$s=n/(2K+1)$. Thus $s\in[1,2]$.
 Write $\R^d$ as the union of cubes,
\[\R^d=\cup_{j\in\Z^d} B_j\text{ where }B_j:=2sj+B(s).\]
Let $B(n)=\cup_{j\in\EL}B_j$. Clearly $\ell:=|\mathcal{L}|=\Theta(n^d)$. Fix $\alpha\in(0, 1)$ and let
\begin{align}\label{eqn:def-B-j-tilde}
\tilde B_j:=B(2sj,n^{\alpha}).
\end{align}
Further, define
\begin{align}\label{eqn:def-B-j-star}
B_j^\star:=
\left\{
\begin{array}{l}
B(2sj, a_n), \text{ if }d\geq 3,\\
B(2sj, \alpha\log n),\text{ if }d=2,
\end{array}
\right.\nonumber
\end{align}
where $a_n$ is a sequence increasing to infinity in a way so that $a_n\leq (\log\log n)^{1/(d-1/2)}$. (We will choose the sequence $a_n$ appropriately later in the proof.) We first prove a result that will be crucial in the proof.
\subsection{Preliminary estimates}\label{sec:preliminary-estimate}
Let $\poi$ be a Poisson process in $\bR^d$ having intensity one, and let $B_j, \tilde B_j$ and $B_j^\star$ be as above. Define the event $E_j$ as follows:
\begin{equation}\label{eqn:def-E}
E_j:=\big\{\poi\text{ contains a wall around }B_j\text{ in }B_j^\star\big\}.
\end{equation}
\begin{proposition}\label{prop:preliminary-estimate}
For any bounded subset $A$ of $\ \bR^d$, set $\cH(A)=M(\poi\cap A)$. Then the following hold.
\begin{enumeratei}
\item For every $j$ with $\|2sj\|_{\infty}\leq n-n^{\alpha}$,
\begin{align}\label{eqn:55}
&\E\bigg[\I_{E_j}\cdot\bigg|\big(\cH(B(n))-\cH(B(n)\setminus B_j)\big)-\big(\cH(\tilde B_j)-\cH(\tilde B_j\setminus B_j)\big)\bigg|\bigg]\nonumber\\
&\hskip60pt\leq
\left\{
\begin{array}{l}
c\ {\exp(c' a_n^{d-1})}{\left(\log n\right)^{-d/2}},\text{ if }d\geq 3,\\
c{(\log n)^3}{n^{-\alpha\beta}},\text{ if }d=2,
\end{array}
\right.
\end{align}
where $\beta$ is as in Lemma \ref{lem:->connectivity bound}.

\item Lower bound on variance:
\begin{equation}\label{eqn:lower bound on sigma^2}
\liminf_n\ \frac{1}{n^d}\ \E\big(\cH(B(n))-\E \cH(B(n))\big)^2>0.
\end{equation}
\end{enumeratei}
\end{proposition}
\noindent{\bf Proof of \eqref{eqn:55}:} We first deal with the case $d\geq 3$. Note that
\begin{align}\label{eqn:telescope0}
&\E\bigg[\I_{E_j}\cdot\bigg|\big(\cH(B(n))-\cH(B(n)\setminus B_j)\big)-\big(\cH(\tilde B_j)-\cH(\tilde B_j\setminus B_j)\big)\bigg|\bigg]\\
&\hskip20pt=\E\ \E_{\eta}
\bigg[\I_{E_j}\cdot\bigg|\big(\cH(B(n))-\cH(B(n)\setminus B_j)\big)-\big(\cH(\tilde B_j)-\cH(\tilde B_j\setminus B_j)\big)\bigg|\bigg],\nonumber
\end{align}
where $\E_{\eta}$ denotes expectation conditional on the event $\{\poi\cap B_j^\star=\eta\}$.

Fix realizations $\eta$, $\omega_1$ and $\omega_2$ of $\poi$
in $B_j^\star$, $\tilde B_j\setminus B_j^\star$ and $B(n)\setminus \tilde B_j$ respectively for which
the event ${E_j}$ is true.
If $|\eta\cap B_j|=0$, then $\cH(B(n))-\cH(B(n)\setminus B_j)$ and $\cH(\tilde B_j)-\cH(\tilde B_j\setminus B_j)$ are both zero.
So let us assume $|\eta\cap B_j|>0$, and write
\[\eta\cap B_j=\{v_1,\hdots,v_m\}\text{ and }\eta\cap\big(B_j^\star\setminus B_j\big)=\{p_1,\hdots,p_r\}.\]
Let $\mathfrak{J}_0=\emptyset$ and $\mathfrak{J}_i=\{v_1,\hdots,v_i\}$ for $1\leq i\leq m$.
Then
\begin{align}\label{eqn:telescopeI}
\big(\cH(B(n))-\cH(B(n)\setminus B_j)\big)-\big(\cH(\tilde B_j)-\cH(\tilde B_j\setminus B_j)\big)=\sum\nolimits_{i=1}^{m}\delta_i,
\end{align}
where
\begin{align*}
\delta_i:=&\big[M\big(\mathfrak{J}_i\cup(\poi\cap(B(n)\setminus B_j))\big)-M\big(\mathfrak{J}_{i-1}\cup(\poi\cap(B(n)\setminus B_j))\big)\big]\\
&\hskip20pt-\big[M\big(\mathfrak{J}_i\cup(\poi\cap(\tilde B_j\setminus B_j))\big)-M\big(\mathfrak{J}_{i-1}\cup(\poi\cap(\tilde B_j\setminus B_j))\big)\big].
\end{align*}
To keep the notation simple, let us focus on $\delta_1$.
Note that since $\eta$ contains a wall around $B_j$ in $B_j^\star$, by Lemma \ref{lem:wall property A}, an MST on the
complete graph on $\{v_1,p_1,\hdots,p_r\}\cup\{\omega_1\cup\omega_2\}$
(resp. $\{v_1,p_1,\hdots,p_r\}\cup\omega_1$)
cannot contain an edge of the form $\{v_1,p\}$ with $p\in \omega_1\cup\omega_2$ (resp. $p\in \omega_1$).
Thus, an MST on the
complete graph on $\{v_1,p_1,\hdots,p_r\}\cup\{\omega_1\cup\omega_2\}$
(resp. $\{v_1,p_1,\hdots,p_r\}\cup\omega_1$)
can be obtained from an MST on $\{p_1,\hdots,p_r\}\cup\{\omega_1\cup\omega_2\}$ (resp. $\{p_1,\hdots,p_r\}\cup\omega_1$)
by introducing the edges $\{v_1,p_j\}$ one by one and deleting the edge with maximum weight in the resulting cycle to make sure
all paths in the new tree are minimax, that is, by repeatedly using the add and delete algorithm (Section \ref{sec:add-and-delete-algo}).
We start with an MST $T_0$ (resp. $\tilde T_0$) on $\{p_1,\hdots,p_r\}\cup\{\omega_1\cup\omega_2\}$ (resp. $\{p_1,\hdots,p_r\}\cup\omega_1$)
with edge set $E$ (resp. $\tilde{E}$) and proceed in the following manner.

Set $E_0=E\ (\text{resp. }\tilde{E}_0=\tilde{E}),\ Y_0=d(v_1,p_1)\ (\text{resp. }\tilde{Y}_0=d(v_1,p_1))$ and let $w_0\ (\text{resp. }\tilde{w_0})$ be the weight of $T_0\ (\text{resp. }\tilde T_0)$. For $k=1,\hdots,r,$
\begin{enumeratei}
\item Introduce the edge $\{v_1,p_k\}$. If $k=1$, there will be no cycles in $E_{0}\cup\{v_1, p_1\}\ (\text{resp. }\tilde{E}_0\cup\{v_1,p_1\})$. In this case, set $E_1=E_0\cup\{v_1,p_1\}\ (\text{resp. }\tilde{E}_1=\tilde{E}_0\cup\{v_1,p_1\})$. Otherwise there will be a unique cycle in $E_{k-1}\cup\{v_1,p_k\}\ (\text{resp. }\tilde{E}_{k-1}\cup\{v_1,p_k\})$ having $\{v_1,p_k\}$ as one of its edges. Delete the edge in this cycle with maximum weight and set $E_k\ (\text{resp. }\tilde{E}_k)$ to be the resulting set of edges. If $k\leq r-1$,  let $Y_k\ (\text{resp. }\tilde Y_k)$ be the maximum edge weight in the path connecting $v_1$ and $p_{k+1}$ in the resulting tree, $T_k\ (\text{resp. }\tilde T_k)$ and let $w_k\ (\text{resp. }\tilde w_k)$ be the total weight of $T_k\ (\text{resp. }\tilde T_k)$.
\item If $k=r$, stop. Otherwise increase $k$ by one and repeat step {\upshape(i)}.
\end{enumeratei}
A consequence of Proposition \ref{prop:add-and-delete-algo} is that the
tree we get at the end of this process is an MST on the graph which has
$\{v_1,p_1,\hdots,p_r\}\cup\{\omega_1\cup\omega_2\}$
(resp. $\{v_1,p_1,\hdots,p_r\}\cup\omega_1$) as its vertex set and contains every possible edge between these vertices
except the ones of the form $\{v_1,p\}$ with $p\in \omega_1\cup\omega_2$ (resp. $p\in \omega_1$).
It is easy to see that the resulting tree is actually an MST
on the complete graph on $\{v_1,p_1,\hdots,p_r\}\cup\{\omega_1\cup\omega_2\}$
(resp. $\{v_1,p_1,\hdots,p_r\}\cup\omega_1$), because as argued before, an edge
of the form $\{v_1, x\}$ with $x\notin B_j^\star$ cannot be present in an an MST
since $\eta$ contains a wall around $B_j$ in $B_j^\star$.

Hence
\begin{align}\label{eqn:delta_1 telescope}
\delta_1&=(w_r-w_0)-(\tilde w_r-\tilde w_0)
=\sum_{k=1}^r\big[(w_k-w_{k-1})-(\tilde w_{k}-\tilde w_{k-1})\big].
\end{align}
Now,
\begin{align}\label{eqn:delta_1 telescope I}
w_k-w_{k-1}=
\left\{
\begin{array}{l}
d(v_1,p_1), \text{ if }k=1,\\
d(v_1,p_k)-\max(Y_{k-1},d(v_1,p_k)),\text{ if }2\leq k\leq r.
\end{array}\right.
\end{align}
A similar statement holds for $\tilde w_k$ with $\tilde Y_{k-1}$ replacing $Y_{k-1}$.
Proposition \ref{prop:add-and-delete-algo} shows that $T_{k-1}$ (resp. $\tilde T_{k-1}$) is an MST
on the graph with vertex set $\mathcal{V}=(\poi\cap(B(n)\setminus B_j))\cup\{v_1\}$
(resp. $\tilde{\mathcal{V}}=(\poi\cap(\tilde B_j\setminus B_j))\cup\{v_1\}$)
and edge set $\mathcal{E}_{k-1}=\cup_{i=1}^{k-1}\{v_1,p_i\}\cup\{\text{edges in the}$ $\text{ complete graph on }\poi\cap(B(n)\setminus B_j)\}$
(resp. $\tilde{\mathcal{E}}_{k-1}=\cup_{i=1}^{k-1}\{v_1,p_i\}\cup\{$ edges in
the complete graph on $\poi\cap(\tilde B_j\setminus B_j)\}$) for $k\geq 2$. Hence $Y_{k-1}$ (resp. $\tilde Y_{k-1}$)
is the maximum edge-weight in a minimax path connecting
$v_1$ and $p_k$ in $(\mathcal{V},\mathcal{E}_{k-1})$ (resp. $(\tilde{\mathcal{V}},\tilde{\mathcal{E}}_{k-1})$).
This gives $Y_{k-1}\leq \tilde Y_{k-1}$. From (\ref{eqn:delta_1 telescope I}),
\begin{equation}\label{eqn:delta_1 telescope II}
0\leq (w_k-w_{k-1})-(\tilde w_{k}-\tilde w_{k-1})\leq \tilde Y_{k-1}-Y_{k-1}.
\end{equation}
Consider a random variable $U$ uniformly distributed on $(0,2\sqrt{d}a_n)$ which is
independent of $\poi$. We have
\begin{align}\label{eqn:introducing auxilaiary U}
&\E_{\eta}|(w_k-w_{k-1})-(\tilde w_{k}-\tilde w_{k-1})|\\
&\hskip15pt\leq\E_{\eta}(\tilde Y_{k-1}-Y_{k-1})
=2\sqrt{d}a_n\cdot\PR_{\eta}\big(Y_{k-1}<U<\tilde Y_{k-1}\big)\nonumber\\
&\hskip30pt=\int_{0}^{2\sqrt{d}a_n}\PR_{\eta}\big(Y_{k-1}<u<\tilde Y_{k-1}\big)\ du
\leq \int_{0}^{2\sqrt{d}a_n}\PR\big(B_j^\star\underset{u/2}{\overset{2}{\longrightarrow}}\tilde B_j\big)\ du.\nonumber
\end{align}

The last inequality holds because of the following reason. Assume that $Y_{k-1}<u<\tilde Y_{k-1}$ and let
$(v_1=z_0,z_1,\hdots,z_{\ell}=p_k)$ be a minimax path connecting $v_1$ and $p_k$ in $(\mathcal{V},\mathcal{E}_{k-1})$.
 Since $Y_{k-1}<\tilde Y_{k-1}$,
$z_i\in \tilde B_j^c$ for some $i\leq\ell$. Let $k_1+1:=\min\{i\leq\ell:z_i\in\tilde B_j^c\}$
and $k_2-1:=\max\{i\leq\ell:z_i\in\tilde B_j^c\}$. Then the $u/2$-clusters
in $\tilde{B}_j\setminus B_j$ containing $\{z_1,\hdots,z_{k_1}\}$
and $\{z_{k_2},\hdots,z_{\ell}\}$ are disjoint, since otherwise we could find a path
$(z_i=y_0,y_1,\hdots,y_t=z_{i'})$ for some $i\leq k_1$, $i'\geq k_2$
such that $y_p\in\tilde{\mathcal{V}}\setminus \{v_1\}$ and $d(y_p,y_{p+1})\leq u$
for every $p\leq t-1$. But this would mean that
$(z_0,\hdots,z_i,y_1,\hdots,y_{t-1},z_i',\hdots,z_{\ell})$ is a path in $(\tilde{\mathcal{V}},\tilde{\mathcal{E}}_{k-1})$
connecting $v_1$ and $p_k$ with maximum edge-weight strictly smaller than $\tilde Y_{k-1}$, a contradiction.
Then the restrictions of the (disjoint) $u/2$-clusters
in $\tilde{B}_j\setminus B_j$ containing $\{z_1,\hdots,z_{k_1}\}$
and $\{z_{k_2},\hdots,z_{\ell}\}$ to $\tilde B_j\setminus B_j^\star$ will contain two disjoint $u/2$-clusters
which will satisfy the criteria for $\{B_j^\star\underset{u/2}{\overset{2}{\longrightarrow}}\tilde B_j\}$
to hold.

Combining \eqref{eqn:delta_1 telescope} and \eqref{eqn:introducing auxilaiary U},
\begin{align}\label{eqn:bound-delta-1}
\E_{\eta}\big[\delta_1\big]
\leq
2\sqrt{d}a_n
\sup_{0<u<2\sqrt{d}a_n} \PR\bigg(B_j^\star\underset{u/2}{\overset{2}{\longrightarrow}}\tilde B_j\bigg)\cdot
\bigg(|\poi\cap B_j^\star|\bigg).
\end{align}

Inductively, having obtained an MST on $\fJ_i\cup(\poi\cap(B(n)\setminus B_j))$
(resp. $\fJ_i\cup(\poi\cap(\tilde B_j\setminus B_j))$), $1\leq i\leq m-1$, an MST on $\fJ_{i+1}\cup(\poi\cap(B(n)\setminus B_j))$
(resp. $\fJ_{i+1}\cup(\poi\cap(\tilde B_j\setminus B_j))$) can be obtained
by introducing the edges $\{v_{i+1},p_j\}$, $1\leq j\leq r$, and $\{v_{i+1},v_s\}$, $1\leq s\leq i$, one by one and again using the {\it add and delete algorithm}. Thus $\delta_{i+1}$ will have a decomposition similar to \eqref{eqn:delta_1 telescope} that has $r+i\leq |\poi\cap B_j^\star|$ terms, and each of these terms will obey the bound on the right side of \eqref{eqn:introducing auxilaiary U}.
Hence, for each $i\leq m$, \eqref{eqn:bound-delta-1} will continue to hold for $\E_{\eta}\big[\delta_i\big]$.

Combining this observation with \eqref{eqn:telescope0} and \eqref{eqn:telescopeI},
we get
\begin{align}\label{eqn:E-I_E}
&\E\bigg[\I_{E_j}\cdot\bigg|\big(\cH(B(n))-\cH(B(n)\setminus B_j)\big)-\big(\cH(\tilde B_j)-\cH(\tilde B_j\setminus B_j)\big)\bigg|\bigg]\\
&\hskip15pt\leq 2\sqrt{d}a_n
\sup_{0<u<2\sqrt{d}a_n} \PR\bigg(B_j^\star\underset{u/2}{\overset{2}{\longrightarrow}}\tilde B_j\bigg)\cdot
\E\bigg(|\poi\cap B_j|\cdot |\poi\cap B_j^\star|\bigg)
\nonumber\\
&\hskip15pt\leq c a_n^{d+1}\sup_{0<u<2\sqrt{d}a_n}
\PR\bigg(B_j^\star\underset{u/2}{\overset{2}{\longrightarrow}}\tilde B_j\bigg)
\leq c \frac{\exp(c' a_n^{d-1})}{\left(\log n\right)^{d/2}},\nonumber
\end{align}
where the last step follows from Lemma \ref{lem:->connectivity bound}. This completes the proof for the case $d\geq 3$.

When $d=2$, we can proceed in the exact same manner and the only difference is the percolation estimate from Lemma \ref{lem:->connectivity bound}. Thus when $d=2$,
\begin{align}\label{eqn:E-I_E-d=2}
&\E\bigg[\I_{E_j}\cdot\bigg|\big(\cH(B(n))-\cH(B(n)\setminus B_j)\big)-\big(\cH(\tilde B_j)-\cH(\tilde B_j\setminus B_j)\big)\bigg|\bigg]\\
&\hskip15pt\leq 2\sqrt{2}\cdot\alpha\log n
\sup_{0<u<2\sqrt{2}\cdot\alpha\log n} \PR\bigg(B_j^\star\underset{u/2}{\overset{2}{\longrightarrow}}\tilde B_j\bigg)\cdot
\E\bigg(|\poi\cap B_j|\cdot|\poi\cap B_j^\star|\bigg)
\nonumber\\
&\hskip15pt\leq c (\log n)^3\sup_{0<u<2\sqrt{2}\cdot\alpha\log n}
\PR\bigg(B_j^\star\underset{u/2}{\overset{2}{\longrightarrow}}\tilde B_j\bigg)
\leq c' \frac{(\log n)^3}{n^{\alpha\beta}}.\nonumber
\end{align}
This completes the proof of \eqref{eqn:55}.\hfill$\blacksquare$

\medskip

\noindent{\bf Proof of \eqref{eqn:lower bound on sigma^2}:}
This is implicit in the work of Kesten and Lee in \cite{kesten}. Let us write $\mathcal{L}=\{j_1,\hdots,j_l\}$
(recall the definition of $\cL$ from around \eqref{eqn:def-B-j-tilde}).
Define the sigma-fields $\mathcal{F}_k:=\sigma\{\poi\cap B_{j_i}:i\leq k\}$ for $k=1,\hdots,\ell$,
and let $\mathcal{F}_0$ be the trivial sigma-field. Then we can express $\cH(B(n))-\E\cH(B(n))$
as a sum of martingale differences:
\[\cH(B(n))-\E\cH(B(n))=\sum_{k=1}^{\ell}Z_k,\text{ where }
Z_k:=\E\big(\cH(B(n))\mid\mathcal{F}_k\big)-\E\big(\cH(B(n))\mid\mathcal{F}_{k-1}\big).\]
From \cite[Equation (4.27)]{kesten}, it will follow that
\[\frac{1}{\ell}\sum_{k=1}^{\ell}Z_k^2\overset{P}{\rightarrow}\zeta,\]
for a positive constant $\zeta$. An application of Fatou's lemma together with fact $\ell=\Theta(n^d)$ yields
\begin{equation*}
\liminf_n\ \frac{1}{n^d}\E\big(\cH(B(n))-\E\cH(B(n))\big)^2>0,
\end{equation*}
as desired.\hfill$\blacksquare$

\subsection{Proof of Theorem \ref{thm:poissonmst}}\label{sec:actual proof of poissonmst}
At this point we ask the reader to recall the notation used in Section \ref{sec:stein's-method}.
Consider two independent Poisson process $\poi$ and $\poi'$ having intensity
one in $\R^d$. We will apply \eqref{eqn:chatterjee} and \eqref{eqn:lrp}
with
\[X_j:=\poi\cap B_j,\quad X_j':=\poi'\cap B_j,\quad X:=(X_j : j\in\EL),\quad X':=(X_j' : j\in\EL),\]
and the function $f:\prod_{i\in\EL}\mathfrak{X}(B_i)\to\R$ given by
\[f\big(\{\omega_i: i\in\EL\}\big)=M\bigg(\displaystyle{\bigcup_{i\in\EL}}\omega_i\bigg).\]

By definition, for any $A\subset\cL$, $X^A$ is a random vector whose $i$-th coordinate is a configuration in $B_i$, $i\in\cL$, but
there is also a natural way of identifying $X^A$ with a configuration in $B(n)$, and we will often blur the distinction between the two to simplify notation.
In particular, with this convention, $X\cap R$ will represent a configuration in $R$ for any $R\subset\bR^d$, and $M(X)$ will be synonymous with $M(\cup_{i\in\cL} X_i)$.
We will use the shorthand $\Delta_j f(X^A):=\Delta_j f(X^A,X')$. Thus
\[\Delta_j f(X^A):=f(X^A)-f(X^{A\cup\{j\}}),\]
for $A\subset\cL$.

We first focus on proving the bounds on the Kantorovich-Wasserstein distance. Bounds of the same order in the Kolmogorov distance can be obtained in an almost identical fashion, and we will briefly comment on this at the end.

\medskip

\noindent{\bf Bounds on the Kantorovich-Wasserstein distance.}
We will use Theorem \ref{thm:chatterjee} to  prove bounds on the Kantorovich-Wasserstein distance.
Note that $X\cap(B(n)\setminus B_j)=X^j\cap(B(n)\setminus B_j)$, and hence
\begin{align*}
\Delta_j f(X)&=\bigg[M(X)-M\big(X\cap(B(n)\setminus B_j)\big)\bigg]
-\bigg[M(X^j)-M\big(X^j\cap(B(n)\setminus B_j)\big)\bigg]
\end{align*}
for every $j\in\EL$. Lemma \ref{lem:moments are bounded} and the fact $s\in[1,2]$ imply that for every $j\in\EL$ and $q\geq 1$,
\begin{equation}\label{eqn:moments of Delta_j  are bounded}
\E|\Delta_j f(X)|^q\leq C_q',
\end{equation}
for constants $C_q'$ depending only on $d$ and $q$.
Here we make note of two direct consequences of \eqref{eqn:moments of Delta_j  are bounded}. First,
\begin{equation}\label{eqn:cov are bounded}
\bigg|\mathrm{Cov}\big(\Delta_j f(X)\Delta_j f(X^A),\ \Delta_{j'} f(X)\Delta_{j'} f(X^{A'})\big)\bigg|
\leq C_{\ref{eqn:cov are bounded}}
\end{equation}
for any $j,j'\in\cL$ and $A, A'\subset\cL$,
where $C_{\ref{eqn:cov are bounded}}$ is a finite constant.
Secondly, \eqref{eqn:moments of Delta_j  are bounded} combined with \eqref{eqn:lower bound on sigma^2}
and the fact $\ell=|\cL|=\Theta(n^d)$ yields
\begin{equation}\label{eqn:second term upper bound}
\frac{1}{\mathrm{Var}(f(X))^{3/2}}\sum_{j=1}^{\ell}\E|\Delta_j f(X)|^3\leq \frac{c}{n^{d/2}}.
\end{equation}
This gives us control over the second term on the right side of \eqref{eqn:chatterjee}. Our aim in the remainder of the proof is
to bound $\var(\E(T|W))$.
We first focus on the case $d\geq 3$.

\medskip
\noindent{\bf Proof of \eqref{eqn:theorem 1 d>2}: }
We plan to show that the covariance term appearing in the numerator on
the right side of (\ref{eqn:var expansion term bounded}) is small when
$j$ and $j'$ are ``far away.'' With this in mind, we break up the sum
on the right side of (\ref{eqn:var expansion term bounded}) into two parts
$\sum_1$ and $\sum_2$; $\sum_1$ denotes the sum over all $(j, j', A, A')\in \mathfrak{E}^{(\alpha)}$
(for some $\alpha\in(0,1)$), where
\begin{align}\label{eqn:def-E-alpha}
\mathfrak{E}^{(\alpha)} := \big\{&(j, j', A, A'):\
A,A'\subsetneq\EL;\ j\in\mathcal{L}\setminus A,\ j'\in\mathcal{L}\setminus A'\text{ and either}\\
&\hskip20pt\|j-j'\|_{\infty}\leq n^{\alpha}\text{ or } \|2sj\|_{\infty}>(n-n^{\alpha})\text{ or }
\|2sj'\|_{\infty}>(n-n^{\alpha})\big\}\nonumber
\end{align}
and $\sum_2$ denotes the sum over the remaining terms, i.e.,
all $(j,j',A,A')\in\mathfrak{F}^{(\alpha)}$ where
\[\mathfrak{F}^{(\alpha)}:=\big\{(j, j', A, A'):\
A,A'\subsetneq\EL;\ j\in\mathcal{L}\setminus A,\ j'\in\mathcal{L}\setminus A'\big\}\setminus\mathfrak{E}^{(\alpha)}.\]


Let $\mathfrak{E}_{1,2}^{(\alpha)}$ be the collection of all $(j,j')$
for which $(j, j', \emptyset, \emptyset)\in\mathfrak{E}^{(\alpha)}$.
Then from (\ref{eqn:cov are bounded}),
\begin{align}\label{eqn:sigma_1 is small}
&\sum\nolimits_1
\frac{\mathrm{Cov}(\Delta_j f(X)\Delta_j f(X^A),\Delta_{j'} f(X)\Delta_{j'} f(X^{A'}))}
{\dbinom{\ell}{|A|}(\ell-|A|)\dbinom{\ell}{|A'|}(\ell-|A'|)}\\
&\leq C_{\ref{eqn:cov are bounded}}\sum_{(j,j')\in\mathfrak{E}_{1,2}^{(\alpha)}}\sum_{\substack{A\not\ni j\\ A'\not\ni j'}}
\left(\dbinom{\ell}{|A|}(\ell-|A|)\dbinom{\ell}{|A'|}(\ell-|A'|)\right)^{-1}\nonumber\\
&= C_{\ref{eqn:cov are bounded}}\sum_{(j,j')\in\mathfrak{E}_{1,2}^{(\alpha)}}\sum_{k,k'=0}^{\ell-1}
\sum_{\substack{A\not\ni j, A'\not\ni j'\\|A|=k,|A'|=k'}}
\left(\dbinom{\ell}{|A|}(\ell-|A|)\dbinom{\ell}{|A'|}(\ell-|A'|)\right)^{-1}\nonumber\\
&=C_{\ref{eqn:cov are bounded}} |\mathfrak{E}_{1,2}^{(\alpha)}|\leq  c (n^{2d-1}\cdot n^{\alpha}+n^d\cdot n^{\alpha d})
\leq c' n^{2d-1+\alpha}.\nonumber
\end{align}

We now turn to the sum $\sum_2$. Note that for $(j,j')\notin\mathfrak{E}_{1,2}^{(\alpha)}$, $\|2sj-2sj'\|_{\infty}>2s n^{\alpha}\geq 2n^{\alpha}$, and so
the cubes $\tilde B_j$ and $\tilde B_{j'}$ are disjoint (recall the definition from \eqref{eqn:def-B-j-tilde}). As a result, the restrictions
of $X$ (and of $X'$) to these cubes are independent.
Let us now define
\begin{equation}\label{eqn:numberit}
\tilde\Delta_j f(X^A):=M\big(X^A\cap\tilde B_j\big)-M\big(X^{A\cup\{j\}}\cap\tilde B_j\big)
\end{equation}
for every $j$ with $\|2sj\|_{\infty}\leq (n-n^{\alpha})$ and $A\subset\cL$.
Whenever $(j,j',A,A')\in\mathfrak{F}^{(\alpha)}$, we have
\begin{align}\label{eqn:decomposition of covariance terms}
&\mathrm{Cov}\left(\Delta_j f(X)\Delta_j f(X^A),\ \Delta_{j'} f(X)\Delta_{j'} f(X^{A'})\right)\\
&\hskip20pt=\mathrm{Cov}\left(\big[\Delta_j f(X)-\tilde\Delta_j f(X)\big]\Delta_j f(X^A),\
\Delta_{j'} f(X)\Delta_{j'} f(X^{A'})\right)\nonumber\\
&\hskip25pt+\mathrm{Cov}\left(\tilde\Delta_j f(X)\big[\Delta_j f(X^A)-\tilde\Delta_j f(X^A)\big],\
\Delta_{j'} f(X)\Delta_{j'} f(X^{A'})\right)\nonumber\\
&\hskip25pt+\mathrm{Cov}\left(\tilde\Delta_j f(X)\tilde\Delta_j f(X^A),\
\big[\Delta_{j'} f(X)-\tilde\Delta_{j'} f(X)\big]\Delta_{j'} f(X^{A'})\right)\nonumber\\
&\hskip25pt+\mathrm{Cov}\left(\tilde\Delta_j f(X)\tilde\Delta_j f(X^A),\
\tilde\Delta_{j'} f(X)\big[\Delta_{j'} f(X^{A'})-\tilde\Delta_{j'} f(X^{A'})\big]\right).\nonumber
\end{align}
We will give an upper bound for the first term on the right side
of \eqref{eqn:decomposition of covariance terms}. The other terms can be dealt with in
a similar fashion. Note that
\begin{align}\label{eqn:bound on covariance}
&\mathrm{Cov}\bigg(\big[\Delta_j f(X)-\tilde\Delta_j f(X)\big]\Delta_j f(X^A),\ \Delta_{j'} f(X)\Delta_{j'} f(X^{A'})\bigg)\\
&\hskip20pt\leq \E\bigg(\bigg|\big(\Delta_j f(X)-\tilde\Delta_j f(X)\big)\Delta_j f(X^A)\Delta_{j'} f(X)\Delta_{j'} f(X^{A'})\bigg|\bigg)\nonumber\\
&\hskip35pt+\E\bigg(\bigg|\big(\Delta_j f(X)-\tilde\Delta_j f(X)\big)\Delta_j f(X^A)\bigg|\bigg)\cdot
\E\bigg(\bigg|\Delta_{j'} f(X)\Delta_{j'} f(X^{A'})\bigg|\bigg)\nonumber\\
&\hskip35pt=:T_1+T_2.\nonumber
\end{align}
Then for any $p,\ q>1$ satisfying $p^{-1}+q^{-1}=1$, we have from (\ref{eqn:moments of Delta_j  are bounded}) that
\begin{align}
T_2\leq & C_2'\left(\E|(\Delta_j f(X)-\tilde\Delta_j f(X))|\right)^{1/p}
\cdot \left(\E|(\Delta_j f(X)-\tilde\Delta_j f(X))|\Delta_j f(X^A)|^q|\right)^{1/q}\nonumber\\
&\hskip20pt\leq c \left(\E|(\Delta_j f(X)-\tilde\Delta_j f(X))|\right)^{1/p}.\nonumber
\end{align}
A similar bound holds for $T_1$. We plug all these estimates into (\ref{eqn:decomposition of covariance terms})
to get
\begin{align}\label{eqn:bound on covariance I}
&\mathrm{Cov}\left(\Delta_j f(X)\Delta_j f(X^A),\ \Delta_{j'} f(X)\Delta_{j'} f(X^{A'})\right)\\
&\hskip20pt\leq c\left(\left(\E|(\Delta_j f(X)-\tilde\Delta_j f(X))|\right)^{\frac{1}{p}}
+\left(\E|(\Delta_{j'} f(X)-\tilde\Delta_{j'} f(X))|\right)^{\frac{1}{p}}\right),\nonumber
\end{align}
for $(j,j',A,A')\in\mathfrak{F}^{(\alpha)}$ .
Let $E_j$ be the event in \eqref{eqn:def-E}.
Then \eqref{eqn:moments of Delta_j  are bounded} and Lemma \ref{lem:moments are bounded} yield
\begin{align}\label{eqn:E-I-E^c}
\E\left(\I_{E_j^c}|(\Delta_{j} f(X)-\tilde\Delta_{j} f(X))|\right)
&\leq \left(\E|(\Delta_{j} f(X)-\tilde\Delta_{j} f(X))|^2\right)^{\frac{1}{2}}\PR(E_j^c)^{\frac{1}{2}}\nonumber\\
&\leq c\exp(-c_{\ref{eqn:E-I-E^c}}a_n^d)
\end{align}
for every $j$ with $\|2sj\|\leq n-n^{\alpha}$. Hence, for every $j$ with $\|2sj\|\leq n-n^{\alpha}$,
\begin{align}\label{eqn:wall exists + no wall exists}
\E\big|(\Delta_{j} f(X)-\tilde\Delta_{j} f(X))\big|&\leq  c \exp(-c_{\ref{eqn:E-I-E^c}}a_n^d)
+\E\left(\I_{E_j}\cdot|(\Delta_{j} f(X)-\tilde\Delta_{j} f(X))|\right).
\end{align}
Since the restrictions of the vectors $X$ and $X^j$ to $B(n)\setminus B_j$ (resp. $\tilde B_j\setminus B_j$) are the same,
\[M\big(X\cap (B(n)\setminus B_j)\big)=M\big(X^j\cap(B(n)\setminus B_j)\big),\
\text{ and }\
M\big(X\cap (\tilde B_j\setminus B_j)\big)=M\big(X^j\cap(\tilde B_j\setminus B_j)\big).\]
Hence we can write, for every $j$ with $\|2sj\|\leq n-n^{\alpha}$,
\begin{align}
\Delta_{j} f(X)-\tilde\Delta_{j} f(X)
&=\bigg[\big(M(X)-M\big(X\cap (B(n)\setminus B_j)\big)\big)-\big(M(X\cap\tilde B_j)
-M\big(X\cap (\tilde B_j\setminus B_j)\big)\big)\bigg]\nonumber\\
&\hskip10pt-\bigg[\big(M(X^j)-M\big(X^j\cap(B(n)\setminus B_j)\big)\big)
-\big(M(X^j\cap\tilde B_j)-M(X^j\cap(\tilde B_j\setminus B_j))\big)\bigg].\nonumber
\end{align}
Therefore, for every $j$ with $\|2sj\|\leq n-n^{\alpha}$,
\begin{align}\label{eqn:telescope0-'}
&\E\bigg[\I_{E_j}|\Delta_{j} f(X)-\tilde\Delta_{j} f(X)|\bigg]\\
&\leq 2\E \bigg[\I_{E_j}\cdot\bigg|\big(M(X)-M\big(X\cap (B(n)\setminus B_j)\big)\big)-\big(M(X\cap\tilde B_j)
-M\big(X\cap (\tilde B_j\setminus B_j)\big)\big)\bigg|\bigg].\nonumber
\end{align}
Using \eqref{eqn:55}, we conclude that
\begin{equation}
\E\left[\I_{E_j}|\Delta_{j} f(X)-\tilde\Delta_{j} f(X)|\right]
\leq
c \frac{\exp(c' a_n^{d-1})}{\left(\log n\right)^{d/2}}.
\end{equation}
In view of (\ref{eqn:wall exists + no wall exists}), we choose $a_n$ so that
$c_{\ref{eqn:E-I-E^c}} a_n^d=\frac{d}{2}\log\log n$ to get
\begin{align}\label{eqn:numerator boundI}
\E|\Delta_{j} f(X)-\tilde\Delta_{j} f(X)|
&\leq
c \frac{\exp(c'' (\log\log n)^{\frac{d-1}{d}})}{\left(\log n\right)^{d/2}}
\end{align}
for every $j$ with $\|2sj\|\leq n-n^{\alpha}$. Hence
\begin{align}\label{eqn:numerator boundII}
\sum\nolimits_2
\frac{\mathrm{Cov}\big(\Delta_j f(X)\Delta_j f(X^A),\Delta_{j'} f(X)\Delta_{j'} f(X^{A'})\big)}
{\dbinom{\ell}{|A|}(\ell-|A|)\dbinom{\ell}{|A'|}(\ell-|A'|)}\hskip60pt &\\
\leq c n^{2d}\max_{\mathfrak{F}^{(\alpha)}}
\ \mathrm{Cov}\big(\Delta_j f(X)\Delta_j f(X^A),\Delta_{j'} f(X)\Delta_{j'} f(X^{A'})\big)&\nonumber\\
\leq c n^{2d}
\frac{\exp(c''\cdot (\log\log n)^{\frac{d-1}{d}}/p)}{\left(\log n\right)^{\frac{d}{2p}}},\hskip70pt &\nonumber
\end{align}
where the last inequality is a consequence of \eqref{eqn:bound on covariance I}
and \eqref{eqn:numerator boundI}. Combining \eqref{eqn:var expansion term bounded},
\eqref{eqn:sigma_1 is small}, and \eqref{eqn:numerator boundII} and observing that \eqref{eqn:numerator boundII} is true for
any $p>1$, we get
\begin{align}\label{eqn:66}
\var\big(\E(T|W)\big)\leq
cn^{2d}\left(\log n\right)^{-\frac{d}{2p}}.
\end{align}

Combining
\eqref{eqn:chatterjee},
\eqref{eqn:lower bound on sigma^2},
\eqref{eqn:66}, and
\eqref{eqn:second term upper bound}, we see that there exists a positive
constant $c$ depending on $p$ and $d$ such that
\begin{align}\label{eqn:bound with log in the expression}
\mathcal{W}(\mu_n,\gamma)\leq
c\left(\log n\right)^{-\frac{d}{4p}},
\end{align}
which is the bound claimed in \eqref{eqn:theorem 1 d>2}.

Let us now turn to the case $d=2$.

\medskip

\noindent{\bf Proof of \eqref{eqn:theorem 1 d=2}:} Let ${\mathfrak{E}}^{(\alpha)}, {\mathfrak{F}}^{(\alpha)}, \sum_1$ and $\sum_2$ be as defined around $\eqref{eqn:def-E-alpha}$. (Later we will make a suitable choice of $\alpha$.)
The calculation in (\ref{eqn:sigma_1 is small}) gives
\begin{equation}\label{eqn:sigma_1 is small in d=2}
\sum\nolimits_1
\frac{\mathrm{Cov}(\Delta_j f(X)\Delta_j f(X^A),\Delta_{j'} f(X)\Delta_{j'} f(X^{A'}))}
{\dbinom{\ell}{|A|}(\ell-|A|)\dbinom{\ell}{|A'|}(\ell-|A'|)}\leq
c n^{3+\alpha}.
\end{equation}

We now bound the sum $\sum_2$. First
recall the definition of $\tilde B_j$ from \eqref{eqn:def-B-j-tilde} and let ${E_j}$ be the event in \eqref{eqn:def-E}.
With $\tilde\Delta_j f(X)$ as in \eqref{eqn:numberit} (defined for $j$ with $\|2sj\|_{\infty}\leq(n-n^{\alpha})$),
\eqref{eqn:decomposition of covariance terms}, \eqref{eqn:bound on covariance} and \eqref{eqn:bound on covariance I}
continue to hold.
Further, the bound \eqref{eqn:E-I-E^c} now reads
\begin{align}\label{eqn:E-I-E^c-d=2}
\E\left(\I_{{E_j}^c}|(\Delta_{j} f(X)-\tilde\Delta_{j} f(X))|\right)
\leq c \exp(-c'(\log n)^2)
\end{align}
for every $j$ with $\|2sj\|_{\infty}\leq(n-n^{\alpha})$,
and \eqref{eqn:55} combined with \eqref{eqn:telescope0-'} gives
\begin{align}\label{eqn:E-I-E-d=2}
\E\left[\I_{E_j}|\Delta_{j} f(X)-\tilde\Delta_{j} f(X)|\right]
\leq c (\log n)^{3}n^{-\alpha\beta}.
\end{align}
Combining \eqref{eqn:E-I-E^c-d=2} and \eqref{eqn:E-I-E-d=2}, we arrive at
\begin{align}
\E\bigg|\Delta_{j} f(X)-\tilde\Delta_{j} f(X)\bigg|
\leq c (\log n)^{3}n^{-\alpha\beta}\nonumber
\end{align}
for every $j$ with $\|2sj\|_{\infty}\leq(n-n^{\alpha})$.
Arguments similar to the ones used previously for $d\geq 3$ (see \eqref{eqn:numerator boundII} and \eqref{eqn:bound on covariance I}) now yield
\begin{equation}\label{eqn:sigma_2 is small in d=2}
\sum\nolimits_2
\frac{\mathrm{Cov}(\Delta_j f(X)\Delta_j f(X^A),\Delta_{j'} f(X)\Delta_{j'} f(X^{A'}))}
{\dbinom{\ell}{|A|}(\ell-|A|)\dbinom{\ell}{|A'|}(\ell-|A'|)}\leq
cn^4/n^{\frac{\alpha\beta}{p}}.
\end{equation}
Combining \eqref{eqn:sigma_1 is small in d=2} and \eqref{eqn:sigma_2 is small in d=2} and taking $\alpha=p/(\beta+p)$, we get
\begin{equation}\label{eqn:final step for thm 1}
\mathrm{Var}(\E(T|W))\leq cn^{{\frac{4p+3\beta}{\beta+p}}}.
\end{equation}

Combining \eqref{eqn:lower bound on sigma^2},
\eqref{eqn:second term upper bound} (with $d=2$),
\eqref{eqn:final step for thm 1}, and
\eqref{eqn:chatterjee}, and noting that $(4p+3\beta)/(\beta+p)<4$, we get the bound in \eqref{eqn:theorem 1 d=2} in the Kantorovich-Wasserstein distance.

\medskip

\noindent{\bf Bounds on the Kolmogorov distance.}
We can use Theorem \ref{thm:lrp} to prove bounds on the Kolmogorov distance in an almost identical fashion.
The difference in the bound in \eqref{eqn:lrp} comes from the terms $T_A'$, which are sums of terms of the form
$\Delta_j f(X, X')|\Delta_j f(X^A, X')|$ (instead of $\Delta_j f(X, X')\Delta_j f(X^A, X')$ as in Theorem \ref{thm:chatterjee}).
To take this into account, we modify \eqref{eqn:decomposition of covariance terms} as follows:
\begin{align*}\label{eqn:decomposition of covariance terms-modified}
&\mathrm{Cov}\left(\Delta_j f(X)\left|\Delta_j f(X^A)\right|,\ \Delta_{j'} f(X)\left|\Delta_{j'} f(X^{A'})\right|\right)\\
&\hskip20pt=\mathrm{Cov}\left(\big(\Delta_j f(X)-\tilde\Delta_j f(X)\big)\left|\Delta_j f(X^A)\right|,\
\Delta_{j'} f(X)\left|\Delta_{j'} f(X^{A'})\right|\right)\nonumber\\
&\hskip25pt+\mathrm{Cov}\left(\tilde\Delta_j f(X)\big(\left|\Delta_j f(X^A)\right|-\left|\tilde\Delta_j f(X^A)\right|\big),\
\Delta_{j'} f(X)\left|\Delta_{j'} f(X^{A'})\right|\right)\nonumber\\
&\hskip25pt+\mathrm{Cov}\left(\tilde\Delta_j f(X)\left|\tilde\Delta_j f(X^A)\right|,\
\big(\Delta_{j'} f(X)-\tilde\Delta_{j'} f(X)\big)\left|\Delta_{j'} f(X^{A'})\right|\right)\nonumber\\
&\hskip25pt+\mathrm{Cov}\left(\tilde\Delta_j f(X)\left|\tilde\Delta_j f(X^A)\right|,\
\tilde\Delta_{j'} f(X)\big(\left|\Delta_{j'} f(X^{A'})\right|-\left|\tilde\Delta_{j'} f(X^{A'})\right|\big)\right),\nonumber
\end{align*}
whenever $(j,j',A,A')\in\mathfrak{F}^{(\alpha)}$. Noting that
\[\left|\left|\Delta_j f(X^A)\right|-\left|\tilde\Delta_j f(X^A)\right|\right|
\leq
\left|\Delta_j f(X^A)-\tilde\Delta_j f(X^A)\right|,
\]
it is easy to see that a bound similar to \eqref{eqn:bound on covariance I} continues to hold:
\begin{align*}
&\mathrm{Cov}\left(\Delta_j f(X)\left|\Delta_j f(X^A)\right|,\ \Delta_{j'} f(X)\left|\Delta_{j'} f(X^{A'})\right|\right)\\
&\hskip20pt\leq c\left(\left(\E|(\Delta_j f(X)-\tilde\Delta_j f(X))|\right)^{\frac{1}{p}}
+\left(\E|(\Delta_{j'} f(X)-\tilde\Delta_{j'} f(X))|\right)^{\frac{1}{p}}\right).\nonumber
\end{align*}
The rest of the analysis can be carried out in the exact same way to get bounds on the Kolmogorov distance.
This completes the proof of Theorem \ref{thm:poissonmst} .

\section{Percolation estimates in the lattice setup}\label{sec:proofs-lattice}
We will now give an analogue of Lemma \ref{lem:->connectivity bound}.
\begin{lemma}\label{lem:->connectivity bound lattice}
Assume that $d\geq 2$, $p_1\in(0,p_c(\Z^d))$, $p_2\in(p_c(\Z^d),1)$ and $n\geq 1$. Then we have the following estimates:
\begin{equation}\label{eqn:->connectivity bound d=3 lattice}
\PR\big(\{0,e_1\}\underset{p}{\overset{2}{\leftrightsquigarrow}}B(n)\big)\leq\left\{
\begin{array}{l}
c_{19}\exp(-c_{20}n) \text{ if }p\leq p_1,\\
c_9\left(\log n/n\right)^{1/2}\text{ if }p\in[p_1,p_2],\\
c_{21}\exp(-c_{22} n)\text{ if }p\geq p_2.
\end{array}
\right.
\end{equation}
The constants appearing here depend only on $p_1$, $p_2$ and $d$. The same bounds hold for $\PR(B(1)\underset{p}{\overset{2}{\longleftrightarrow}}B(n))$. Further,
\begin{equation}\label{eqn:->connectivity bound d=3 lattice-in-cube}
\PR\big(B(1)\underset{p}{\overset{2}{\leftrightsquigarrow}}B(n)\text{ in }Q\big)
\leq\left\{
\begin{array}{l}
c_{19}\exp(-c_{20}n) \text{ if }p\leq p_1,\\
c_{21}\exp(-c_{22} n)\text{ if }p\geq p_2,
\end{array}
\right.
\end{equation}
whenever $Q$ is a cube containing the origin and $\partial^{\inn} B(n)$ has a vertex in $Q$.
\end{lemma}
\noindent{\bf Proof:}
The bounds in the subcritical regime follow from Menshikov's Theorem (see e.g. \cite{grimmett}).
When $d\geq 3$ and $p\geq p_2$, exponential decay will follow from the proof
of \cite[Lemma 7.89]{grimmett}. When $d=2$ and $p\geq p_2$, the stated bound follows from arguments similar to the ones used in the proof of Proposition 10.13 in \cite{penrosebook}.
The bound for $p\in[p_1,p_2]$ is just the content of Lemma \ref{lem:cerf}.
\hfill$\blacksquare$
\section{Rate of convergence in the CLT in the lattice setup}
We will prove Theorem \ref{thm:latticemst} in this section.
Let $u_1,\hdots,u_{\ell}$ be the edges of $\Z^d$ having both endpoints in $B(n)$, and let
$X_1,\hdots,X_{\ell}$ be the weights associated with them. Define $X=(X_1,\hdots,X_{\ell})$
and let $X'=(X_1',\hdots,X_{\ell}')$ be an independent copy of $X$. Write $F_{\mu}$ for the distribution function of $X_1$.
Fix $\alpha\in(0,1)$. We will make an appropriate choice of $\alpha$ later. Let
\begin{align}\label{eqn:def-L}
\cJ:=\big\{j\mid \text{ both endpoints of }u_j\text{ are in }B(n-n^{\alpha})\big\},\ \text{ and }\ \cL:=\{1,\hdots,\ell\}.
\end{align}
For each $j\in\cL$, choose and fix an endpoint $x_j$ of $u_j$, and let
\begin{align}\label{eqn:def-B-j-tilde-lattice}
B_j:=B(x_j, 1)\cap B(n),\ \text{ and }\ \tilde B_j:=B(x_j,n^{\alpha})\cap B(n).
\end{align}
Thus $\tilde B_j=B(x_j, n^{\alpha})$ if $j\in\cJ$.

We will apply \eqref{eqn:chatterjee} with
\[f(X)=M(B(n), X).\]
As in the proof of Theorem \ref{thm:poissonmst}, we will use the shorthand
\[
\Delta_j f(X^A):=\Delta_j f(X^A, X')
\]
for any $A\subset\cL$ and $j\in\cL$. We further define
\begin{align}\label{eqn:def-tilde-delta-lattice}
\tilde\Delta_j f(X^A):= M\big(\tilde B_j, X^A\big)-M\big(\tilde B_j, X^{A\cup\{j\}}\big)
\end{align}
for every $j\in\cL$ and $A\subset\cL$.

\subsection{Preliminary estimates}\label{sec:preliminary-estimate-lattice}
In this section we give an analogue of Proposition \ref{prop:preliminary-estimate}.
\begin{proposition}\label{prop:preliminary-estimate-lattice}
The following hold.
\begin{enumeratei}
\item Let $Z_j$ be the maximum of the weights associated with the edges of $B_j-u_j$. Let $\PR_1$ denote probability
conditional on the weights associated with the edges of $B_j-u_j$. Then for $j\in\cL$,
\begin{align}\label{eqn:47}
\E|\Delta_j f(X)-\tilde{\Delta}_j f(X)|
\leq 2\E\bigg[Z_j\cdot\int_0^1\PR_1\bigg(B_j\underset{F_{\mu}(uZ_j)}{\overset{2}{\leftrightsquigarrow}}\tilde B_j\text{ in }B(n)\bigg)du\bigg].
\end{align}
\item Order of variance:
\begin{align}\label{eqn:48}
\var\big(M\big(B(n), X\big)\big)=\Theta(n^d).
\end{align}
\end{enumeratei}
\end{proposition}

\noindent{\bf Proof:}
For $j\in\cL$, define $Y_j$ to be the maximum edge-weight in the path connecting the two endpoints
of $u_j$ in an MST of $B(n)-u_j$, when
the edge-weights are given by the appropriate subvector of $X$. From the add and delete algorithm (Section \ref{sec:add-and-delete-algo}),
it follows that
\begin{align}\label{eqn:41}
M(B(n), X)=M(B(n)-u_j, X)+X_j-\max(X_j,Y_j)\quad \text{for }\ j\in\cL,
\end{align}
and a similar assertion is true when $X$ is replaced by $X^j$.
Similarly define $\tilde{Y}_j$ to be the maximum edge-weight in the path connecting the two endpoints
of $u_j$ in an MST of $\tilde B_j-u_j$. Then \eqref{eqn:41} holds if we replace $B(n)$ by $\tilde B_j$ and $Y_j$ by $\tilde Y_j$.

Note also that for $j\in\cL$
\begin{align}\label{eqn:42}
M\big(B(n)-u_j, X\big)=M\big(B(n)-u_j, X^j\big),\ \text{ and }\
M\big(\tilde B_j-u_j, X\big)=M\big(\tilde B_j-u_j, X^j\big).
\end{align}
Hence
\begin{align}\label{eqn:Z}
\bigg|\Delta_j f(X)-\tilde{\Delta}_j f(X)|&=|\max(X_j,Y_j)-\max(X_j,\tilde{Y}_j)
-\max(X_j',Y_j)+\max(X_j',\tilde{Y}_j)\bigg|\\
&\leq 2|Y_j-\tilde{Y}_j|.\nonumber
\end{align}

From Lemma \ref{lem:mst minimax criterion}, it follows that $Y_j\leq \tilde{Y}_j$. Combining this with the definition of $Z_j$,
we get $Y_j\leq \tilde{Y}_j\leq Z_j$. Thus
\begin{equation}\label{eqn:W}
\E|Y_j-\tilde{Y}_j|=\E\bigg(Z_j\cdot\E_1\frac{(\tilde{Y}_j-Y_j)}{Z_j}\bigg),
\end{equation}
where $\E_1$ denotes expectation conditional on the weights associated with the edges of
$B_j-u_j$. Then for a random variable $U$ following $\mathrm{Uniform}[0,1]$ distribution
that is independent of $X,X'$,
\begin{align}\label{eqn:U}
\E_1\frac{(\tilde{Y}_j-Y_j)}{Z_j}&=\PR_1\left(Y_j <U Z_j<\tilde{Y}_j\right)\\
&=\int_0^1\PR_1\left(Y_j <u Z_j<\tilde{Y}_j\right)du
\leq\int_0^1 \PR_1\bigg(B_j\underset{F_{\mu}(uZ_j)}{\overset{2}{\leftrightsquigarrow}}\tilde B_j\text{ in }B(n)\bigg)du,\nonumber
\end{align}
where the last inequality follows from an argument identical to the one given right after \eqref{eqn:introducing auxilaiary U}.
\eqref{eqn:47} follows upon combining \eqref{eqn:Z}, \eqref{eqn:W}, and \eqref{eqn:U}.

The conclusion in \eqref{eqn:48} is included in the more general Theorem \ref{thm:general graphs} whose proof will be given in Section \ref{sec:general graphs}. \hfill$\blacksquare$

\subsection{Proof of Theorem \ref{thm:latticemst}}\label{sec:actual proof of latticemst}
The proof can be divided into two parts.

\medskip

\noindent\textbf{Proof of \eqref{eqn:d=2}:}
Recall the definition of the sets $\cJ$ and $\cL$ from \eqref{eqn:def-L}.
Recall also that for every $j\in\cL$, we have chosen and fixed an endpoint $x_j$ of $u_j$. Mimicking the proof of Theorem \ref{thm:poissonmst}, we define the sets
\begin{align}
\mathfrak{E}^{(\alpha)} = \big\{&(j, j', A, A'):\ j,j'\in\cL,\
A,A'\subsetneq\cL;\ j\notin A,\ j'\notin A'\text{ and}\nonumber\\
&\text{  either }j\notin\cJ,\text{ or }j'\notin\cJ,\text{ or }\|x_j-x_{j'}\|_{\infty}\leq 2n^{\alpha}\big\},\nonumber
\end{align}
and
\[\mathfrak{F}^{(\alpha)} =
\big\{(j, j', A, A'):\ j,j'\in\cL,\
A,A'\subsetneq\cL;\ j\notin A,\ j'\notin A'\big\}\setminus \mathfrak{E}^{(\alpha)}.\]
From \eqref{eqn:41} and \eqref{eqn:42}, it is clear that under the assumption of finite $(4+\delta)$-th moment on $\mu$,
\begin{equation}\label{eqn:delta_j finite moments lattice}
\E|\Delta_j f(X)|^{(4+\delta)}\leq C,\text{ for every }j\leq \ell.
\end{equation}
Hence \eqref{eqn:cov are bounded} remains true in our present setup. In view of \eqref{eqn:48}, \eqref{eqn:second term upper bound}continues to hold as well.

If we split the sum appearing in (\ref{eqn:var expansion term bounded})
into two parts $\Sigma_1$ (the sum over $\mathfrak{E}^{(\alpha)}$) and $\Sigma_2$ (the sum
over $\mathfrak{F}^{(\alpha)}$), then (\ref{eqn:sigma_1 is small}) continues to hold, i.e.,
\begin{align}\label{eqn:sigma_1 is small_lattice}
&\sum\nolimits_1
\frac{\mathrm{Cov}\big(\Delta_j f(X)\Delta_j f(X^A),\Delta_{j'} f(X)\Delta_{j'} f(X^{A'})\big)}
{\dbinom{\ell}{|A|}(\ell-|A|)\dbinom{\ell}{|A'|}(\ell-|A'|)}
\leq c' n^{2d-1+\alpha}.
\end{align}
Further, \eqref{eqn:decomposition of covariance terms}
and \eqref{eqn:bound on covariance} apply whenever $(j,j',A,A')\in \mathfrak{F}^{(\alpha)}$.
Let $T_1$ and $T_2$ be as in \eqref{eqn:bound on covariance}.
If $\mu$ has bounded support, then
\begin{equation}\label{eqn:X}
T_1+T_2\leq c\E|\Delta_j f(X)-\tilde{\Delta}_j f(X)|,
\end{equation}
and if $\mu$ has unbounded support and finite $(4+\delta)$-th moment, then
\begin{align}\label{eqn:Y}
T_1 &\leq \left(\E|\Delta_j f(X)-\tilde{\Delta}_j f(X)|\right)^{1/q'}\\
&\phantom{m}\cdot\left[\E\left( |\Delta_j f(X)-\tilde{\Delta}_j f(X)|
\left(|\Delta_j f(X^A)\Delta_{j'} f(X) \Delta_{j'} f(X^{A'})|\right)^q \right)\right]^{1/q}\nonumber\\
&=C_{\ref{eqn:Y}}\left(\E|\Delta_j f(X)-\tilde{\Delta}_j f(X)|\right)^{1/q'},\nonumber
\end{align}
where $q=1+\delta/3$ and $q'=1+3/\delta$. That $C_{\ref{eqn:Y}}$ is finite is ensured by
\eqref{eqn:delta_j finite moments lattice}. An application of H\"{o}lder's inequality
will give a similar bound for $T_2$. Let
\begin{align}\label{eqn:def-overline-q}
\overline{q}=
\left\{
\begin{array}{l}
1, \text{ if }\mu\text{ satisfies Property }B,\\
1+3/\delta, \text{ if }\mu\text{ satisfies Property }A_{\delta}.
\end{array}
\right.
\end{align}
We have thus shown
\begin{align}\label{eqn:77}
\mathrm{Cov}\left(\Delta_j f(X)\Delta_j f(X^A),\ \Delta_{j'} f(X)\Delta_{j'} f(X^{A'})\right)
\leq c\left(\E|\Delta_j f(X)-\tilde{\Delta}_j f(X)|\right)^{1/\overline{q}},
\end{align}
whenever $(j,j',A,A')\in\mathfrak{F}^{(\alpha)}$.

We now consider two possibilities separately.

\medskip
\noindent{\bf If $\mu$ satisfies Property C.}
If there exists a unique $x\in\R$ such that $\mu[0,x]=p_c(\Z^d)$, then there are two possibilities. The first possibility is that the distribution function of $\mu$, namely $F_{\mu}$, is continuous at $x$, and the second possibility is $F_{\mu}(x-)<F_{\mu}(x)=p_c$.

Assume first that $F_{\mu}$ is continuous at $x$ where $x$ is the unique point such that $F_{\mu}(x)=p_c(\Z^d)$. Choose a small enough positive $\eps_0$ so that $F_{\mu}(x-\eps_0)>0$ and $F_{\mu}(x+\eps_0)<1$.
For $\epsilon>0$, define the functions
\[p_1(\epsilon)=F_{\mu}(x-\epsilon)\ \text{ and }\ p_2(\epsilon)=F_{\mu}(x+\epsilon).\]
Note that when $j\in\cJ$, the integral on the right side of \eqref{eqn:47} can be written as
$\int_0^1\PR_1\bigg(B_j\underset{F_{\mu}(uZ_j)}{\overset{2}{\leftrightsquigarrow}}\tilde B_j\bigg)du$.
For any $\eps\in(0,\eps_0)$, we can break up this integral into
\[
\int_0^{\min((x-\eps)/Z_j,1)},\ \  \int_{{\min((x-\eps)/Z_j,1)}}^{{\min((x+\eps)/Z_j,1)}} \ \ \text{and} \ \
\int_{{\min((x+\eps)/Z_j,1)}}^1,
\]
to get
\begin{equation}\label{eqn:T}
\int_0^1\PR_1\bigg(B_j\underset{F_{\mu}(uZ_j)}{\overset{2}{\leftrightsquigarrow}}\tilde B_j\bigg)du
\leq c_{23}\exp(-c_{24}n^{\alpha})+\frac{2\eps}{Z_j}\cdot c_9\left(\frac{\alpha\log n}{n^{\alpha}}\right)^{1/2}
\end{equation}
by an application of Lemma \ref{lem:->connectivity bound lattice}.
The constants $c_{23}$ and $c_{24}$ depend on $c_{19}$, $c_{20}$, $c_{21}$ and $c_{22}$
as in Lemma \ref{lem:->connectivity bound lattice} corresponding to the choices
$p_i=p_i(\eps),\ i=1,2$, and the constant $c_{9}$ is the one from
Lemma \ref{lem:->connectivity bound lattice}  corresponding to the choices
$p_i=p_i(\eps_0),\ i=1,2$.

From \eqref{eqn:47} and \eqref{eqn:T}, we get
\begin{align}\label{eqn:78}
\E|\Delta_j f(X)-\tilde{\Delta}_j f(X)|\leq
2c_{23}\E(Z_j)\exp(-c_{24}n^{\alpha})+ 4\eps c_9\left(\frac{\alpha\log n}{n^{\alpha}}\right)^{1/2}
\end{align}
for every $j\in\cJ$. Combining \eqref{eqn:77} with \eqref{eqn:78}, we get
\begin{align*}
&\sum\nolimits_2
\frac{\mathrm{Cov}(\Delta_j f(X)\Delta_j f(X^A),\Delta_{j'} f(X)\Delta_{j'} f(X^{A'}))}
{\dbinom{\ell}{|A|}(\ell-|A|)\dbinom{\ell}{|A'|}(\ell-|A'|)}\\
&\hskip100pt\leq c\cdot n^{2d}\bigg(2c_{23}\ \E(Z_j)\exp(-c_{24}n^{\alpha})+ 4\eps c_9\bigg(\frac{\alpha\log n}{n^{\alpha}}\bigg)^{1/2}\bigg)^{1/\overline{q}}.
\end{align*}
The last inequality combined with
\eqref{eqn:sigma_1 is small_lattice},
\eqref{eqn:48},
\eqref{eqn:chatterjee},
\eqref{eqn:var expansion term bounded}, and
\eqref{eqn:second term upper bound} (we have already observed that the last inequality holds in our present setup) yields
\[\cW(\nu_n,\gamma)\leq
c'\bigg[\frac{1}{n^{(1-\alpha)/2}}
+\bigg(2c_{23}\E(Z_j)\exp(-c_{24}n^{\alpha})+ 4\eps c_9\bigg(\frac{\alpha\log n}{n^{\alpha}}\bigg)^{1/2}\bigg)^{\frac{1}{2\overline{q}}}
+\frac{1}{n^{d/2}}\bigg],\]
where $c'$ is a constant free of $\eps$. We take $\alpha=2\bar q/(1+2\bar q)$ in the last inequality. It then follows that
\[\limsup_n \frac{n^{\frac{1}{2(1+2\overline q)}}}{(\log n)^{1/(4\overline q)}}\mathcal{W}(\nu_n,\gamma)
\leq c' \bigg(4 \eps c_9\bigg(\frac{2\bar q}{1+2\bar q}\bigg)^{1/2}\bigg)^{\frac{1}{2\overline{q}}}.\]
This inequality is true for any $\eps>0$, and recall that $c'$ and $c_9 (=c_9(p_1(\eps_0), p_2(\eps_0)))$ do not depend on $\eps$.
This shows that \eqref{eqn:d=2} holds in this case.

The argument is similar if ({\it i}) $\mu[0,x]=p_c(\Z^d)$ for some unique $x\in\R$ and $F_{\mu}(x-)<F_{\mu}(x)$ or ({\it ii}) $\mu[0,x)=p_c(\Z^d)$ for some unique $x\in\R$, so we do not repeat it.

\medskip
\noindent{\bf If $\mu$ does not satisfy Property C.}
Combining the bound in \eqref{eqn:47} with \eqref{eqn:77}, and
Lemma \ref{lem:->connectivity bound lattice}, we get
\begin{align}\label{eqn:V}
&\mathrm{Cov}\bigg(\Delta_j f(X)\Delta_j f(X^A),\ \Delta_{j'} f(X)\Delta_{j'} f(X^{A'})\bigg)\\
&\hskip50pt
\leq \sup_{0<p<1} c\bigg[\PR\bigg(B_j\underset{p}{\overset{2}{\leftrightsquigarrow}}\tilde B_j\bigg)\bigg]^{1/\bar q}
\leq c'\left(\frac{\log n}{n^{\alpha}}\right)^{1/2\bar q},\notag
\end{align}
whenever $(j,j',A,A')\in\mathfrak{F}^{(\alpha)}$, and hence
\begin{align*}
\sum\nolimits_2
\frac{\mathrm{Cov}(\Delta_j f(X)\Delta_j f(X^A),\Delta_{j'} f(X)\Delta_{j'} f(X^{A'}))}
{\dbinom{\ell}{|A|}(\ell-|A|)\dbinom{\ell}{|A'|}(\ell-|A'|)}
\leq c\cdot n^{2d}\left(\frac{\log n}{n^{\alpha}}\right)^{1/2\bar q}.
\end{align*}
Combining this inequality with \eqref{eqn:sigma_1 is small_lattice}, and taking $\alpha=2\bar q/(1+2\bar q)$, we get
\begin{equation}\label{eqn:first-term-bound-lattice}
\mathrm{Var}\big(\E(T| W)\big)\leq cn^{2d}\frac{(\log n)^{1/2\bar q}}{n^{1/(1+2\bar q)}}.
\end{equation}
The last bound together with
\eqref{eqn:48},
\eqref{eqn:chatterjee}, and
\eqref{eqn:second term upper bound}
yields the bound in \eqref{eqn:d=2}.

\medskip

\noindent\textbf{Proof of \eqref{eqn:L3}:}
We introduce
\begin{align}
\overline{\mathfrak{E}}^{(\alpha)} :=& \big\{(j, j', A, A'):
\ j,j'\in\cL,\
A,A'\subsetneq\cL;\ j\notin A,\ j'\notin A'\text{ and }\|x_j-x_{j'}\|_{\infty}\leq 2n^{\alpha}\big\},\text{ and}\nonumber
\end{align}
\[
\overline{\mathfrak{F}}^{(\alpha)} := \big\{(j, j', A, A'):
\ j,j'\in\cL,\
A,A'\subsetneq\cL;\ j\notin A,\ j'\notin A'\big\}\setminus \bar{\mathfrak{E}}^{(\alpha)}
\]
for $0<\alpha<1$.
We split the sum appearing in (\ref{eqn:var expansion term bounded})
into $\overline\Sigma_1$, the sum over $(j, j', A, A')\in\overline{\mathfrak{E}}^{(\alpha)}$
and $\overline\Sigma_2$, the sum over $(j, j', A, A')\in\overline{\mathfrak{F}}^{(\alpha)}$.
Then similar to \eqref{eqn:sigma_1 is small},
\begin{align}\label{eqn:sigma_1 is small lattice d>2}
\overline{\sum}_1\
&\frac{\mathrm{Cov}(\Delta_j M(X)\Delta_j M(X^A),\Delta_{j'} M(X)\Delta_{j'} M(X^{A'}))}
{\dbinom{\ell}{|A|}(\ell-|A|)\dbinom{\ell}{|A'|}(\ell-|A'|)}\\
&\hskip50pt\leq c |\{(j, j'): (j, j', \emptyset, \emptyset)\in\overline{\mathfrak{E}}^{(\alpha)}\}|
\leq  c n^{d+\alpha d}.\notag
\end{align}
Further, the argument leading to \eqref{eqn:V} yields
\begin{align}\label{eqn:V-1}
\mathrm{Cov}\bigg(\Delta_j f(X)\Delta_j f(X^A),\ \Delta_{j'} f(X)\Delta_{j'} f(X^{A'})\bigg)
\leq \sup_{0<u<1}
c\bigg[\PR_1\bigg(B_j\underset{F_{\mu}(uZ_j)}{\overset{2}{\leftrightsquigarrow}}\tilde B_j\text{ in }B(n)\bigg)\bigg]^{1/\bar q},
\end{align}
whenever $(j, j', A, A')\in\overline{\mathfrak{F}}^{(\alpha)}$, where $\bar q$ is as in \eqref{eqn:def-overline-q}.
Since $\mu$ satisfies Property $D$ by assumption, $\mathrm{Range}(F_{\mu})\subset(p_c-\eps,p_c+\eps)^c$ for some $\eps>0$.
It thus follows from \eqref{eqn:V-1} and Lemma \ref{lem:->connectivity bound lattice} that
\begin{align}
\mathrm{Cov}\bigg(\Delta_j f(X)\Delta_j f(X^A),\ \Delta_{j'} f(X)\Delta_{j'} f(X^{A'})\bigg)
&\leq \sup_{p\notin(p_c-\eps,p_c+\eps)}
c\bigg[\PR\bigg(B_j\underset{F_{\mu}(uZ_j)}{\overset{2}{\leftrightsquigarrow}}\tilde B_j\text{ in }B(n)\bigg)\bigg]^{1/\bar q}\\
&\leq c'\exp(-c''n^{\alpha})\nonumber
\end{align}
whenever $(j, j', A, A')\in\overline{\mathfrak{F}}^{(\alpha)}$. Hence
\begin{align}\label{eqn:numerator boundII lattice I}
\overline\sum_2\
\frac{\mathrm{Cov}(\Delta_j M(X)\Delta_j M(X^A),\Delta_{j'} M(X)\Delta_{j'} M(X^{A'}))}
{\dbinom{\ell}{|A|}(\ell-|A|)\dbinom{\ell}{|A'|}(\ell-|A'|)}
\leq c n^{2d}\exp(-c' n^{\alpha})
\end{align}
As before, we combine
\eqref{eqn:48},
\eqref{eqn:sigma_1 is small lattice d>2},
\eqref{eqn:numerator boundII lattice I},
\eqref{eqn:second term upper bound}, and
\eqref{eqn:chatterjee}
to conclude that
\[\mathcal{\cW}(\nu_n,\gamma)\leq c/n^{\frac{d(1-\alpha)}{2}}.\]
We get the bound in \eqref{eqn:L3} once we replace $d(1-\alpha)/2$ by $\eta$.

\medskip

\noindent{\bf Bounds on the Kolmogorov distance.}
Bounds on the Kolmogorov distance can be obtained by using Theorem \ref{thm:lrp} and following the same line of arguments (see the discussion at the end of Section \ref{sec:actual proof of poissonmst}). Note the presence of the term $\E|\Delta_j f(X, X')|^6$ in \eqref{eqn:lrp}.
We require $\mu$ to satisfy either Property $B$ or Property $A_{\delta}$ with $\delta\geq 2$ to show that $\E|\Delta_j f(X, X')|^6<\infty$. Rest of the argument goes through verbatim.

This concludes the proof of Theorem \ref{thm:latticemst}.

\section{General graphs: Proof of Theorem \ref{thm:general graphs}}\label{sec:general graphs}
To fix ideas, we first assume that $G$ is symmetric, i.e., for every two pairs of adjacent vertices $v_1, v_2$ and $v_1', v_2'$, there exists a graph automorphism $f$ of $G$ such that $f(v_i)=v_i'$, $i=1, 2$.

If $G$ is symmetric and deletion of an edge of $G$ creates two components,
then $G$ is a regular tree. Hence
all our claims follow trivially. So we can assume that
this is not the case.

Let $E_n=\{u_1,\hdots,u_{\ell_n}\}$ and let $X_1,\hdots,X_{\ell_n}$ be the associated edge weights.
Let $X=(X_1,\hdots,X_{\ell_n})$ and let $X'=(X_1',\hdots,X_{\ell_n}')$ be an independent copy of $X$.
Then $M_n=M(G_n, X)$. As before, we want to apply Theorem \ref{thm:chatterjee} with
\[f(X):=M(G_n, X).\]
As in the proof of Theorem \ref{thm:poissonmst}, we will use the shorthand
\[\Delta_j f(X^A):=\Delta_j f(X^A, X'),\]
for any $A\subset[\ell_n]$ and $1\leq j\leq\ell_n$.

Define
\[\mathcal{I}_r^n:=\big\{i\leq \ell_n:\ S(v,r)\subset G_n\text{ for each endpoint }v\text{ of }u_i\big\}.\]
For large $r$ and $i\in\mathcal{I}_r^n$, fix an endpoint $v_i$ of $u_i$ and let $Y_i^n$ (resp. $Y_i^n(r)$)
be the maximum edge weight in a path connecting the endpoints of  $u_i$ in an MST of
$G_n-u_i$ (resp. $S(v_i,r)-u_i$) with the edge weights being the appropriate subvector
of $X$. We will suppress the dependence on $n$ and
simply write $\mathcal{I}_r$, $Y_i$ and $Y_i(r)$.

An application of Lemma \ref{lem:conditional var bound} yields, with our usual notation,
\begin{align}\label{eqn:general graph variance lower bound I}
\mathrm{Var}\big(f(X)\big)&\geq \sum_{i=1}^{\ell_n}\var\big(\E(f(X)|X_i)\big)
=\frac{1}{2}\sum_{i=1}^{\ell_n}\E\bigg[\E\big(f(X)|X_i\big)-\E\big(f(X^{i}\big)|X_i')\bigg]^2\\
&\geq \frac{1}{2}\sum_{i\in\mathcal{I}_r}\E\bigg[\E\big(f(X)|X_i\big)-\E\big(f(X^{i})|X_i'\big)\bigg]^2
=\frac{1}{2}\sum_{i\in\mathcal{I}_r}\E\bigg[\E\big(f(X)-f(X^{i})\mid X_i, X_i'\big)\bigg]^2.\nonumber
\end{align}
By the add and delete algorithm (Section \ref{sec:add-and-delete-algo}),
for $i\in\mathcal{I}_r$,
\[f(X)=M(G_n-u_i, X)+X_i-\max(X_i,Y_i),\]
and hence
\begin{align}\label{eqn:900}
f(X)-f(X^{i})=\min(X_i,Y_i)-\min(X_i',Y_i).
\end{align}
Since $\mu$ is non-degenerate, we can find real numbers $b>a$ such that
$\mu[0,a]>0$ and $\mu[b,\infty]>0$. Going back to (\ref{eqn:general graph variance lower bound I}),
\begin{align}\label{eqn:general graph variance lower bound II}
\var\big(f(X)\big)&\geq \frac{1}{2}\sum_{i\in\mathcal{I}_r}
\E\bigg[\E\big(\min(X_i,Y_i)-\min(X_i',Y_i)\mid X_i, X_i'\big)\bigg]^2\\
&\hskip20pt\geq \frac{1}{2}\sum_{i\in\mathcal{I}_r}\E\bigg[\big((b-a)\PR(Y_i\geq b)\big)^2\I\{X_i\leq a, X_i'\geq b\}\bigg]\nonumber\\
&\hskip40pt\geq \frac{1}{2}|\mathcal{I}_r|(b-a)^2\cdot p^2\cdot\mu[0,a]\cdot\mu[b,\infty),\nonumber
\end{align}
where
\begin{align*}
p:=\PR\big(\text{The weight associated with each edge sharing one vertex with }u_i\text{ is at least }b\big).
\end{align*}
Note that $p$ does not depend on the edge $u_i$ since $G$ is symmetric. By assumption (III),
\begin{align}\label{eqn:general graph variance lower bound 5}
|\mathcal{I}_r|=\Theta(|V_n|).
\end{align}
From
\eqref{eqn:general graph variance lower bound II} and
\eqref{eqn:general graph variance lower bound 5}, it follows that
\[\mathrm{Var}\big(f(X)\big)\geq c|V_n|.\]
The upper bound is a simple consequence of the Efron-Stein inequality:
\[\var\big(f(X)\big)\leq \frac{1}{2}\sum_{j=1}^{\ell_n}\E\big(\Delta_j f(X)\big)^2.\]
Thus we have proven that $\var(M_n)=\var(f(X))=\Theta(|V_n|)$.

Turning toward the proof of the central limit theorem, define, for large $r$,
\begin{align}
\mathfrak{E}_n(r)&:=\big\{(j,j',A,A'):\ j,j'\leq\ell_n,\ A,A'\subsetneq\{1,\hdots,\ell_n\},\ j\notin A,\  j'\notin A'\nonumber\\
&\text{ and either } d_G(x_j,x_{j'})\leq 2r\text{ or } S(x_j,r)\not\subset G_n
\text{ or } S(x_{j'},r)\not\subset G_n\nonumber\\
&\text{ for some endpoints }x_j, x_{j'}\text{ of }u_j\text{ and }u_{j'}\text{ respectively}\big\},\text{ and}\nonumber\\
\mathfrak{F}_n(r)&=\big\{(j,j',A,A'):\ j,j'\leq\ell_n,\ A,A'\subsetneq\{1,\hdots,\ell_n\}, j\notin A,\  j'\notin A'\big\}\setminus
\mathfrak{E}_n(r).\nonumber
\end{align}
Proceeding as before, we split the sum in (\ref{eqn:var expansion term bounded}) into
$\Sigma_1$, the sum over all $(j,j',A,A')\in\mathfrak{E}_n(r)$ and $\Sigma_2$,
the sum over the rest of the terms. It follows from \eqref{eqn:900} that
$|\Delta_j f(X)|\leq |X_j-X_j'|$. Further, $\E(X_j^4)<\infty$. Thus, a computation similar to (\ref{eqn:sigma_1 is small})
will yield
\begin{align}\label{eqn:sigma_1 is small general graph}
&\sum\nolimits_1
\frac{\mathrm{Cov}\big(\Delta_j f(X)\Delta_j f(X^A),\ \Delta_{j'} f(X)\Delta_{j'} f(X^{A'})\big)}
{\dbinom{\ell_n}{|A|}(\ell_n-|A|)\dbinom{\ell_n}{|A'|}(\ell_n-|A'|)}\\
&\hskip70pt\leq c|V_n|\big(a_r+|\{v\in V_n:\ S(v,r)\not\subset G_n\}|\big),\nonumber
\end{align}
where $a_r:=|\{v'\in V:\ d_G(v,v')\leq 2r\}|$ for some (and hence all, by symmetry) $v\in V$.

For $j\in\mathcal{I}_r$, define
\[\tilde{\Delta}_j f(X)=M\big(S(v_j,r), X\big)-M\big(S(v_j,r), X^j\big).\]
 With this definition of $\tilde{\Delta}_j f(X)$,
\eqref{eqn:decomposition of covariance terms} and \eqref{eqn:bound on covariance} hold
for $(j,j',A,A')\in\mathfrak{F}_n(r)$.
As in \eqref{eqn:X} and \eqref{eqn:Y}, we get
\[T_1\leq c\left(\E|\Delta_j f(X)-\tilde{\Delta}_j f(X)|\right)^{\frac{1}{1+3/\delta}}\]
for some $\delta\geq 0$, where $\delta=0$ if $\mu$ satisfies Property $B$ and $\delta>0$ if $\mu$ satisfies Property $A_{\delta}$.
A similar bound holds for $T_2$. A calculation similar to (\ref{eqn:Z}) yields
\begin{equation}\label{eqn:|delta-tilde delta| for general graphs}
|\Delta_j f(X)-\tilde{\Delta}_j f(X)|\leq 2(Y_j(r)-Y_j)
\end{equation}
for $j\in\mathcal{I}_r$.

Fix a vertex $v$ of $G$ and let $e$ be an edge incident to $v$.
Let $Y(v,e,r)$ be the maximum edge weight in the path connecting the endpoints
of $e$ in an MST of $S(v,r)-e$,
clearly $Y(v,e,r)$ is decreasing in $r$. Define
\[Y(v,e):=\lim_{r\to\infty}Y(v,e,r).\]
The above convergence also holds in $L^1$ as a consequence of dominated convergence theorem.
Since $G$ is symmetric, $Y_i^n(r)$ has the same distribution as $Y(v,e,r)$
and $Y_i^n$ dominates $Y(v,e)$ stochastically for every $i\in\mathcal{I}_r^n$.
Hence
\begin{align}
&\lim_{r\to\infty}\limsup_{n\to\infty}\bigg[\max_{i\in\mathcal{I}_r^n}\ \E\big(Y_i^n(r)-Y_i^n\big)\bigg]
\leq \lim_{r\to\infty}\E\big(Y(v,e,r)-Y(v,e)\big)=0.
\end{align}
Thus we have
\begin{align}\label{eqn:general graph variance lower bound 3}
\lim_{r\to\infty}\limsup_{n\to\infty}\max_{\substack{(j,j',A,A')\\ \in\mathfrak{F}_n(r)}}
\mathrm{Cov}\left(\Delta_j f(X)\Delta_j f(X^A),\Delta_{j'} f(X)\Delta_{j'} f(X^{A'})\right)=0,
\end{align}
which gives us control over $\sum_2$.
Further,
\[\frac{1}{\mathrm{Var}\big(f(X)\big)^{3/2}}\sum_{j=1}^{\ell_n}\E|\Delta_j f(X)|^3\leq \frac{c}{|V_n|^{1/2}}.\]
The last inequality together with (\ref{eqn:sigma_1 is small general graph}),
(\ref{eqn:general graph variance lower bound 3}),
(\ref{eqn:chatterjee}) and the fact that
$\mathrm{Var}(M_n)=\Theta(|V_n|)$ yields
\[\limsup_n \mathcal{W}(\mu_n,\gamma)=0,\]
where $\mu_n$ is the law of $(M_n-\E M_n)/\sqrt{\mathrm{Var}(M_n)}$.

Assume now that $G$ is vertex-transitive, so that there are two kinds of edges.
Call an edge $e\in E$ of type A if deletion of $e$ results in the creation
of two disjoint components. We say $e$ is of type B if it is not of type A.
Define $\tilde{\mathcal{I}}_r^n:=\{i\in \mathcal{I}_r:\ i\text{ is of type B}\}$.
Define $Y_i^n$ and $Y_i^n(r)$ as before for each $i\in\tilde{\mathcal{I}}_r^n$.
Then as in \eqref{eqn:general graph variance lower bound I},
\begin{align*}
\mathrm{Var}\big(f(X)\big)&\geq
\frac{1}{2}\sum_{i\in\tilde{\mathcal{I}}_r}\E\bigg[\E\big(f(X)|X_i\big)-\E\big(f(X^{i})|X_i'\big)\bigg]^2.
\end{align*}
Note also that $|\tilde{\mathcal{I}}_r|=\Theta(|V_n|)$ if $G$ is not a tree. So we can argue
as before to conclude that $\mathrm{Var}(f(X))=\Theta(|V_n|)$.

Next, note that if $j\in\mathcal{I}_r^n$ and $u_j$ is of type A,
then $\Delta_j f(X)-\tilde{\Delta}_j f(X)=0$. Further, our previous arguments show that
\begin{align}
&\lim_{r\to\infty}\limsup_{n\to\infty}\bigg[\max_{i\in\tilde{\mathcal{I}}_r^n}\ \E\big(Y_i^n(r)-Y_i^n\big)\bigg]
\leq \lim_{r\to\infty}\sum\nolimits_{\ast}\E\big(Y(v,e,r)-Y(v,e)\big)=0,
\end{align}
where $\sum_{\ast}$ is the sum over all type B edges $e$ incident to $v$.
The rest of the arguments remain the same. This finishes the proof of the central limit theorem.

\appendix
\section{} \label{Appendix}
\subsection{Completing the proof of Lemma \ref{lem:->connectivity bound}}
The following proposition fills in the gap in the proof of Lemma \ref{lem:->connectivity bound}.
\begin{proposition*} Assume that $n\geq 2$, $a\in [1/2,\log n]$, and $r_c\leq r\leq (\log n)^2$. Then there exists positive universal constants $c_{12}$ and $\beta$ such that
\[\PR\big(B_{\bR^2}(a)\underset{r}{\overset{2}{\longrightarrow}}B_{\bR^2}(n)\big)\leq c_{12}/n^{\beta}.\]
\end{proposition*}
\noindent{\bf Proof:}
As usual $\poi$ will denote a Poisson process of intensity one.
Let $\sigma((a,b);r,j)$ denote the probability of an occupied crossing of the rectangle $[0,a]\times[0,b]$
at level $r$ in the $j$-th direction, $j=1,2$; that is
\begin{align*}
\sigma((a,b);r,1)=&\PR(\poi^{(r)}\text{ contains a curve }\gamma\subset[0,a]\times[0,b]\\
&\phantom{\PR m}\text{such that }\gamma\text{ intersects both }S_1\text{ and }S_2)
\end{align*}
where $S_1=\{0\}\times[0,b]$ and $S_2=\{a\}\times[0,b]$
and define $\sigma((a,b);r,2)$ similarly.
First we note that
\[\sigma((m,3m);r_c,1)\geq \kappa_0:=(9e)^{-122}\text{ whenever }m>r_c.\]
(We can prove this assertion by observing that $\sigma((m,3m);r,1)$ is a continuous function of $r$ and then using arguments similar to the ones given right after \eqref{eqn:vacant-crossing} and \cite[Lemma 3.3]{royII}.)

Now, the proof of Lemma 4.4 of \cite{royII} applies to occupied crossings as well. Since
$\sigma((m,3m);r_c,1)\geq \kappa_0$ for $m>r_c$, the arguments of Lemma 4.4 of \cite{royII}
would furnish positive constants $f(t)$ for each $t>0$ such that
$$\sigma((m,(1+t)m);r_c,1)\geq f(t).$$

Applying Theorem 2.1 of \cite{alexander} with the parameters
$h=\ell/(1+t)$ and $b=\ell/(1+t)^2$ with $t$ small enough so that $2/(1+t)^2-1/2>1+\eps$
(for some positive $\eps$) and $(1+t)^2<4/3$ and $\ell$ large so that $h>4r_c$ and $b>\ell/2+2r_c$,
we get
\begin{align}
&\sigma\bigg(\bigg(\ell\left[\frac{2}{(1+t)^2}-\frac{1}{2}\right]+r_c,\frac{\ell}{1+t}-2r_c\bigg); r_c,1\bigg)\nonumber\\
&\geq c
\sigma\bigg(\bigg(\frac{\ell}{(1+t)^2}+r_c,\frac{\ell}{1+t}-4r_c\bigg); r_c,1\bigg)^4\times
\sigma\bigg(\bigg(\ell,\frac{\ell}{1+t}+3r_c\bigg); r_c,2\bigg)^2\nonumber
\end{align}
for large $\ell$. Hence
\begin{align}
&\sigma\bigg(\bigg(\ell(1+\eps),\ell\bigg); r_c,1\bigg)\nonumber\\
&\geq c
\sigma\bigg(\bigg(\frac{\ell}{(1+3t/4)^2}+r_c,\frac{\ell}{1+5t/4}\bigg); r_c,1\bigg)^4\times
\sigma\bigg(\bigg(\ell,\frac{\ell}{1+t/2}\bigg); r_c,2\bigg)^2\nonumber\\
&\geq cf\bigg(\frac{(1+3t/4)^2}{1+5t/4}-1\bigg)^4\times f\left(t/2\right)^2\nonumber
\end{align}
for every $\ell$ bigger than a fixed threshold $\ell_0$.
Hence Lemma 3.1 of \cite{alexander} yields
\begin{equation}\label{eqn:RSW}
\sigma\left((3\ell,\ell); r_c,1\right)\geq \kappa_1
\end{equation}
for a positive constant $\kappa_1$ and $\ell\geq\ell_0$.

Let $A_k$ be the event that there is an occupied circuit at level $r_c$
in the annulus $B_{\bR^2}(3\ell_k/2)\setminus B_{\bR^2}(\ell_k/2)$, where $\ell_k=3\ell_{k-1}+4r_c$
and $\ell_1=\max(2a+2r,\ell_0)$. FKG inequality and (\ref{eqn:RSW})
gives $\PR(A_k)\geq \kappa_1^4$. Hence
\begin{align}
\PR\big(B_{\bR^2}(a)\underset{r}{\overset{2}{\longrightarrow}}B_{\bR^2}(n)\big)&\leq \PR(A_1^c\cap\hdots\cap A_t^c)
=\prod_{k=1}^t\PR(A_k^c)\leq (1-\kappa_1^4)^t\nonumber
\end{align}
where $3\ell_t/2+r_c\leq n-r<3\ell_{t+1}/2+r_c$. This yields the desired bound.\hfill$\blacksquare$



\subsection{ Proof of Lemma \ref{lem:critical connectivityl lattice}}
Fix $p\in[p_1,p_2]$.
Let $u_1,\hdots,u_m$ be the edges of $\Z^d$ both of whose endpoints lie
in $B(n)$ and let $X_1,\hdots,X_m$ be i.i.d. Bernoulli$(p)$ random variables
(i.e. $\PR(X_1=1)=p=1-\PR(X_1=0)$) associated to them.
Let $X:=(X_1,\hdots,X_m)$ and let $X':=(X_1',\hdots,X_m')$ be an independent copy
of $X$. As earlier we define the event
\begin{align}
E:=\big\{\text{there is exactly one }p\text{-cluster in }B_{\bZ^d}(n)
\text{ that intersects both }B_{\bZ^d}(a_n)\text{ and }\partial^{\inn} B_{\bZ^d}(n)\big\}\nonumber
\end{align}
for some $a_n\to\infty$ in a way so that $a_n=o(n)$. Define the function $f$ by
$f(X):=\I_E(X).$ Then an application of Lemma \ref{lem:conditional var bound} yields
\begin{equation}\label{eqn:conditional var bound application lattice}
\mathrm{Var}(f(X))\geq \sum_{i\in\mathcal{I}} \mathrm{Var}\left[\E(f(X)|X_i)\right]
\end{equation}
where $\mathcal{I}:=\{i\leq m:\ \text{both endpoints of }u_i\text{ lie in }B(a_n/3)\}$.
Fix $i\in\mathcal{I}$, denote the endpoints of $u_i$ by $v_1$ and $v_2$.
With our usual notation,
\begin{align}\label{eqn:lower bound on variance I lattice}
\mathrm{Var}[\E\left(f(X)|X_i\right)] &=
\frac{1}{2}\E\left[\left(\E\bigl(f(X)|X_i\bigr)-\E\bigl(f(X^{ i})|X_i'\bigr)\right)^2\right]\\
&\geq \frac{1}{2}\E\left[\PR(A_i)^2 \I\{X_i=1, X_i'=0\}\right],\nonumber
\end{align}
where
\begin{align}
A_i=&\{u_i\underset{p}{\overset{2}{\leftrightsquigarrow}}B_{\bZ^d}(n)-u_i,\text{ any }p\text{-cluster in }B_{\bZ^d}(n)-u_i
\text{ that intersects}\nonumber\\
&\text{ both }\partial^{\inn} B_{\bZ^d}(n)\text{ and }B_{\bZ^d}(a_n)\text{ contains either }v_1\text{ or }v_2\}.\nonumber
\end{align}
Now,
\begin{align}
\PR(A_i)&\geq\PR\big(u_i\underset{p}{\overset{2}{\leftrightsquigarrow}}B_{\bZ^d}(v_1, 2n)-u_i,
\text{ if }\C\text{ is a }p\text{-cluster in }B_{\bZ^d}(v_1,2n)-u_i\nonumber\\
&\phantom{mm}\text{then every connected component of }\C\cap B_{\bZ^d}(v_1,n/2)\nonumber\\
&\phantom{mm}\text{that intersects }\text{ both }\partial^{\inn} B_{\bZ^d}(v_1,n/2)\text{ and }B_{\bZ^d}(v_1,2a_n)\nonumber\\
&\phantom{mm}\text{contains either }v_1\text{ or }v_2\big)\nonumber\\
&= \PR(F)/(1-p),\nonumber
\end{align}
where
\begin{align}
F:=&\big\{\{0,e_1\}\underset{p}{\overset{2}{\leftrightsquigarrow}}B_{\bZ^d}(2n),
\text{ if }\C\text{ is a }p\text{-cluster in }B_{\bZ^d}(2n)-\{0,e_1\}\nonumber\\
&\ \text{ then every connected component of }\C\cap B_{\bZ^d}(n/2)\nonumber\\
&\ \text{ that intersects }\text{ both }\partial^{\inn} B_{\bZ^d}(n/2)\text{ and }B_{\bZ^d}(2a_n)\text{ contains either }0\text{ or }e_1\big\}.\nonumber
\end{align}
From \eqref{eqn:conditional var bound application lattice} and \eqref{eqn:lower bound on variance I lattice},
we conclude that
\begin{equation}\label{eqn:1}
\PR(F)\leq c/a_n^{d/2}.
\end{equation}
Note that
\begin{equation}\label{eqn:decomposition into two terms lattice}
\PR\big(\{0,e_1\}\underset{p}{\overset{2}{\leftrightsquigarrow}}B_{\bZ^d}(2n)\big)
\leq \PR(F)+\PR\big(B_{\bZ^d}(2a_n)\underset{p}{\overset{3}{\leftrightsquigarrow}}B_{\bZ^d}(n/2)\big).
\end{equation}
We now define a cube $Q\subset B_{\bZ^d}(n/2)$ to be a trifurcation box in $B_{\bZ^d}(n/2)$
at level $p$, if
\begin{enumeratei}
\item there is a $p$-cluster $\C$ in $B_{\bZ^d}(n/2)$ with $\C\cap Q\neq\emptyset$, and
\item the vertices of $\C$ contained in $B_{\bZ^d}(n/2)-Q$ contain at least three
$p$-clusters in $B_{\bZ^d}(n/2)-Q$ each of which intersects $\partial^{\inn} B_{\bZ^d}(n/2)$.
\end{enumeratei}
We can then apply the arguments in the proof of Lemma \ref{lem:burton-keane} (see the arguments leading up to \eqref{eqn:burton-keane IV}) to show that
\[\PR\big(B_{\bZ^d}(2a_n)\text{ is a trifurcation box in }B_{\bZ^d}(n/2)\text{ at level }p\big)\leq \frac{ca_n^d}{n},\]
from which it will follow that
\begin{equation}\label{eqn:number}
\PR\big(B_{\bZ^d}(2a_n)\underset{p}{\overset{3}{\leftrightsquigarrow}}B_{\bZ^d}(n/2)\big)\leq c\exp(c' a_n)\frac{a_n^d}{n}.
\end{equation}
Combining \eqref{eqn:1}, \eqref{eqn:decomposition into two terms lattice} and \eqref{eqn:number},
we choose $c'a_n=\log n/2$ to get the desired bound.\hfill$\blacksquare$

\vskip.3in
\noindent{\bf Acknowledgments.} The authors are indebted to  Larry Goldstein and \"{U}mit I\c{s}lak for their valuable help in improving the manuscript and checking the proof.
The authors also thank two anonymous referees whose careful reading and detailed reports improved the presentation considerably.
The work of SC was partially supported by NSF grant DMS-1005312.
The work of SS was supported in part by  NSF grant DMS-1007524 and Netherlands Organisation for Scientific Research (NWO) through the Gravitation Networks grant 024.002.003.

\end{document}